# ALGEBRAIC STRUCTURES USING NATURAL CLASS OF INTERVALS

W. B. Vasantha Kandasamy
Florentin Smarandache

**2011**

# ALGEBRAIC STRUCTURES USING NATURAL CLASS OF INTERVALS

**W. B. Vasantha Kandasamy**
**Florentin Smarandache**

**2011**

# CONTENTS

# PREFACE

Authors in this book introduce a new class of intervals called the natural class of intervals, also known as the special class of intervals or as natural intervals. These intervals are built using increasing intervals, decreasing intervals and degenerate intervals. We say an interval [a, b] is an increasing interval if a < b for any a, b in the field of reals R.

An interval [a, b] is a decreasing interval if a > b and the interval [a, b] is a degenerate interval if a = b for a, b in the field of reals R. The natural class of intervals consists of the collection of increasing intervals, decreasing intervals and the degenerate intervals.

Clearly R is contained in the natural class of intervals. If R is replaced by the set of modulo integers $Z_n$, n finite then we take the natural class of intervals as [a, b] where a, b are in $Z_n$ and we do not say a < b or a > b for such ordering does not exist on $Z_n$.

The authors extend all the arithmetic operations without any modifications on the natural class of intervals. The natural class of intervals is closed under the operations addition, multiplication, subtraction and division.

In this book we build algebraic structures using this natural intervals. This book has eleven chapters. Chapter one describes all types of natural class of intervals and the arithmetic operations on them. Chapter two introduces the semigroup of natural class of intervals using R or $Z_n$ and study the properties associated with them.

Chapter three studies the notion of rings constructed using the natural class of intervals. Matrix theory using the special class of intervals is analyzed in chapter four of this book.

Chapter five deals with polynomials using interval coefficients. New types of rings of natural intervals are introduced and studied in chapter six. The notion of vector space using natural class of intervals is built in chapter seven. In chapter eight fuzzy natural class of intervals are introduced and algebraic structures on them is built and described. Algebraic structures using natural class of neutrosophic intervals are developed in chapter nine. Chapter ten suggests some possible applications. The final chapter proposes over 200 problems of which some are at research level and some difficult and others are simple.

One of the features of this book is it gives over 330 examples to make the book interesting. We thank Dr. K.Kandasamy for proof reading.

<div style="text-align: right;">
W.B.VASANTHA KANDASAMY<br>
FLORENTIN SMARANDACHE
</div>

**Chapter One**

# INTRODUCTION TO INCREASING AND DECREASING INTERVALS

In this chapter we for the first time introduce the notion of decreasing and increasing intervals and discuss the properties enjoyed by them. These notions given in this chapter will be used in the following chapters of this book.

*Notation:* Z the set of positive and negative, integers with zero. Q the set of rationals, R the set of reals and C the complex numbers.

We first define the notion of increasing intervals.

**DEFINITION 1.1**: *Let $[x, y]$ be an interval from Q or R or Z which is a closed interval and if $x < y$ (that is x is strictly less than y) then we define $[x, y]$ to be a closed increasing interval or increasing closed interval. That is both x and y are included.*
    *On the other hand $(x, y)$ if $x < y$ is the increasing open interval or open increasing interval. $[x, y)$ $x < y$ is the half*

*closed and half open increasing interval and (x, y] x < y is the half open and half closed increasing interval.*

We will first illustrate this situation by an example.

*Example 1.1:* [3, 8], [-8, 0], [-5, -2] and [0, 12] are closed increasing intervals. (9, 12), (-7, 4), (-5, 0) and (0, 2) are open increasing intervals or increasing open intervals. [3, 7), (0, 5], (-7, -2], (-3, 0] and [-3, 8) are half open-closed (closed-open) increasing intervals.

Now we proceed onto describe decreasing intervals.

**DEFINITION 1.2:** *Let [x, y] be an interval where x, y belongs to Z or Q or R with x > y (x is strictly greater than y) then we define [x, y] to be a decreasing closed interval or closed decreasing interval.*

*If (x, y) is considered with x > y then we say (x, y) is a decreasing open interval. Similarly [x, y) and (x, y], x > y are decreasing half open intervals.*

We will illustrate this situation by an example.

*Example 1.2:* (3, 0), (-7, -15), (0, -2) and (15, 9) are decreasing open intervals. [7, 0], [5, 2], [8, 0], [-1, -8] and [0, -4] are examples of decreasing closed intervals. (5, 3], (12, 0], [0, -7), [5, 2) and (-1, -7) are examples of half open or half closed decreasing intervals.

We define an interval to be a degenerate one if that interval is not an increasing one or a decreasing one. That is all intervals [x, y], (or (x, y) or [x, y) or (x, y)) in which x = y. Thus by using the concept of degenerate intervals we accept the totality of all numbers; Z or Q or R to be in the collection of intervals.

Let us call the collection of all increasing intervals, decreasing intervals and degenerate intervals as natural class of intervals.

But $N_c(Z_n) = \{[a, b] \mid a, b \in Z_n\}$ denotes the collection of all closed natural intervals as we can not order $Z_n$. Hence we have $Z_n \subseteq N_c(Z_n)$.

$N_o(Z_n) = \{(a, b) \mid a, b \in Z_n\}$ denotes the collection of all open intervals which forms the subclass of natural intervals. Clearly $Z_n \subseteq N_o(Z_n)$ $N_{oc}(Z_n) = \{(a, b] \mid a, b \in Z_n\}$ denotes the collection of all open-closed intervals forming the subclass of natural intervals.

$N_{co}(Z_n) = \{[a, b) \mid a = b\}$ denotes the class of closed open intervals forming a subclass of natural intervals.

We will denote the increasing interval [a, b] (or (a, b) or (a, b] or [a, b)) by $[a, b] \uparrow$ ($(a, b) \uparrow$, $(a, b] \uparrow$ or $[a, b) \uparrow$) and decreasing interval [a, b] (or (a, b) or (a, b] or [a, b)) by $[a, b] \downarrow$ (or $(a, b) \downarrow$ or $(a, b] \downarrow$ or $[a, b) \downarrow$). If a = b we just denote it by 'a' or 'b' as the degenerate intervals. Now we can on similar lines define $N_{oc}(Z)$, $N_{oc}(Q)$, $N_{oc}(R)$ and $N_o(Z)$, $N_o(Q)$, $N_o(R)$ and $N_c(Z)$, $N_c(Q)$, $N_c(R)$ and $N_{co}(R)$, $N_{co}(Z)$ and $N_{co}(Q)$. Now we proceed to give operations on them.

$a[x, y] \uparrow = [ax, ay] \uparrow$ if a > 0, x and y are greater than zero or less than zero.

If a < 0 and x, y both less than zero than $a[x, y] \uparrow = [ax, ay] \downarrow$.
For if a = -5; $[-3, -2] \uparrow$ then $-5[-3, -2] = [15, 10] \downarrow$.

Thus if y > x than a > 0; ay > ax; if a < 0 then ay < ax. Now $[x, y] \downarrow$ that is x > y if a > 0 than ax > ay so $[ax, ay] \downarrow$ if a < 0 than ax < ay so $[ax, ay] \uparrow$.

Thus if we consider only increasing (or decreasing) intervals then compatibility does not exist with respect to multiplication by degenerate intervals. So we take in the natural class of intervals increasing intervals decreasing intervals, and degenerate intervals. Now we show multiplication of intervals which are not in general degenerate.

Let $[x, y] \uparrow$ and $[a, b] \uparrow$ be two increasing intervals such that 0 < x < y and 0 < a < b, then $[x, y] \uparrow [a, b] \uparrow = [xa, yb] \uparrow$.

Clearly $[a, b] \uparrow [x, y] \uparrow = [ax, by] \uparrow = [xa, yb] \uparrow = [x, y] \uparrow [a, b] \uparrow$.

Since if 0 < x < y and 0 < a < b then
    xa < ya if a > 0
    xb < yb as 0 < a < b
    a < b so if x > 0
    xa < xb.

Hence xa < xb < yb, so [xa, yb] ↑ if y > x > 0 and b > a > 0.

On the other hand [-3, 8] ↑ and [-10, -2] ↑ but their product [-3, 8] ↑ [-10, -2] ↑ = [30, -16] ↓
So the assumption b > a > 0 and y > x > 0 is necessary for [ax, by] ↑. Also [-3, 8] ↑ and [-10, 2] ↑ but [-3, 8] ↑ [-10, 2] ↑ = [30, 16] ↓. So our first theorem is stated below.

**THEOREM 1.1**: *Let [a, b] ↑ and [c, d] ↑ be two intervals. In general the product [a, b] ↑ [c, d] ↑ is not an increasing interval.*

*Corollary 1.1:* If [a, b] ↑ and [c, d] ↑ such that b > a > 0 and d > c > 0 then [a, b] ↑ [c, d] ↑ = [ac, bd] ↑.

The proof is direct and hence is left as an exercise for the reader to prove.

This is true in case of open increasing intervals and half open and half closed increasing intervals.

Now we proceed onto study the decreasing closed intervals.

Let [x, y] ↓ that is x > y if a > 0 then a[x, y] ↓ = [ax, ay] ↓. Suppose a < 0 and [x, y] ↓, x > y be a decreasing interval than a[x, y] ↓ = [ax, ay] ↑. For instance a = 7 and [3, 0] ↓; 7[3, 0] ↓ = [21, 0] ↓. If a = -7 then -7 [3, 0] ↓ = [-21, 0] ↑. If [-7, -8] ↓ interval and a = 2 then 2[-7, -8] ↓ = [-14, -16] ↓.

If a = -1 then -1 [-7, -8] ↓ = [7, 8]↑.

If [4, -2] ↓ and a = 3 then 3 [4, -2] ↓ = [12, -6] ↓, suppose a = -4 then -4[4, -2] ↓ [-16, 8]↑.

Now we proceed onto see how interval product of two decreasing interval looks like.

Let [-3, 0] ↓ and [7, 2] ↓ be two decreasing intervals then [-3, 0] ↓ [7, 2] ↓ = [-21, 0]↓. [-3, -7] ↓ and [-10, -12]↓ be two decreasing intervals. Now [-3, -7] ↓ [-10, -12] ↓ = [30, 84] ↑.

In view of this we have the following theorem.

**THEOREM 1.2:** *Let [a, b]↓ and [c, d]↓ be two decreasing intervals, their product in general is not a decreasing interval.*

The proof is direct and hence is left as an exercise to the reader.

Also [-3, -7] ↓ and [7, 2]↓ are decreasing intervals, we see [-3, -7] ↓ [7, 2]↓ = [-21, -14] ↑. Also if [0, -3] ↓ and [8, 0] ↓ then we see [0, -3] ↓ [8, 0] ↓ = [0, 0] = 0 a degenerate interval.

Thus we see the product of two increasing intervals or two decreasing intervals can also be a degenerate interval. Consider [0, 7] ↑ and [-2, 0] ↑ are two increasing intervals.

Now [0, 7] ↑ [-2, 0] ↑ = [0, 0] a degenerate interval.

We also consider the general resultant of an increasing interval with a decreasing interval.

Consider [8, 2] ↓ and [-3, -10] ↑, the product is [8, 2] ↓ [-3, -10] ↓ = [-24, 20]↑. Suppose [8, 2] ↓ and [6, 9] ↑ be two intervals [8, 2] ↓ [6, 9] ↑ = [48, 18] ↓. Consider [-2, 0]↑ and [2, -2]↓ be two intervals.

$$[-2, 0] \uparrow [2, -2] \downarrow = [-4, 0]\uparrow.$$

Thus we see the product of an increasing and decreasing intervals can be a increasing interval or a decreasing interval.

Interested reader can put conditions on a, b, c and d where [a, b]↓ and [c, d]↑ so that the product is increasing interval or a decreasing interval.

Now we can derive all the results in case of open decreasing intervals and half open and half closed decreasing intervals.

Now we can add two increasing intervals, open or closed or half open or half closed.

Consider [a, b] ↑ and [x, y] ↑ where b > a and y > x be are two increasing intervals their sum [a, b] ↑ + [x, y] ↑ = [a+x, b+y]↑ as a < b and x < y then a+x < b+y. On the same lines if [a, b] ↓ and [c, d]↓ then [a, b]↓ +[c, d]↓ = [a+c, b+d] ↓ evident from the fact a>b and c>d then a+c>b+d hence the claim.

However we cannot say anything about the sum of a decreasing and an increasing interval.

For consider [-3, 7] ↑ and [3, -7]↓ intervals, their sum [-3, 7]↑ + [3, -7] ↓ = [0, 0] = 0 is a degenerate.

Consider [-3, 7]↑ and [2-5]↓ intervals; their sum [-3, 7]↑ + [2, -5]↓ = [-3+2,

7-5] = [-1, 2]↑ is an increasing interval. However their product [-3, 7] ↑ [2, -5]↓ = [-6, -35]↓ is a decreasing interval.

Consider [1, 3]↑ and [3, 1] ↓ intervals; their sum [1, 3] ↑ + [3, 1]↓ = 4 = [4, 4] is the degenerate interval.

Now consider [1, -3]↓ and [2, 20]↑ intervals then sum [1, -3]↓ + [2, 20]↑ = [1+2, -3+20] = [3, 17]↑ is an increasing interval.

Thus it is an interesting task to find conditions when a sum of an increasing and a decreasing interval is an increasing interval (and a decreasing interval).

**THEOREM 1.3:** *Let [a, b]↑ be an increasing interval (a>b) then [b, a]↓ is a decreasing interval (b>a). Now their sum is a degenerate interval a+b.*

The proof is straight forward and hence is left as an exercise to the reader.

Now all the results discussed in case of closed intervals are true in case of open increasing and decreasing intervals and half open and half closed intervals.

However (a, b) ↓ ⊆ [a, b] ↓, (a, b)↑ ⊆ [a, b]↑, (a, b]↑ ⊆ [a, b]↑ (a, b)↑ ⊆ [a, b] ↑ (a, b]↓ ⊆ [a, b]↓, [a, b) ↑ ⊆ [a, b]↑ and [a, b) ↓ ⊆ [a, b]↓. Thus [a, b] ↓, the closed interval is the largest interval for any given a and b for all the other three types of intervals are properly contained in [a, b].

Now a natural question would be subtraction of natural intervals.

It is pertinent to mention here that we are all the time talking of intervals got from Z or Q or R. For the case of $Z_n$ happens to be an entirely different one which will be discussed separately.

Let [a, b]↑ and [c, d]↑ be two increasing closed intervals a < b and c < d, then a – c < b – d or c – a < d – b one of them is true and other does not hold good.

Consider [-3, 7]↑ and [8, 10]↑ two intervals.

[-3, 7]↑ - [8, 10]↑ = [-3 -8, 7-10] = [-11 -3] ↑

But

$$[8, 10]\uparrow - [-3, 7]\uparrow = [8+3, 10-7] = [11, 3]\downarrow$$

Thus the difference of two increasing intervals can be an increasing interval or a decreasing interval or a degenerate.

For consider $[3, 14]\uparrow$ and $[2, 3]\uparrow$ two intervals $[3, 4]\uparrow - [2, 3]\uparrow = [1, 1] = 1$ and $[2, 3]\uparrow - [3, 4]\uparrow [-1, -1] = -1$ a degenerate interval.

Further we see $[a, b]\uparrow - [c, d]\uparrow \neq [c, d]\uparrow - [a, b]\uparrow$.

Now the same / similar results hold good in case of open increasing intervals and half open half closed open closed intervals which are increasing.

We will now discuss about decreasing intervals.

Consider $[-3, -20]\downarrow$ and $[3, 0]\downarrow$ intervals, $[-3, -20]\downarrow - [3, 0]\downarrow = [-6, -20]\downarrow$.

But $[3, 0]\downarrow - [-3, -20]\downarrow = [6, 20]\uparrow$. Thus we see the operation '-' on intervals is non commutative further the difference can give either a decreasing interval or an increasing interval.

We see the difference can also be a degenerate interval.

Thus we see the difference of two decreasing intervals can be a decreasing interval or a degenerate interval or an increasing interval.

We have discussed only with the closed decreasing intervals however the result is true in case of open decreasing intervals or open-closed interval or closed-open decreasing intervals.

In view of the we have the following theorem.

**THEOREM 1.4:** *Let $[a, b]\uparrow$ and $[c, d]\uparrow$ be any two increasing intervals. Their difference can be a increasing interval or a decreasing interval or a degenerate interval.*

The proof is direct and hence the reader is left with the task of proving it.

Now the closed increasing interval in the theorem 1.4 can be replaced by open increasing interval or half open – closed increasing interval or half closed-open increasing interval, and the result or conclusion of the theorem remains true.

**THEOREM 1.5:** *Let [a, b]↓ and [c, d]↓ be two decreasing intervals, their difference can be a decreasing interval or an increasing interval or a degenerate interval.*

The proof is left as an exercise to the reader.

We can replace the closed decreasing interval in theorem 1.5 by a open decreasing interval or half open-closed decreasing interval and half closed – open interval the result of theorem continues to hold good.

Now we discuss about the degenerate intervals.

**THEOREM 1.6:** *Let [x, x] and [y, y] be any two degenerate intervals their difference is always a degenerate interval.*

The reader is expected to supply the proof.

Now we will proceed onto find the sum and difference of a decreasing interval with an increasing interval or with a generate interval.

Consider [9, -7]↓ and [8, 12]↑ be a decreasing and an increasing interval respectively. Then sum [9, -7]↓ + [8, 12]↑ = [17, 5]↓ is a decreasing interval.

Now we find the difference [9, -7]↓ - [8, 12]↑ = [1, -19]↓ is a decreasing interval [8, 12]↑ - [9, -7]↓ = [-1, 19]↑ is an increasing interval.

Consider [2, 7]↑ and [13, 8]↓ two increasing and decreasing intervals respectively.

Their sum [2, 7]↑ + [13, 8]↓ = [15, 15] is a degenerate interval.

Their difference [2, 7]↑ - [13, 8]↓ = -[11, -1]↑ is an increasing interval.

Now [13, 8]↓ - [2, 7]↑ = [11, 1]↓ is a decreasing interval.

Thus we see the sum of an increasing interval with that of a decreasing interval can be a increasing interval or a decreasing interval or a degenerate interval. Likewise their difference consider [-13, -17]↓ and [-13, -2]↑ be two intervals. Their sum [-13, -17]↓ + [-13, -2]↑ = [-26, -19]↑ is an increasing interval.

Their difference $[-13, -17]\downarrow - [-13, -2]\uparrow = [0, -15]\downarrow$ is a decreasing interval $[-13, -2]\uparrow - [-13, -17]\downarrow = [0, 15]\uparrow$ is an increasing interval.

$10 = [10, 10]$, $-[3, -7]\uparrow$ be two intervals. $10 + [3, 9]\uparrow = [13, 17]\uparrow$ and $10 = [3, 7]\uparrow + [7, 3]\downarrow$.

We can have several interesting properties in this direction.

Now we will define the division of increasing / decreasing intervals. In case of degenerate interval we have division defined provided the denominator is not '0' = [0, 0].

Let $[a, b]\uparrow$ and $[c, d]\uparrow$ be two increasing intervals where $a \neq 0$, $b \neq 0$, $c \neq 0$ and $d \neq 0$.

Now

$$\frac{[a,b]\uparrow}{[c,d]\uparrow} = \left[\frac{a}{c}, \frac{b}{d}\right]$$

it may be an increasing interval or a decreasing interval or a degenerate interval.

Also

$$\frac{[c,d]\uparrow}{[a,b]\uparrow} = \left[\frac{c}{a}, \frac{d}{b}\right]$$

may be a decreasing interval or a degenerate interval or an increasing interval. We will illustrate this situation by some examples.

Let $[7, 13]\uparrow$ and $[-3, 2]\uparrow$ be two increasing intervals.

Now

$$\frac{[7,13]\uparrow}{[-3,2]\uparrow} = \left[\frac{-7}{3}, \frac{13}{2}\right]\uparrow$$

is an increasing interval.

Consider $[-2, 3]\uparrow$ and $[-7, 42]\uparrow$ any two increasing intervals

$$\frac{[-2,3]\uparrow}{[-7,42]\uparrow} = \left[\frac{2}{7}, \frac{3}{42}\right]\downarrow$$

is a decreasing interval.

Consider $[3, 5]\uparrow$ and $[3, 5]\uparrow$ clearly

$$\frac{[3,5]\uparrow}{[3,5]\uparrow} = [1, 1]$$

is a degenerate interval. That is in $[a, b]\uparrow$ then

$$\frac{[a,b]\uparrow}{[a,b]\uparrow} = [1,\ 1];$$

$a \neq 0$ $b \neq 0$. Consider $[0, 5]\uparrow$ and $[2, 4]\uparrow$ be two intervals. We see

$$\frac{[0,5]\uparrow}{[2,4]\uparrow} = \left[0, \frac{5}{4}\right]\uparrow.$$

Clearly

$$\frac{[2,4]\uparrow}{[0,5]\uparrow}$$

is not defined. Now consider $[-5, 0]\uparrow$ and $[-7, 9]$ intervals

$$\frac{[-5,0]\uparrow}{[-7,9]\uparrow} = \left[\frac{5}{7}, 0\right]\downarrow$$

is defined, where as

$$\frac{[-7,9]\uparrow}{[-5,0]\uparrow}$$

is not defined.

Now we will work with decreasing intervals.
    Let $[a, b]\downarrow$ and $[c, d]\downarrow$ be decreasing intervals

$$\frac{[a,b]\downarrow}{[c,d]\downarrow}$$

is defined if and only if $c \neq 0$ and $d \neq 0$.
    We will give examples of them. Let $[-5, -12]\downarrow$ and $[-2, -10]\downarrow$ is decreasing intervals.

Consider $\dfrac{[-5,-12]\downarrow}{[-2,-10]\downarrow} = \left[\dfrac{5}{2}, \dfrac{12}{10}\right]\downarrow$ is a decreasing interval.

Consider

$$\frac{[-2,-10]\downarrow}{[-5,-12]\downarrow} = \left[\frac{2}{5}, \frac{10}{12}\right]\downarrow$$

interval. Take $[3, 0]\downarrow$ and $[7, 2]\downarrow$ intervals

$$\frac{[3,0]\downarrow}{[7,2]\downarrow} = \left[\frac{3}{7}, 0\right]\downarrow$$

is a decreasing interval, however

$$\frac{[7,2]\downarrow}{[3,0]\downarrow}$$

is not defined.

Now consider $[8, 2]\downarrow$ and $[10, 5]\downarrow$ two decreasing intervals

$$\frac{[8,2]\downarrow}{[10,5]\downarrow} = \left[\frac{8}{10}, \frac{2}{5}\right]\downarrow$$

is a decreasing interval.
But

$$\frac{[10,5]\downarrow}{[8,2]\downarrow} = \left[\frac{10}{8}, \frac{5}{2}\right]\uparrow$$

is an increasing interval. Thus we can say if a quotient x/y is an increasing interval then the y/x is a decreasing interval, where x = [a, b] and y = [c, d]; $a \neq 0$, $b \neq 0$, $c \neq 0$ and $d \neq 0$.

Consider $[3, 1]\downarrow$ and $[-7, -10]\downarrow$ intervals.

$$\frac{[3,1]\downarrow}{[-7,-10]\downarrow} = \left[\frac{-3}{7}, \frac{-1}{10}\right]\uparrow$$

is an increasing interval.
Consider

$$\frac{[-7,-10]\downarrow}{[3,1]\downarrow} = \left[\frac{-7}{3}, 10\right]\downarrow$$

is a decreasing interval. Thus if $[a, b]\uparrow$ interval b > a then

$$\frac{1}{[a,b]\uparrow} = \left[\frac{1}{a}, \frac{1}{b}\right]\downarrow$$

is a decreasing interval.

Also if $[x, y]\downarrow$ is a decreasing interval that is x > y then

$$\frac{1}{[x,y]\downarrow} = \left[\frac{1}{x}, \frac{1}{y}\right]\uparrow$$

is a increasing interval.

All the results we have discussed for closed increasing or decreasing intervals are true in case of open increasing or decreasing intervals and degenerate intervals.

Now we will see whether the different types of operations are distributive, associative etc. We see if $[x, y]\uparrow$, $[a, b]\uparrow$ and $[c, d]\uparrow$ are increasing intervals then is $[x, y]\uparrow + ([a, b]\uparrow + [c, d]\uparrow) = ([x, y]\uparrow + [a, b]\uparrow) + [c, d]\uparrow$.

As operation addition is associative on reals so is the addition of these intervals which are defined as natural intervals are associative under addition.

Consider $[-7, 3]\uparrow$, $[-2, 0]\uparrow$ and $[8, 12]\uparrow$ be these increasing intervals

$([-7, 3]\uparrow + [-2, 0]\uparrow) + [8, 12]\uparrow$
$\quad = [-9, 3]\uparrow + [8, 12]\uparrow$
$\quad = [-1, 15]\uparrow$ \hfill (1)

Take
$[-7, 3]\uparrow + ([-2, 0]\uparrow + [8, 12]\uparrow)$
$\quad = [-7, 3]\uparrow + [6, 12]\uparrow$
$\quad = [-1, 15]\uparrow$ \hfill (2)

(1) and (2) are identical. Thus the three intervals are associative under the operation '+'. The same result holds good in case of degenerate intervals and decreasing intervals. Further the same is true in case of open and half open closed intervals.

Consider $[-7, -19]\downarrow$, $[3, 0]\downarrow$ and $[8, 0]\downarrow$ decreasing intervals $([-7, -19]\downarrow + [3, 0]\downarrow) + [8, 0]\downarrow = [-4, -19]\downarrow + [8, 0]\downarrow = [4, -19]\downarrow$. Now $[-7, -19]\downarrow + ([3, 0]\downarrow + [8, 0]\downarrow = [-7, -19]\downarrow + [11, 0]\downarrow = [4, -19]\downarrow$. Hence '+' is associative on decreasing intervals.

Since degenerate is nothing but reals or rationals or integers we see the operation '+' is trivially associative.

The subtraction operation '-' need not in general be associative in case of natural intervals. Consider $[-7, 2]\uparrow$, $[8, 10]\uparrow$ and $[-9, -20]\downarrow$ increasing and decreasing intervals; $([-7, 2]\uparrow + [8, 10]\uparrow) + [-9, -20]\downarrow = [1, 12]\uparrow + [-9, -20]\downarrow = [-8, -8]$ is degenerate interval $[-7, 2]\uparrow + ([8, 10]\uparrow + [-9, -20]\downarrow) = [-7, 2] + [-1, -10]\downarrow = [-8, -8]$ is non degenerate. Hence '+' is associative on all types of intervals.

Consider [3, 7]↑, [8, 9]↑ 0 and [-3, 0]↑ intervals.
[3, 7]↑ - ([8, 9]↑ - [-3, 0]↑) = [3, 7] ↑ - [11, 9]↑ = [-8, -2]↑.
Now ([3, 7]↑ - [8, 9]↑) – [-3, 0]↑ = [-5, -2]↑ - [-3, 0]↑ = [-2, -2] degenerate.

Thus the operation ' –' on natural intervals is non associative. We see as in case of reals the product of natural intervals is associative. We will just give a few examples as we have for [a, b]↑ [c, d]↑ and [e, f] ↑ ([a, b] ↑ [c, d]↑) [e, f] ↑ = [ac, bd]↑ [e, f]↑ or ↓ = [(ac)e, (bd)f] ↑ or ↓.

Likewise [a, b]↑ ([c, d]↑ [e, f]↑) = [a, b]↑ ([ce, df] ↑ or ↓) = [a (ce), b(df)] ↑ or ↓.

Since (ac) e = a(ce) and (bd) f = b(df) in reals the product of natural interval is associative. The only problem lies in determining whether it is increasing or decreasing or degenerative interval.

Consider
$$[3, 8] ↑, [-2, 0]↑ \text{ and } [7, 1]↓.$$
$$([3, 8]↑ [-2, 0]↑) [7, 1] ↓ = [-42, 0]↓$$
and
$$[3, 8]↑ ([-2, 0]↑ [7, 1↓) = [-42, 0]↓.$$

Thus the product happens to be a decreasing one. Another natural question would be does the product distributes over sum. The answer is yes

[a, b]↓ ([c, d] ↑ + [e, f]↓ = [ac + ae, bd + bf] (↓ or ↑)
  = [ac, bd] (↓ or ↑) + [ae, bf] (↓ or ↑).

Now we will work with $N(Z_n) = \{[a, b] \mid a, b \in Z_n\}$.

We see the working is little different.

Consider $Z_{12}$, [7, 5], 2[7, 5] = [2, 10]. 3[7, 5] = [21, 15] = [9, 3]. Also for $3 \in Z_{12}$ we see 3[3, 4] = [9, 0]. Now let [3, 4] and [6, 3] intervals. Consider [3, 4], [6, 3] = [6, 0].

Now [3, 4] + [6, 3] = [9, 7]. We see since we have no negative integers we have only modulo numbers. Only this collection is finite and we can have closed intervals $N_c [Z_n] = \{[a, b] \mid a, b \in Z_n\}$.

For example $N_c(Z_2) = \{0, 1, [0, 1], [1, 0]\}$, $N_o(Z_2) = \{0, 1, (0, 1), (1, 0)\}$; $N_{oc}(Z_2) = \{0, 1, (0, 1], (1, 0]\}$ and $N_{co} = \{0, 1, [0, 1), [1, 0)\}$. Thus we see the cardinality of all the collection is just four. $N_c(Z_3) = \{0, 1, 2, [0, 1], [0, 2], [1, 2], [2, 1], [2, 0], [1, 0]\}$ and so on and $|N_c(Z_3)| = 3^2$.

It is easily verified that $o(N_c(Z_n)) = o(N_o(Z_n)) = o((N_{oc}(Z_n)) = o(N_{co}(Z_n)) = n^2$; $n \geq 2$.

Consider $[3, 7]$, $[5, 2]$ and $[8, 1]$ intervals from $N_c(Z_9)$. $[3, 7]$ $([5, 2] [8, 1]) = [3, 7] ([4, 2]) = [3, 5]$

Also $([3, 7] [5, 2]) ([8, 1]) = [6, 5] [8, 1] = [3, 5]$. It is easy to verify that associative law with respect to multiplication is true in case of $N_c(Z_n)$ and other types of intervals.

It is further interesting to note that $Z_n \subseteq N_c(Z_n)$, $Z_n \subseteq N_o(Z_n)$ and $Z_n \subseteq N_{oc}(Z_n)$. Further we see always $[a, b]$ covers all open, half closed half open and half open half closed intervals.

Throughout this book we shall call these special type of intervals as a Natural class of intervals for ordering is not possible in $Z_n$.

Increasing, decreasing are meaningful only in case Z or Q or R and not in case of $Z_n$.

**Chapter Two**

# SEMIGROUPS OF NATURAL CLASS OF INTERVALS

In this chapter we discuss those natural intervals which have semigroup structure and also at the same time construct semigroups using this natural class of intervals. We give illustrations of them.

$N_c(Z) = \{[a, b]\downarrow, [a, b]\uparrow, a \mid a, b \in Z\}$ be the natural class of closed intervals $(N_c(Z), +)$ is a semigroup $(N_c(Z), +)$ is a commutative semigroup $[a, b]\downarrow + c = [a, b]\downarrow + [c, c] = [a + c, b + c] \downarrow .0 = [0, 0]$ acts as the additive identity. In fact $N_c(Z)$ is an infinite commutative semigroup with identity or is an infinite commutative monoid known as the natural class of closed intervals monoid.

Further $Z \subseteq N_c(Z)$. Thus the semigroup of integers is a subsemigroup of $N_c(Z)$. In view of this the following theorems are left as an exercise for the reader.

**THEOREM 2.1**: *$N_c(Z)$ is a semigroup of natural class of closed intervals under addition (a commutative monoid).*

**THEOREM 2.2**: *($N_o(Z)$, +) is a commutative monoid of infinite order.*

**THEOREM 2.3**: *($N_{oc}(Z)$, +) is a commutative monoid.*

**THEOREM 2.4**: *($N_{co}(Z)$, +) is a commutative monoid.*

If Z in Theorems 2.1 to 2.4 is replaced by $Z_n$ or Q or R, these results continue to be true.

***Example 2.1***: Let $N_o(Z_5)$ = {0, 1, 2, 3, 4, (0, 1), …, (0, 4) … (4, 3)} be the natural class of open intervals. We see ($N_o(Z_5)$, +), addition modulo 5 is a semigroup with identity 0. Clearly $o(N_o(Z_5)) = 5^2$.

*Note 1:* Since from the very context one can understand or know whether the given interval is a decreasing one or an increasing one we shall not make a specific mention of it.

*Note 2:* Since by very observation we know whether the interval is an open interval closed interval or half open half closed interval or half closed half open interval we will not make a specific mention of it. Our intervals can be increasing or decreasing or degenerate in Z or Q or R only.

***Example 2.2***: Let ($N_{oc}(Q)$, +) be the semigroup of half open half closed intervals under addition. $N_{oc}(Q)$ is of infinite order. In fact $N_{oc}(Q)$ is a commutative monoid. Clearly $Q \subseteq (N_{oc}(Q))$ is a subsemigroup of ($N_{oc}(Q)$, +).

***Example 2.3***: Let $N_c(Z_4)$ = {0, 1, 2, 3, [0, 1] [0, 2], [0, 3] [1, 2] [1, 3] [2, 3] [3, 0] [2, 0] [1, 0] [2, 1], [3, 1] [3, 2] be a semigroup under addition modulo 4.
    Consider P = {0, 2, [0, 2], [2, 0]} $\subseteq N_c(Z_4)$, P is a subsemigroup of $N_c(Z_4)$.

**DEFINITION 2.1**: *Let $N_c(Z)$ (or $N_o(Z)$ or $N_{oc}(Z)$ or $N_{co}(Z)$) be a semigroup under addition, $Z \subseteq (N_c(Z)$ is a sub-semigroup under* *addition known as the inherited or natural subsemigroup of $N_c(Z)$.*

All structures $N_o(Z)$, $N_{oc}(Z)$ and $N_{co}(Z)$ have inherited or natural subsemigroup.

Consider $N_o(2Z) \subseteq N_o(Z)$, clearly $N_o(2Z)$ is only a subsemigroup of $N_o(Z)$ and is not an inherited subsemigroup of $N_o(Z)$.

We see the definition 2.1 is true if we replace Z by Q or R or $Z_n$; $n < \infty$.

***Example 2.4**:* Let $(N_o(R), +)$ be a semigroup of natural class of open intervals. $(N_o(Q), +) \subseteq N_o(R)$ is a subsemigroup and $R \subseteq (N_o(R), +)$ is an inherited subsemigroup $Z \subseteq (N_o(R), +)$, is also an inherited subsemigroup of $(N_o(R), +)$. Infact $(N_o(R), +)$ has infinitely many inherited subsemigroups and subsemigroups.

***Example 2.5***: Let $(N_{oc}(Z_{15}), +)$ be a semigroup. Clearly $(N_o(Z_{15})) = 15^2$. Consider $Z_{15} \subseteq N_{oc}(Z_{15})$ is an inherited subsemigroup. In fact $N_{oc}(Z_{15})$ has only 3 inherited subsemigroups given by $Z_{15}$, $H = \{0, 3, 6, 9, 12\} \subseteq N_{oc}(Z_{15})$ and $P = \{0, 5, 10\} \subseteq N_{oc}(Z_{15})$. Also $((N_o(Z_{15}), +)$ has only finite number of subsemigroups. $(N_{oc}(P), +) = \{0, 5, 10, (0, 5], (0, 10], (5, 0], (10, 5], (5, 10] (10, 0]\} \subseteq N_{oc}(Z_{15})$. Here $o(N(P)) \mid o(N_{oc}(Z_{15}))$.

Similarly $N_{oc}(H) = \{0, 3, 6, 9, 12, (0, 3], (0, 6], (0, 9], (0, 12], (3, 0], (6, 0], (9, 0] (12, 0], (3, 6], (6, 3], (3, 9], (9, 3], (3, 12], (12, 3], (6, 9], (9, 6], (6, 12], (12, 6], (9, 12], (12, 9]\} \subseteq (N_o(Z_{15})$ is a subsemigroup of $N_{oc}(Z_{15})$. Clearly $o(N_{oc}(H)) \mid o(N_{oc}(Z_{15}))$.

Now we will proceed onto define semigroup under multiplication on natural class of intervals $N_o(Z)$, $(N_c(Z)$ and $N_{oc}(Z))$.

**DEFINITION 2.2**: *Let $N_c(Q)$ be the set of increasing, decreasing and degenerate closed intervals with entries from Q. $(N_c(Q), \times)$ is a semigroup called the multiplicative semigroup of natural closed intervals.*

Clearly $N_o(Q)$, $N_{oc}(Q)$, $N_{co}(Q)$ are all multiplicative semigroups of natural respective intervals.

If we replace Q by Z or R or $Z_n$, $n < \infty$ still we get them to be a multiplicative semigroup of natural intervals.

We will first illustrate them by some examples.

*Example 2.6*: Let $\{N_o(Z_3), \times\} = S = \{0, 1, 2, (0, 1), (0, 2), (1, 2)$ $(1, 0)$ $(2, 0), (2, 1), \times\}$ be a semigroup of open intervals under multiplication modulo three. S has zero divisors. 1 acts as the identity; for 1. $(0, 1) = (0, 1) \cdot (1, 2)$ $(1 \cdot 2) = (1, 1) = 1$; $(2, 1)$ $(2 \cdot 1) = (1, 1) = 1$ and $(1 \cdot 2)$ $(2 \cdot 1) = 2$ and so on.

*Example 2.7*: Let $\{N_c(Z_5), \times\} = S = \{0, 1, 2, 4, 3, [0, 1], [0, 2], [0, 3], [0, 4], [1, 2], [1, 3], [1, 4]$ , …, $[3, 4], [4, 3]\}$ is a semigroup of order 25, we have zero divisors and units. Infact S is a monoid under multiplication modulo 5.

*Example 2.8*: Let $S = \{N_o(Z), \times\}$ be a semigroup under multiplication $1 \in S$ is the identity element of S. S has zero divisors.

We now proceed onto define subsemigroups and ideals.

**DEFINITION 2.3**: *Let S be a semigroup of natural intervals under multiplication. $P \subseteq S$ is said to be a subsemigroup of natural intervals under product, if P itself is a semigroup under the operations of S.*

If in addition for every $p \in P$ and $s \in S$ ps and $sp \in P$ then we define P to be an ideal of S.

We will illustrate this situation by some examples.

*Example 2.9*: Let $S = \{N_c(Z), \times\}$ be a semigroup of special intervals under multiplication. Consider $P = \{N_c(3Z), \times\} \subseteq S$, P is a subsemigroup of special intervals of S. Infact P is an ideal of S. In fact S has infinite number of ideals.

***Example 2.10***: Let $S = \{N_o(Z_{12}), \times\}$ be a semigroup of special open intervals. Consider $W = \{(N_o(\{0, 2, 4, 6, 8, 10\}), \times\} \subseteq \{N_o(Z_{12}) \times\}$, W is a subsemigroup of S; in fact an ideal of S.

Take $T = \{N_o(\{0, 3, 6, 9\}), \times\} = \{0, 3, 6, 9, (0, 3), (0, 6), (0, 9), (3, 0) (6, 0), (9, 0), (3, 6), (6, 3) (3, 9) (9, 3), (6, 9), (9, 6), \times\} \subseteq S$, T is an ideal of S.

We see $Z_{12} \subseteq S$ is only a subsemigroup of S and is not an ideal of S. For $(1, 1) = 1 \in Z_{12} \subseteq S$ so $S \cdot 1 = S \subseteq Z_{12}$ a contradiction so $Z_{12}$ cannot be an ideal of S.

***Example 2.11***: Let $S = (N_o(Z), \times)$ be a semigroup under multiplication. $Z \subseteq S$; Z is only a subsemigroup which is the inherited subsemigroup of S. We see Z is not an ideal of S for $(0, 9) \in S$ and $2 \in Z$ we see $2(0, 9) = (0, 18) \notin Z$.

In view of this we have the following theorem.

**THEOREM 2.5**: *Let $S = (N_o(Z), \times)$ (or $N_c(Z)$ or $N_{oc}(Z)$ or $N_{co}(Z)$ or Z replaced by Q or R or $Z_n$) be a special class of open intervals. S is a semigroup under multiplication. No inherited subsemigroup of S is an ideal of S.*

Proof is direct and hence is left as an exercise for the reader to prove.

Now having seen examples of ideals we can see every semigroup of this type have ideals. To this effect we give the following example.

***Example 2.12***: Let S $\{N_o(Z_3), \times\}$ be a semigroup under multiplication. $S = \{0, 1, 2, (0, 1), (0, 2), (1, 0) (2, 0), (1, 2), (2, 1)\}$.

Clearly $\{0, 1, 2\} \subseteq S$ is an inherited subsemigroup of S. $P = \{(0, 1), (0, 2)\} \subseteq S$ is a subsemigroup of S.

$T = \{(1, 0), (2, 0)\} \subseteq S$ is also a subsemigroup of S.

$M = \{0, (0, 1), (0, 2) (2, 0), (1, 0)\} \subseteq S$ is also a subsemigroup of S. None of P or T or $Z_3$ are ideals of S as $(1, 2) \notin Z_3$, $0 \in P$ or T so P and T are not ideals $0(0, 1) = 0 \notin P$; $0(2,$

0) = 0 ∉ T. M is an ideal of S. V = {0, (0, 1), (0, 2)} ⊆ S is also an ideal of S.

W = {0, (1, 0), (2, 0)} ⊆ S is also an ideal of S. We see o(V) = 3 = o(W) and o(V)/o(S) and o(W) / o(S). But as subsemigroup o(P) / o(S) o(T) $X_o(S)$. Also o(M) = 5 and o(M) / o(S) that is 5 / 9. We see S has ideals and subsemigroups which are not ideals.

In view of this we have the following theorem.

**THEOREM 2.6**: *Let S be a natural class of decreasing, increasing and degenerate intervals. Every ideal of S is a subsemigroup but a subsemigroup need not in general be an ideal.*

This theorem is direct hence left as an exercise to the reader.

*Example 2.13*: Let S = ($N_o$ (Q), ×) = {(a, b), the collection of all increasing or decreasing or degenerate intervals with a, b ∈ Q under multiplication} be a commutative semigroup.

Consider H = {(0, a), (a, 0), 0| a ∈ Q} ⊆ S is a subsemigroup of S also H is an ideal of S. All subsemigroups are not ideals. For take P = {(a, b) where a, b ∈ Z} ⊆ S, P is a subsemigroup and not an ideal of S.

*Example 2.14*: Let S = {$N_{oc}$ ($Z_{11}$), ×} = {all intervals of the form (a, b] where a, b ∈ $Z_{11}$} be a semigroup.

S has ideals and subsemigroup. Further S has zero divisors and units and no idempotents.

Now we have seen that these intervals form a semigroup under '+' or under ×; or used in the mutually exclusive sense.

The following observations are important.
1) Semigroups of natural class of intervals are commutative with identity be it under addition or be it under multiplication.
2) All natural class of semigroups under multiplication have zero divisors.
3) All natural class of semigroups under multiplication has subsemigroups which are not ideals.

4) All natural class of semigroups under multiplication have ideals.
5) All natural class of semigroups have the inherited subsemigroups.
6) All natural class of semigroups built using $Z_n$ are of finite order.
7) Natural class of semigroups built using $Z_n$ (n not a prime) have nontrivial idempotents.
8) Natural class of semigroups built using Z or R or Q or $Z_p$ (p a prime) have no nilpotents but have zero divisors.
9) If we take only intervals of the form [0, a] (or (0, a) or [0, a] or [0, a] a ∈ $Z_n$ or $Q^+ \cup \{0\}$ or $R^+ \cup \{0\}$ or $Z^+ \cup \{0\}$ still they form a semigroup.

We will just give examples before we proceed to define more structures.

***Example 2.15***: Let S = {$N_o$ (Z), ×} be a semigroup of natural class of open increasing and decreasing intervals under multiplication.
    Consider H = {0, (0, 1), (0, 2), (0, 3), (0, 4)} ⊆ ($N_o$ ($Z_5$), ×) = S; H is a subsemigroup of S.
    Also P = {0, (1, 0), (2, 0), (3, 0), (4, 0)} ⊆ ($N_o$ ($Z_5$), ×) = S is a subsemigroup of S.

***Example 2.16***: Let S = {$N_c$ (Z), ×} be a semigroup of natural class of open increasing and decreasing intervals under multiplication.
    Consider H = {0, [0, a] | a ∈ $Z^+ \cup \{0\}$ ⊆ S; H is a subsemigroup of S. Also P = {0, [a, 0] | a ∈ $Z^+ \cup \{0\}$} ⊆ S is a subsemigroup of S. However H and P are not ideals of S.
    Now we proceed onto study stronger structure on the special class of intervals. In the first place none of these special class of intervals can be a group under multiplication. We can only study group structure under addition.

**DEFINITION 2.4**: *Let {$N_o$ (Z), ×} (or {$N_c$(Z), or $N_{oc}$ (Z) or $N_{co}$ (Z)} be a collection of natural class of intervals. For (a, b) be an increasing interval or decreasing interval or a degenerate*

*($N_o$ (Z), +) is a group 0 = (0, 0) acts as the additive identity. For every (a, b) ∈ $N_o$(Z) we take (-a, -b) ∈ $N_o$(Z) then (a, b) + (-a, -b) = (0, 0).*

Thus {$N_o$ (Z), +} is a group. Clearly ($N_o$ (Z), +) is an abelian group of infinite order.

It is important to mention here all these groups of infinite order and is torsion free.

We can in the definition 2.4 replace Z by Q or R or $Z_n$ and still the results hold good. (In $Z_n$ the concept of increasing or decreasing has no meaning).

***Example 2.17***: Let S = {$N_o$ ($Z_{12}$), +} be the group of natural class of intervals and o(S) = $12^2$.

***Example 2.18***: Let P = {$N_c$ (Z), +} be group of infinite order and is commutative.

***Example 2.19***: Let B = {$N_{oc}$ (R), +} be an infinite commutative group.

Now the task of defining a subgroup is left to the reader. However we will illustrate this concept by some examples.

***Example 2.20***: Let S = {$N_o$ (Z), +} be the group. Consider P = {$N_o$ (3Z), +} ⊆ S; P is a subgroup of S.

***Example 2.21***: Let P = {$N_{oc}$ (R), +} be the group. Consider B = {$N_{oc}$(Z), +} ⊆ P. B is a subgroup of P.

***Example 2.22***: Let M = {$N_{co}$ (Q), +} be the group. W = {$N_{co}$(5Z), +} ⊆ M is subgroup of M.

***Example 2.23***: Let T = {$N_c$ ($Z_9$), +} be a group. Consider M = {$N_c$ ({0, 3, 6}), + } = {0, 3, 6, [0, 3], [0, 6], [3, 6], [3, 0], [6, 0], [6, 3]} ⊆ T is a subgroup of T.

Now we have the concept of subsemigroups in groups. Thus we can find weaker substructures in these strong structure. Thus we find subsemigroup in natural interval groups under addition.

*Example 2.24*: Let $W = \{N_c(Z), +\}$ be a natural class of intervals. W is a group.
Consider $S = \{N_c(Z^+ \cup \{0\}), +\} \subseteq W$, S is a semigroup of W.

*Example 2.25*: Let $T = \{N_o(Q), +\}$ be a group. Take $M = \{N_o(Q^+ \cup \{0\} +\} \subseteq T$ is a semigroup of T.

Now these are given special Smarandache algebraic structures.

**DEFINITION 2.5**: *Let $M = \{N_c(Z), \times\}$ be a semigroup under multiplication. If M has a proper subset P such that P under $\times$ is a group, we define M to be a Smarandache semigroup.*

The result continues to work even if Z is replaced by Q or R or $Z_n$ ($n < \infty$). We will give some examples of them.

*Example 2.26*: Let $W = \{N_o(Z), \times\}$ be a semigroup. Consider $T = \{(1, -1), (-1, 1), (1, 1), (-1, -1), \} \subseteq W$ is a group under $\times$, evident from the following table.

| $\times$ | (1, 1) | (1, -1) | (-1, 1) | (-1, -1) |
|---|---|---|---|---|
| (1, 1) | (1, 1) | (1, -1) | (-1, 1) | (-1, -1) |
| (1, -1) | (1, -1) | (1, 1) | (-1, -1) | (-1, 1) |
| (-1, 1) | (-1, 1) | (-1, -1) | (1, 1) | (1, -1) |
| (-1, -1) | (-1, -1) | (-1, 1) | (1, -1) | (1, 1) |

Thus W is a Smarandache semigroup. $M = \{(1, 1), (-1, -1)\} \subseteq W$ is also a group of W.

*Example 2.27*: Let $T = \{N_c(2Z), \times\}$ be a semigroup. Clearly T is not a Smarandache semigroup.
Thus a natural interval semigroup may be a Smarandache semigrioup or may not be a Smarandache semigroup.

*Example 2.28*: Let $M = \{N_o(Q), \times\}$ be a semigroup under multiplication. $P = \{(0, a) \mid a \in Q^+\} \subseteq M$ is a group with (0, 1) acting as multiplicative identity. Thus M is a Smarandache semigroup (S-semigroup).

***Example 2.29***: Let $T = \{N_{oc}(R), \times\}$ be a semigroup under multiplication. $W = \{(0, a] \mid a \in Q^+\} \subseteq T$ is a group. So T is also a Smarandache semigroup.

***Example 2.30***: Let $T = \{N_c(Q), \times\}$ be a semigroup. It is easily verified that T is a S-semigorup.
None of the semigroups described in the above three examples is a group as they have many zero divisors.

***Example 2.31***: Consider $S = \{N_o(Z_7) \setminus \{0\}, \times\}$, S is a group.

***Example 2.32***: Let $S = \{N_o(Z_3) \setminus \{0\}), \times\} = \{(1, 2), 1, (2\ 1), 2\}$ is a group.

***Example 2.33***: Let $P = \{N_o(Z_5) \setminus \{0\}, \times\} = \{(1, 2)\ (1, 3)\ (1, 4), (4, 1), (4, 2)\ (4, 3), (2, 1)\ (3, 1), (2, 3)\ (3, 2)\ 1, 2, 3, 4, (3, 4), (4, 3)\}$ be a group of order 16. Clearly S is an abelian group.
If we replace the $N_o(Z_5 \setminus \{0\})$ by $N_c(Z_5 \setminus \{0\})$ or $N_{oc}(Z_5 \setminus \{0\})$ or $N_{co}(Z_5 \setminus \{0\})$ the result that it is a multiplicative group of order 16 remains true with 1 as its identity.

In view of this we have the following theorem.

**THEOREM 2.7 :** *Let $\{N_o(Z_p \setminus \{0\}), \times\} = G$ (or $N_{oc}(Z_p \setminus \{0\})$ or $N_c(Z_p \setminus \{0\})$ or $N_{co}(Z_p \setminus \{0\})$) be a group where p is a prime. $o(G) = (p-1)^2$.*

The proof is direct and hence left as an exercise to the reader.

**THEOREM 2.8:** *Let $G = \{N_o(Z_n), \times\}$ be collection of natural open intervals. G is not a group only a S-semigroup (n a composite number).*

This proof is also direct and hence left as an exercise to the reader.
We will illustrate this by some simple examples.

***Example 2.34***: Let $G = \{N_c (Z_{12}), \times\}$ be the semigroup of order $12^2$. Consider $H = \{[1, 11], [11, 1], [1, 1], [11, 11]\} \subseteq G$; H is a subgroup of G.

So G is a S-semigroup.

In view of this we give the following theorem.

**THEOREM 2.9**: *Let $\{N_c (Z_n), \times\} = G$ be a semigroup. $H = \{[1, n-1], [n-1, 1], [1, 1], [n-1, n-1]\} \subseteq G$ is a group and thus G is a S-semigroup.*

The proof is direct and hence left as an exercise to the reader.

It is clear if $N_c (Z_n)$ is replaced by $N_o (Z_n)$ or $N_{oc} (Z_n)$ or $N_{co} (Z_n)$ still the theorem holds good.

We will recall that a group (G, *) is said to be a Smarandache special definite group if G has nonempty of subset $H \subseteq G$ such that (H, *) is a semigroup.

We will illustrate this situation by some examples.

***Example 2.35***: Let $G = \{N_c (Z), +\}$ be a group. Consider $S = \{N_c (Z^+ \cup \{0\}), +\} \subseteq G$ is a semigroup. Hence G is a Smarandache special definite group.

***Example 2.36***: Let $G = \{N_{oc} (Q), +\}$ be a group. Consider $P = \{N_{oc} (Z^+ \cup \{0\}), +\} \subseteq G$, P is a semigroup. So G is a S-special definite group.

***Example 2.37***: Let $T = \{N_c (2Z), +\}$ be a group $M = \{N_{co} (Z^+ \cup \{0\}), +\} \subseteq T$ is a semigroup. Hence T is a S-special definite group.

Let G be a group of natural intervals under addition.

Let (H, +) be a subgroup of G. If, (H, +) is itself a Smarandache special definite group then we call (H, +) to be a S-special definite subgroup of G.

***Example 2.38***: Let $G = \{N_c(Q), +\}$ be a group consider $H = \{N_c (Z), +\}$ be a subgroup of G. Now $P = \{N_c (Z^+ \cup \{0\}), +\}$ is a semigroup under addition +. Thus H is a Smarandache special definite subgroup of G.

In view of this we have the following theorems.

**THEOREM 2.10** : *Let (G, +) be a group of natural class of intervals under addition. Suppose (H, +) be a subgroup of (G, +) which is a Smarandache special definite subgroup of G. Then (G, +) is a S-special definite group.*

This proof is also direct and hence left as an exercise to the reader.

**THEOREM 2.11**: *Let G = {$N_c$ (R), +} (Q or Z) be a group. G is a S-special definite group.*

If $N_c$ (R) is replaced by $N_o$(R) or $N_{oc}$ (R) or $N_{co}$ (R) still the results hold good. Also R can be replaced by Q or Z still the result is true.

Chapter Three

# RINGS OF NATURAL CLASS OF INTERVALS

Now in this chapter we proceed onto give a ring structure to the class of natural intervals obtained using Z or Q or R or $Z_n$, $n < \infty$ and study the special properties enjoyed by them.

**DEFINITION 3.1**: *Let $S = \{N_o(Z), +, \times\}$ be the collection of open decreasing and open increasing intervals together with the degenerate intervals on Z. Clearly $(N_o(Z), +)$ is an abelian group and $(N_o(Z), \times\}$ is a semigroup with identity $1 = (1, 1)$. Further $(a, b) \times [(c, d) + (e, f)] = (ac + ae, bd + bf) = (a, b) \times (c, d) + (a, b) \times (e, f) = (ac, bd) + (ae, bf) = (ac + ae, bd + bf)$.*

*That is multiplication distributes over addition $S = (N_o(Z), +, \times\}$ is defined as the ring of the class of natural open intervals.*

It is important to note that if Z in definition 3.1 is replaced by R or Q still we have the ring structure.
 Also if the open intervals is replaced by closed intervals or half open-half closed intervals or half closed-half open intervals still the structure is a ring.

We will illustrate it by some simple examples.

*Example 3.1*: Let S = {$N_c(R)$, +, ×} be a ring. Clearly S is a commutative ring with unit. S has infinitely many zero divisors but has no idempotents or nilpotents in it. S has also infinite number of elements which are units.

*Example 3.2*: Let S = {$N_{oc}(Z)$, +, ×} be a ring. S is a commutative ring with unit of infinite order S has no idempotents or nilpotent elements in it, however S has infinite number of zero divisors. S has no units.

*Example 3.3*: Let S = {$N_{co}(Q)$, +, ×} be the ring. S is a commutative ring with unit. S has zero divisors but S does not contain nilpotent elements or idempotent elements. In fact infinite number of elements in S are invertible.

*Example 3.4*: Let S = {$N_c(Z_{12})$, +, ×} be a ring. S is of finite order S has units, idempotents, nilpotents and zero divisors. S is a commutative ring with unit. [0, 4] × [0, 4] = [0, 4] is an idempotent in S.

[0, 6] × [0, 6] = [0, 0] = 0 is a nilpotent element in S.
[3, 4] × [4, 3] = [0, 0] is a zero divisor in S.
[1, 11] [1, 11] = [1, 1] = 1 unit.
[11, 11] [11, 11] = [1, 1] = 1 is the degenerate unit.

Also [11, 1] × [11, 1] = [1, 1] is a unit and [0, 8] × [5, 0] = 0 is a zero divisor. Clearly $o(S) = 12^2 = 144$.

*Example 3.5*: Let S = {$N_c(Z_7)$, + ×} be a ring of order $7^2$ elements in it. S has zero divisors but has no idempotents or nilpotents in it. However S has units.

Now we have seen finite and infinite rings of natural intervals, for any given $Z_n$ we have 4 different rings of order $n^2$; all of them are commutative with identity 1 = (1, 1). We can as in case of rings define in case of natural class of interval rings also the notion of ideals and subrings. We will only illustrate this situation and leave the task of defining these concepts to the reader.

***Example 3.6***: Let $S = \{N_o(Z), +, \times\}$ be a ring. $W = \{N_o(3Z), +, \times\} \subseteq S$ is a subring of S. In fact W is an ideal of S.

***Example 3.7***: Let $S = \{N_c(Q), +, \times\}$ be a ring. $T = \{N_c(Z), +, \times\}$ is a subring of S. Clearly T is not an ideal of S.

***Example 3.8***: Let $S = \{N_o(Z_{24}), +, \times\}$ be a ring. $T = \{N_o(0, 2, 4, 6, \ldots, 20, 22\}), +, \times\} \subseteq S$ is a subring of S as well as ideal of S.

***Example 3.9***: Let $S = \{N_{oc}(R), +, \times\}$ be a ring. $T = \{(N_{oc}(Q), +, \times\}$ is a subring and not an ideal of S.

***Example 3.10***: Let $S = \{N_{co}(Z_7), \times, +\}$ be a ring.
   $T = \{[0, 1], [0, 2], [0, 3], [0, 4], [0, 5], [0, 6], 0\} \subseteq S$ is a subring of S and also an ideal of S.
   $P = \{0, 1, 2, 3, 4, 5, 6\} \subseteq S$ is the inherited subring of S. But P is not an ideal of S.

In fact no inherited subring of natural class of intervals is an ideal.

**THEOREM 3.1**: *Let S be a natural class of interval of rings. Every ideal of S is a subring of S but subring of S in general need not be an ideal of S.*

The proof is direct and hence is left as an exercise for the reader.

**THEOREM 3.2**: *Let S be ring of natural intervals (open or closed or open closed or closed open). Let P be the inherited subring of S. P is not an ideal of S.*

This proof is also a left as an exercise to the reader.

***Example 3.11***: Let $S = \{N_o(Z), +, \times\}$ be a ring.
   $P = (Z, +, \times\} \subseteq S$ is the inherited subring of S.

T = {(0, a) } a ∈ Z, +, ×} ⊆ S is a subring S, which is also an ideal of S.

M = {(a, 0) | a ∈ Z, +, ×} ⊆ S is a subring of S which is also an ideal of S.

In fact T is isomorphic to P and P is isomorphic to M and T is isomorphic to M are isomorphic as subrings.

**THEOREM 3.3**: *Let $S = \{N_o(Z), +, \times\}$ ($N_o(R)$ or $N(Z_n)$ or $N_o(Q)$) be a ring S has nontrivial ideals.*

We just give hint of the proof.

T = {(0, a) } a ∈ Z} ({(0, a) } a ∈ R} or {(0, a) | a ∈ $Z_n$} or {(0, a) } a ∈ Q}) ⊆ S is a nontrivial ideal of S · P = {(a, 0)| a ∈ Z} ({(a, 0) | a ∈ Q or {(a, 0) | a ∈ R}, {(a, 0) | a ∈ $Z_n$} are non trivial ideals of S.

Thus these structures are rich as substructures. They are not simple structures.

None of them are integral domains but these substructures contain subrings.

***Example 3.12***: Let $S = \{N_o(Z), +, \times\}$ be a ring, T = {0, (0, a) | a ∈ Z, +, ×} ⊆ S; T is an integral domain of S. P = {(a, 0) } a ∈ Z, +, ×} ⊆ S is again an integral domain of S. W = {a ∈ Z, +, ×} ⊆ S is an integral domain. The integral domains T and P are also ideals of S. M = {(0, a) | a ⊆ nZ, +, ×} ⊆ S is an integral domain of S as well as an ideal of S.

Now we proceed onto give examples of maximal ideals and minimal ideals in S, the ring of natural intervals using Z or Q or R or $Z_n$, n < ∞.

***Example 3.13***: Let $S = \{N_c(Z), +, \times\}$ be a ring. Consider P = {[0, a] | a ∈ Z} ⊆ S, is an ideal of S infact the maximal ideal of S. T = {[a, 0] | a ∈ Z} ⊆ S is an ideal of S which is a maximal ideal of S. M = {[0, a] | a ∈ 2Z} ⊆ S is not a maximal ideal of S as M ⊆ P ⊆ S. W = {[a, 0] | a ∈ 3Z} ⊆ S is not a maximal ideal of S as W ⊆ T ⊆ S.

Clearly {0} is the only minimal ideal of S.

*Example 3.14*: Let $S = \{N_c(Z_{30}), \times, +\}$ be the ring. Consider $I = \{N_c(\{0, 10, 20\}), \times, \} \subseteq S$; I is an ideal of S. $J = \{N_c(\{0, 15\}), \times, +\} \subseteq S$ is an ideal of S. Both I and J are minimal ideals of S.

Consider $K = \{N_c(\{0, 2, 4, \ldots, 28\}), +, \times\} \subseteq S$; K is a maximal ideal of S. Take $P = \{N_c(\{0, 3, 6, \ldots, 27\}), +, \times\} \subseteq S$, P is a maximal ideal of S.

Thus we have seen the concept of minimal and maximal ideals of the ring S of finite order.

*Example 3.15*: Let $M = \{N_o(Z_{19}), +, \times\}$ be a ring.

Consider $P = \{0, 1, 2, \ldots, 18\} \subseteq M$, P is only the inherited subring of M. Take $V = \{(0, 1), (0, 2), \ldots, (0, 18), 0\} \subseteq M$; V is a subring of M as well as ideal of M.

Take $T = \{(0, 0), (1, 0), (2, 0), \ldots, (18, 0)\} \subseteq M$, T is a subring of M as well as ideal of M. The ideals T and V are both maximal and minimal ideals of S.

In view of this we have the following theorem.

**THEOREM 3.4**: *Let $S = \{N_o(Z_p), +, \times\}$ (or $N_c(Z_p)$ or $N_{oc}(Z_p)$ or $N_{co}(Z_p)$; p a prime) be a ring. S has only two ideals $I_1$ and $I_2$ both are both maximal and minimal given by*
  $I_1 = \{(0, a) \mid a \in Z_p\} \subseteq S$ *and*
  $I_2 = \{(a, 0) \mid a \in Z_p\} \subseteq S$.

The proof is direct and hence is left as an exercise to the reader.

The ideals which we have discussed are only principal ideals of S. Now as in case of usual rings, we can define in case of natural class of interval rings also the concept of Smarandache rings.

Clearly the natural class of interval rings can never be fields as they contain zero divisors for we include intervals of the form (0, a) and (a, 0) is always a zero divisor. Thus only the inherited subring contained in the ring of natural class of intervals can be a field.

We will first illustrate this situation by some examples.

*Example 3.16*: Let $S = (N_o(Q), \times, +)$ be a ring. Clearly S is not a field as it has infinite number of zero divisors. Consider P = {(0,

a) | a ∈ Q} ⊆ S; P is a field called the natural interval of field of rationals.

Take V = {(a, 0) } a ∈ Q} ⊆ S; V is also an integral domain we see V and P are isomorphisic as integral domain however V contains decreasing intervals of a special form and P increasing intervals of the special form.

However the natural subsemigroup Q ⊆ S is the prime field of characteristic zero. Thus S is a Smarandache ring.

*Example 3.17*: Let R = {$N_c(Z)$ +, ×} be the ring. Clearly R has no subring which is field. Hence R is not a Smarandache ring.

Thus we have rings of natural intervals which are not Smarandache rings (S-rings).

*Example 3.18*: Let S = {$N_o(3Z)$, +, ×} be a ring which is not a S-ring.

*Example 3.19*: Let V = {$N_{oc}(17Z)$, +, ×} be a ring which is not a S-ring.

*Example 3.20*: Consider P = {$N_{co}(8Z)$, +, ×}, P is ring which is not a S-ring.

*Example 3.21*: Let W = {$N_c(nZ)$, +, ×} (n < ∞) be a ring. W is not a S-ring.

Thus we have an infinite collection of natural interval rings which are not S-rings.

*Example 3.22*: Let S = {$N_c(Z_p)$, +, ×}, p a prime be a ring. Clearly S is a S-ring.

*Example 3.23*: Let V = {$N_c(Z_{11})$, +, ×} be a ring. V is a S-ring.

*Example 3.24*: Let V = {$N_{oc}(Z_{19})$, +, ×} be a ring. V is a S-ring.

In view of this we have the following theorems.

**THEOREM 3.5**: *Every natural interval ring need not in general be a S-ring.*

The proof is direct and hence is left as an exercise for the reader.

**THEOREM 3.6**: *Let $V = \{(N_o(Z_p), +, \times\}$ (or $N_c(Z_p)$ or $N_{oc}(Z_p)$ or $N_{co}(Z_p)\}$ (p a prime) V is a S-ring.*

Proof is simple hence the task of proving it is left for the reader as an exercise.

We see we have infinite number of ring of natural intervals which are S-rings. Now we proceed onto give more examples of S-rings.

*Example 3.25*: Let $V = \{N_o(Z_{12}), +, \times\}$ be a ring. Consider $H = \{0, 4, 8\} \subseteq V$, H is a field in which 4 is the multiplicative identity of H. So V is a S-ring. Consider $P = \{0, (0, 4), (0, 8)\} \subseteq V$, P is a field. Also $T = \{0, (4, 0), (8, 0)\} \subseteq V$, T is a field.

Thus even in $Z_p$ if p is not a prime still the natural class of interval rings can be a S-ring.

We have just shown rings which are S-rings we can also have the notion of Smarandache zero divisors and Smarandache idempotents of these rings. For the definition refer [16], where the elements are replaced by intervals.

This is a matter of routine however we will illustrate it by some examples.

*Example 3.26*: Let $S = \{N_c(Z_{20}) +, \times\}$ be a ring.
Consider $[0, 10], [0, 16] \in S$ we see
$$[0, 5] \times [0, 16] = 0 \pmod{20}$$
$$[0, 6] \times [0, 20] = 0 \pmod{20}$$
$$[0, 6] \times [0, 5] \neq 0 \pmod{20}.$$
Thus [0, 10], [0, 16] is a S-zero divisor of S.

*Example 3.27*: Let $S = \{(N_c(Z_{30}) +, \times\}$ be the ring [0, 6], [0, 10], [0, 15], [0, 16] [0, 21] and [0, 25] are nontrivial idempotents of S. [0, 6] and [6, 0] are S-idempotents of S as [0, 6] × [0, 6] = [0, 6] are S-idempotents of S as [0, 6] × [0, 6] = [0, 6] (mod 30) and [6, 0] [6, 0] = [6, 0] (mod 30). [0, 24] [0, 24] = [0, 6] (mod 30) and [0, 6] [0, 24] = [0, 24] (mod 30). Similarly [24, 0] [24, 0] = [24, 0] mod 30.

Interested reader can find any number of examples and also prove related results [16].

We can define as in case of ring for these natural class of interval rings also the quotient ring. We will illustrate this by some simple examples.

***Example 3.28***: Let $S = \{N_o (Z_3), \times, +\}$ be a ring $I = \{0, (0, 1), (0, 2)\} \subseteq S$ be an ideal of S. $S/I = \{I, 1 + I, (1, 0) + I, (2, 0) + I, (2, 1) + I, 2 + I, (1, 2) + I\}$ is the quotient ring. $(1, 0) + I + (1, 2) + I = (2, 2) + I = 2 + I$. Further $((1, 0) + I) ((1, 2) + I) = (1, 0) + I$.

Thus we see S/I has seven elements and is a ring of characteristic three. Also by taking the ideal $J = \{0, (1, 0), (2, 0)\} \subseteq S$ we the get the quotient $S/J = \{J, 1 + J, 2 + J, (0, 1) + J, (0, 2) + J, (2, 1) + I, (1, 2) + J\}$ is again a ring of characteristic three.

***Example 3.29***: Let $S = \{N_c (Z_6), \times, +\}$ be a ring. $S = \{0, 1, 2, \ldots, 5,$

$$[0, 1], [0, 2], [0, 3], [0, 4], [0, 5],$$
$$[1, 5], [1, 2], [1, 3], [1, 4], [1, 0],$$
$$[2, 0], [2, 1], [2, 3], [2, 4], [2, 5],$$
$$[3, 0], [3, 1], [3, 2], [3, 4], [3, 5],$$
$$[4, 0], [4, 1], [4, 2], [4, 3], [4, 5],$$
$$[5, 1], [5, 2], [5, 3], [5, 4], [5, 0]\}.$$

$P = \{[0, a] \mid a \in Z_6\} \subseteq S$ is an ideal of S.
$S/P = \{P, 1 + P, 2 + P, 3 + P, 4 + P, 5 + P, [1, 2] + P, \ldots, [5, 4] + P\}$ is the ring with 31 elements in it. Clearly S/P is also a ring of characteristic six.

***Example 3.30***: Let $\{N_{oc} (Z_4), \times, +\}$ be a ring. Consider $I = \{0, (0, 1], (0, 2], (0, 3]\} \subseteq S$, I is an ideal of S.
Now $S/I = \{I, 1 + I, 2 + I, 3 + I, (1, 0] + I, (2, 0] + I, (3, 0] + I, (1, 2] + I, (2, 1] + I, (1, 3] + I, (3, 1] + I, (2, 3] + I, (3, 2] + I\}$ is a ring with 13 elements in it but its characteristic is four and has zero divisors.

Thus by the quotient rings we can get several special interval rings of prime order but characteristic prime or nonprime. However all $Z_n$ does not yield a field.

Consider the following examples.

**Example 3.31**: Let $S = \{N_o(Z_5), \times, +\}$ is the ring of order 25. Consider $I = \{0, (0, 1), (0, 2) (0, 3), (0, 4)\} \subseteq S$, I is an ideal of S.

Now $S/I = \{I, 1 + I, 2 + I, 3 + I, 4 + I, (1, 0) + I, (2, 0) + I, (3, 0) + I, (4, 0) + I, (2, 1) + I, (1, 2) + I, (3, 1) + I, (1, 3) + I, (1, 4) + I, (2+3) + I, (4, 1) + I, (3, 2) + I, (4, 2) + I, (2, 4) + I, (4, 3) + I, (3, 4) + I\}$ is the quotient ring and o(S/I). 21 of characteristic five.

**Example 3.32**: Let $S = \{N_o(Z_{10}), +, \times\}$ be the ring. $J = \{0, (0, a) \mid a \in Z_{10}\} \subseteq S$ be an ideal of S. Now $S/J = \{J, (a, 0) + J, (a, b) + J, a, b \in Z_{10} \setminus \{0\}\}$ is a ring of characteristic 10. We see $o(S/J) = 10^2 + 1 - 10 = 91$. Clearly S/J has both zero divisors and idempotents. Consider $((5, 0) + J) ((2, 0) + J) = 0 + J = J$ is a zero divisor. Also $((5, 0) + J) ((5, 0) + J) = (5, 0) + J$ is an idempotent in S/J.

We can get rings of prime order but with different non prime characteristic.

**Example 3.33**: Let $S = \{N_c(Z_9), +, \times\}$ be a ring with 81 elements in it. Let $P = \{0, [0, a] \mid a \in Z_9 \setminus \{0\}\} \subseteq S$ be an ideal. Consider the quotient ring $S/P = V = \{0, 1, 2, \ldots, 8, [a, b] + J \mid a \in Z_9 \setminus \{0\}$ and $b \in Z_9\}$. Clearly $o(S/P) = 73$. V is a ring of characteristic 9. V has zero divisors and units in it $[6, 3] + J \in V$ is a zero divisor of V for $([6, 3] + J) ([6, 3] + J) = ([0, 0] + J) = J$ $([8, 1] + J ([8, 1] + J) [1, 1] + J = 1 + J$ is a unit of V. However V has no idempotents in it.

**Example 3.34**: Let $M = \{N_o(Z_7), +, \times\}$ be a ring of characteristic 7 and order 49. Let $J = \{(0, a) \mid a \in Z_7|$ be an ideal of M. $M/J = V = \{a + J, (a, 0) + J, (x, y) + J$ where $a \in Z_7 \setminus \{0\}$ and $x, y \in Z_7 \setminus \{0\}\}$ is the quotient ring of characteristic 7. V is of order 43.

V has no zero divisors but has units and is commutative. We see $(a, 0) + J$ in V are non invertible elements but also are not zero divisors. We see any $(0, a)$ for $a \in Z_7 \setminus \{0\}$ is not invertible or has no inverse but it does not also contribute to idempotents or zero divisors. Further we see V is of characteristic seven.

***Example 3.35***: Let $M = \{N_c(Z_{11}), +, \times\}$ be a ring. $J = \{[0, a] \mid a \in Z_{11}\} \subseteq M$ be an ideal of M. Now let $V = M/J = \{[a, 0] + J, [a, b] + J, a + J$ where $a, b \in Z_{11} \setminus \{0\}\}$ where M/J is a ring of order $11^2 - 11 + 1 = 111$.

Thus we see when we use the field $Z_3$ we get a ring D of order $3^2 - 3 + 1 = 7$ and characteristic of D is 3. However D is not a field when $Z_4$ is used we get for the ideal $J = \{0, (0, 1), (0, 2), (0, 3) \in \{N_o(Z_4), +, \times\} = S$ the quotient ring $S/J = D$ to be only a ring and not an integral domain of order $4^2 - 4 + 1 = 13$ of characteristic four.

When $\{N_c(Z_5), +, \times\} = P$ is used as a ring we using the ideal $I = \{[0, a], a \in Z_5\}$ we get P/I to be an integral domain of order $5^2 - 5 + 1 = 21$ of characteristic five. When $\{N_{oc}(Z_6), +, \times\} = M$ is used as the ring for the idea $J = \{(0, a] \mid a \in Z_6\}$, $M/J = D$ is a ring of order $6^2 - 6 + 1 = 31$ of characteristic 6. This has zero divisors and also subring of order 6. When we use the ring $K = \{\{N_{co}(Z_7), +, \times\}$ and the ideal $W = \{[0, a) \mid a \in Z_7\} \subseteq K$. The quotient ring K/W is of order $7^2 - 7 + 1 = 43$ and of characteristic 7.

Clearly K/W is ring. When $P = \{N_c(Z_8), +, \times\}$ is used as the ring, for the ideal $M = \{[0, a] \mid a \in Z_8\}$ we get the quotient ring P/M to be a ring with zero divisors and $o(P/M) = 8^2 - 8 + 1 = 57$.

For $S = \{N_o(Z_9), \times, +\}$ we get a ring of order $9^2 = 81$ and using the ideal $J = \{(0, a) \mid a \in Z_9\} \subseteq S$ we have the quotient ring to be of order 73 of characteristic nine and has zero divisors. Further we have elements in $\{N_c(Z_{11}), +, \times\} = S$ such that $a^{n(a)} = a$.

For instance $[0, 5] \in N_c(Z_{11})$ we see $[0, 5]^6 = [0, 5]$. $[0, 3]^6 = [0, 3]$ and so on. Now consider $J = \{[0, a] \mid a \in Z_{11}\} \subseteq S$, J is an ideal of S. Consider the quotient ring $S/J = D = \{[a, 0] + J, [b, a] + J \mid a, b \in Z_{11}\}$. D is ring order 111. Clearly D has no zero divisors further D is characteristic 11 and every element d in D is of the form $d^{n(d)} = d$ for a suitable integer $n(d) > 0$.

**Example 3.36**: Let $V = \{N_o (Z_{13}), \times, +\}$ be a ring. Consider $J = \{(0, a) \mid a \in Z_{13}\} \subseteq V$ be the ideal of V. Let $D = V/J = \{J (a, 0) + J, (a, b) + J \; a, b \in Z_{13} \setminus \{0\}\}$ be the quotient ring of characteristic thirteen. $o(D) = 13^2 - 13 + 1 = 157$. D is ring with 157 elements and has no zero divisor and is of characteristic thirteen. Consider $(9, 0) \in D$ $(9, 0)^4 = (9, 0)$. Consider $(12, 0) \in D$ $(12, 0)^3 = (12, 0)$ $(11, 0)^{13} = (11, 0)$. Consider $(3, 0)^4 = (3, 0)$ and so on.

**Example 3.37**: Let $S = \{N_o (Z_7), +, \times\}$ be a ring of natural open intervals with entries from $Z_7$. S is a commutative ring with unit $1 = (1, 1)$ of finite order. $o(S) = 49$.

$S = \{0, 1, 2, \ldots, 6, (0, 1), (0, 2), \ldots, (0, 6), (1, 0), (2, 0), \ldots, (6, 0), (1, 2), (2, 1), (1, 3), \ldots, (6, 1), \ldots (5, 6), (6, 5)\}$.

We see $\{(0, a) \mid a \in Z_7\} = J$ is an ideal of S. Also $\{(a, 0) \mid a \in Z_7\} = I$ is an ideal of S. Further $(0, a)$ and $(b, 0)$, $a, b \in Z_7 \setminus \{0\}$ are non invertible elements of Z. They are not units of S. Consider the quotient ring $S/J = D = \{J, a + J, 1, + J (a, 0) + J, (a, b) + J / a, b \in Z_7 \setminus \{0\}\}$. It is easily verified $S/J = D$ has only 43 elements. D is commutative ring and has no zero divisors. $(a, 0) + J$ are non invertible but are such that $[(a, 0) + J]^n = (a, 0) + J$ where n is an integer.

Further each $(a, b) + J$ with $a, b \in Z_7 \setminus \{0\}$ has inverse $(c, d) + J$ such that $((a, b) + J) ((c, d) + J) = (1, 1) + J = 1 + J$.

**Example 3.38**: Let $S = \{N_o (Z_5), +, \times\}$ be the ring. $S = \{0, 1, 2, 3, 4, (1, 0), (2, 0), (3, 0), (4, 0), (0, 1), (0, 2), (0, 3), (0, 4) (1, 2), (1, 3), (1, 4), (2, 3), (2, 4), (3, 4), (4, 1), (4, 2), (4, 3), (2, 1), (3, 2), (3, 1)\}$ is a commutative ring with $(1, 1) = 1$ as unit and of order 25.

Consider $J = \{(0, 1), (0, 2), (0, 3), (0, 4), 0\} \subseteq S$ is an ideal of S. $S/J = D = \{J, 1 + J, (1, 0) + J, (2, 0) + J, (3, 0) + J, (4, 0) + J, (1, 2) + J, (1, 3) + J, (1, 4) + J, (2, 1) + J, (2, 3) + J, (2, 4) + J, (3, 4) + J, (3, 1) + J, (4, 1) + J, (3, 2) + J, (4, 2) + J, (4, 3) + J, + J, 3 + J, 4 + J\}$ is a commutative ring with 2l elements in it and of characteristic 5. J acts as the additive identity and $1 + J$ acts as the multiplicative identity.

$(1, 0) + J$ is such that

$$((1, 0) + J) ((1, 0) + J) = (1, 0) + J.$$
$$((2, 0) + J)^5 = (2, 0) + J,$$
$$((3, 0) + J)^5 = (3, 0) + J,$$
$$((4, 0) + J)^3 = (4, 0) + J,$$
$$(2+J)^5 = 2 + J, (3 + J)^5 = 3 + J,$$
$$(4+J)^3 = 4 + J, ((1, 2) + J)^5 = (1, 2) + J, ((1, 3) + J)$$
$$= (1, 3) + J,$$
$$((1, 4) + J)^3 = (1, 4) + J, ((2, 1) + J)^5,$$
$$= (2, 1) + J, ((3, 1) + J)^5$$
$$= (3, 1) + J,$$
$$((4, 1) +J)^3 = (4, 1) + J,$$
$$((2, 3) + J)^5 = (2, 3) + J.$$
$$((3, 2) + J)^5 = (3, 2) + J,$$
$$((2, 4) + J)^5 = (4, 2) + J,$$
$$((3, 4) + J)^5 = (3, 4) + J \text{ and}$$
$$((4, 3) + J)^5 = (4, 3) + J.$$

***Example 3.39***: Let R = {$N_c$ ($Z_3$), +, ×} be a ring R = {0, 1, 2, [0, 1], [0, 2], [1, 0], [2, 0], [1, 2], [2, 1] is a ring with nine elements and characteristic of R is three.

R is a finite commutative ring with unit and has nontrivial zero divisors given by [0, 1] [2, 0] = 0, [0, 2] [2, 0] = 0, [0, 1] [1, 0] = 0 and [1, 0] [0, 2] = 0. Consider P = {0, [0, 1], [0, 2]} ⊆ R P is an ideal of R. Take the quotient ring of R with the ideal P.

R/P = D = {0, 1+P, 2+P, [1, 0] + P, [2, 0] + P, [1, 2] + P, [2, 1] + P}; D is a commutative ring of order seven and of characteristic three. Clearly D is ring such that $(2+P)^3 = 2 +$ , $([2, 0] + P)^3 = [2, 0] + P, ([1, 2] + P) ([2, 1] + P) = 2 + P. ([1, 2] + P)^2 = [1, 1] + P = 1 + P.$

$([2, 1] + P)^2 = [1, 1] + P = 1 + P.$

$([1, 2] + P)^3 = [1, 2] + P$ and $([2, 1] + P)^3 = [2, 1] + P.$ Clearly D is not a field as [1, 0] + P and [2, 0] + P have no inverse. Here also D is a only a ring of characteristic three and is not a field.

We have the following theorems.

**THEOREM 3.7**: *Let $S = \{N_o (Z_p), +, \times\}$ (or $N_c (Z_p)$ or $N_{oc} (Z_p)$ or $N_{co}(Z_p))$, p a prime;*

*(1) S is a commutative ring with unit of order $p^2$.*
*(2) S has zero divisors.*
*(3) $J = \{(0, a)/ a \in Z_p\} \subseteq S$ is an ideal of S.*
*(4) $K = \{(a, 0) / a \in Z_p\} \subseteq S$ is an ideal of S.*
*(5) $S/J = D$ and $V = S/K$ are finite ring of order $p^2 - p + 1$ and of characteristic p.*
*(6) Every element in D and V are of the form $x^{n(x)} = x$, $n(x) > 1$. $n(x) > 1$.*

The proof is direct and hence left for the reader as an exercise.

**THEOREM 3.8** *Let $S = \{N_c (Z_n), +, \times\}$ (or $N_o (Z_n) +, \times\}$ (or $N_o (Z_n)$ or $N_{oc} (Z_n)$ or $N_{co} (Z_n))$, n a composite number, S is a commutative ring with zero divisors, units and idempotents and $o(S) = n^2$.*

*(1) $J = \{[0, a] / a \in Z_n\} \subseteq S$ is an ideal of S.*
*(2) $K = \{[a, 0] / a \in Z_n\} \subseteq S$ is an ideal of S.*
*(3) $D = S/J$ and $V = S/K$ are quotient rings of order $n^2-n + 1$.*
*(4) D and V have zero divisors, units, nilpotents and idempotents.*

The proof is also direct and hence is left as an exercise for the reader to prove.

*Example 3.40*: Let $S = \{N_o (Z_{12}), +, \times\}$ be a ring. Clearly S is a ring with unit and is commutative having zero divisors, units, idempotents and nilpotents. Order of S is $12^2 = 144$. Consider the ideal $J = \{(0, a) / a \in Z_{12}\} \subseteq S$, I is an ideal of S. Now let D $= S/J = \{J, a + J, (a, 0) + J, (a, b) + J \mid a, b \in Z_{12} \setminus \{0\}\}$ be the quotient ring. D has zero divisors, nilpotents and idempotents in it.

Consider $(3, 4) + J, (4, 3) + J$ in D is such that $((3, 4) + J) ((4, 3) + J) = J$, J is the zero of D. Consider $(4, 4) + J$ in D; clearly $((4, 4) + J) ((4, 4) + J) = (4, 4) + J$ is an idempotent of D.

Consider $(6, 0) + J$ in D $((6, 0) + J) ((6, 0) + J) = J$ is a nilpotent element of D. $(4, 0) + J$ in D is an idempotent of D. $(11, 11) + J$ is a unit in D. Thus D has idempotents, units, nilpotents and zero divisors. Clearly D is not an integral domain but a ring of order $12^2 - 12 + 1$.

***Example 3.41***: Let $W = \{N_c (Z_{53}), +, \times\}$ be a ring of order $53^2$. Clearly S is a commutative ring with unit and with zero divisors.

Consider $J = \{(0, a) \mid a \in Z_{53}\} \subseteq S$. J is an ideal of S. $S/J = D$ is a ring of order $53^2 - 53 + 1$. Clearly D is not a field for $(a, 0) + J$ has no inverse for all $a \in Z_{53} \setminus \{0\}$. Every element $(a, b) + J$ has inverse. Thus D is a ring of order $53^2 - 53 + 1$ and of characteristic 53.

Thus we have infinite number of finite rings of order $p^2 - p + 1$ for every prime p and of characteristic p. We call a ring to be a Smarandache ring if D contains a proper subset $P \subseteq D$ with P a field.

***Example 3.42***: Let $S = \{N_o (Z_7), +, \times\}$ be a ring. Consider $J = \{(0, a) \mid a \in Z_7\}$ be an ideal of S. $D = S/J = \{J, a+J, (a, 0) + J, (a, b) + J \mid a, b \in Z_7 \setminus \{0\}\}$ be the quotient ring. Clearly D is a ring $\{a + J\} \mid a \in Z_7\}$ is a field contained in D, hence D is a Smarandache ring.

In view of this we have the following theorem.

**THEOREM 3.9**: *Let $S = \{N_c (Z_p), \times, +\}$ be a ring, p a prime. $J = \{[0, a] \mid a \in Z_p\}$ be an ideal. $S/J = D$ is a Smarandache ring.*

Proof is easy hence left for the reader as an exercise.

***Example 3.43***: Let $S = \{N_c (Z_{10}), +, \times\}$ be a ring. $J = \{[a, 0] \mid a \in Z_{10}\} \subseteq S$ be an ideal of S. Consider the quotient ring $S/J = D = \{a + J, J, [0, a] + J, [a, b] + J \mid a, b \in J\}$.

***Example 3.44***: Let $S = \{N_o (Z), +, \times\}$ be a ring. Let $J = \{(0, a) \mid a \in Z\} \subseteq S$ be an ideal of S. $S/J = \{J, a + J, (a, 0) + J, (a, b) + J \mid a \in Z \setminus \{0\}\}$ is a quotient ring. S/J is a ring.

***Example 3.45***: Let $S = \{N_o(Z), +, \times\}$ be a ring. Let $J = \{N_o(3Z), +, \times\}$ be an ideal of S. $S/J = \{J, 1 + J, 2 + J, (1, 0) + J, (0, 1) + J, (2, 0) + J, (0, 2) + J, (1, 2) + J (2, 1) + J\}$ is a ring S/J is isomorphic with $\{N_o(Z_3), +, \times\} = \{0, 1, 2, (0, 1), (1, 0), (1, 0), (2, 0), (0, 2), (1, 2), (2, 1)\}$.

***Example 3.46***: Let $S = \{N_o(Z), +, \times\}$ be a ring. Let $J = \{N_o(4Z), +, \times\} \subseteq S$; J is an ideal of S. Consider $S/J = \{J, 1 + J, 2 + J, 3 + J, (0, 1) + J, (0, 2) + J, (0, 3) + J, (1, 0) + J, (2, 0) + J, (3, 0) + J, (1, 2) + J, (2, 1) + J, (1, 3) + J (3, 1) + J, (2, 3) + J, (3, 2) + J\}$ is a ring which is commutative of order sixteen and has zero divisors.

Clearly S/J is isomorphic with $\{N_o(Z_4), +, \times\}$ is isomorphic with $\{0, 1, 2, 3, (0, 1), (0, 2), (0, 3), (1, 0), (1, 2), (2, 1), (2, 0) (3, 0) (1, 3) (3, 1), (2, 3), (3, 1)\}$. Thus we can say as in case of Z we have

$$\frac{Z}{\langle nZ \rangle} \cong Z_n$$

like wise

$$\frac{\{N_c(Z), +, \times\}}{\{N_o(nZ), +, \times\}} \cong \{N(Z_n)\}, +, \times\}.$$

This is true even if the open intervals are replaced by closed intervals or half open-closed of half or closed-open intervals.

***Example 3.47***: Let $S = \{N_{oc}(Z), +, \times\}$ be a ring. Consider $J = \{N_{oc}(11Z), +, \times\} \subseteq S$ be an ideal of S. $S/J \cong \{N_{oc}(Z_{11}), +, \times\}$.

Thus we have the following theorem the proof of which is left for the reader.

**THEOREM 3.10:** *Let $S = \{N_{oc}(Z), +, \times\}$ be a ring. $J = \{N_{oc}(nZ), +, \times\} \subseteq S$ ($n < \infty$) is an ideal of S. Clearly S/J is the quotient ring and S/J is isomorphic to $\{N_{oc}(Z_n), +, \times\} = T$ for all $1 < n < \infty$. Thus we see the ring $S = \{N_{oc}(Z), +, \times\}$ is a commutative* *ring with unit S has zero divisors and elements in S are torsion free and S is of infinite order.*

We have set of elements which are not units for $([0, a])^n \neq [0, a]$ for any n.

Now having seen the concept of rings using special intervals we proceed onto built polynomials and matrices using these special intervals in the following chapter.

**Chapter Four**

# MATRIX THEORY USING SPECIAL CLASS OF INTERVALS

In this chapter we build algebraic structures using the matrix of natural class of intervals and describe a few of their associated properties.

**DEFINITION 4.1**: *Let $A = \{(a_1, \ldots, a_n)$ such that $a_i \in N_c(Z)$ (or $N_o(Z)$ or $N_{oc}(Z)$ or $N_{co}(Z)\}$ be a set of row interval matrix. A is known as the collection of natural closed intervals from Z. (Z can be replaced by Q or R or $Z_n$, $n < \infty$).*

We will illustrate this situation by some examples.

*Example 4.1*: Let $A = \{(a_1, a_2, a_3) \mid a_i \in N_{oc}(Z) \; 1 \leq i \leq 3\}$ be the collection of natural open closed intervals. Thus $x = ((0, 5], (7, 2], (3, 5]) \in A$ and $y = ((-7, 2], (8, 10], (-9, 0]) \in A$.
 $([0, 7], [9, 2], [7, 3], [8, 4], [-2, -7], [10, 8])$ is a closed natural $1 \times 6$ row interval.

**DEFINITION 4.2**: *Let* $A = \begin{bmatrix} a_1 \\ a_2 \\ \vdots \\ a_n \end{bmatrix}$ *where* $a_i \in N_o(Z)$; $1 \leq i \leq n$. *We call A the column natural open interval vector with entries from Z.* ($a_i \in N_c(Z)$ *or* $a_i \in N_{oc}(Z)$ *or* $a_i \in N_{co}(Z)$ $1 \leq i \leq n$), *Z can also be replaced by Q or* $Z_n$ *or R* ($n < \infty$).

*Example 4.2:* Let

$$B = \begin{bmatrix} [0,8) \\ [9,2) \\ [7,-1) \\ [8,0) \\ [-11,-14) \\ [0,16) \\ [18,3) \end{bmatrix}$$

is a $7 \times 1$ natural column closed-open interval from $N_{co}(R)$.

*Example 4.3:* Let

$$D = \begin{bmatrix} [8,1] \\ [-2,0] \\ \left[\dfrac{9}{7},0\right] \\ [-3,-5] \\ [20,\sqrt{3}] \\ [8,\sqrt{7}] \end{bmatrix}$$

is a $6 \times 1$ column natural closed interval $N_c(R)$.

Now having seen the column and row natural intervals we now proceed onto define natural interval matrix.

**DEFINITION 4.3**: *Let* $A = (a_{ij})$; $a_{ij} \in N_c(Z)$; $1 \leq i \leq n$, $1 \leq j \leq m$ *be a n × m natural class of closed intervals (or from* $N_o(Z)$ *or* $N_{oc}(Z)$ or $N_{co}(Z)$). A is defined as the $n \times m$ natural closed interval matrix.

We will illustrate this by some examples.

**Example 4.4**: Let

$$A = \begin{pmatrix} [0,5] & [3,0] & [9,2] \\ [-1,-3] & [-1,5] & [6,3] \\ [7,-2] & [0,7] & [3,1] \end{pmatrix}$$

be a $3 \times 3$ natural closed interval matrix.

**Example 4.5**: Let

$$P = \begin{pmatrix} (0,3) & (2,6) \\ (0,-2) & (6,1) \\ (7,1) & (-7,-10) \\ (8,9) & (0,-5) \\ (10,1) & (6,-4) \\ (-1,2) & (7,0) \\ (5,-3) & (8,2) \end{pmatrix}$$

be a $7 \times 2$ natural open interval matrix where entries of P are from $N_o(Z)$.

**Example 4.6**: Let

$$T = \begin{pmatrix} [0,3) & [5,2) & [3,0) \\ [7,1) & [3,4) & [0,-1) \\ [5,-2) & [9,1) & [0,8) \\ [9,10) & [1,-1) & [-8,-10) \end{pmatrix}$$

be a $4 \times 3$ natural closed open interval matrix with entries from $N_{co}(Z)$.

Now we can define operations addition / multiplication whenever there is compatibility of them.

To this end we show by examples as usual intervals are replaced by natural intervals that is intervals which can be decreasing (a, b); a > b, increasing (c, d) c < d and degenerate (a, b) a = b.

Let A = {([$a_1$, $b_1$), [$a_2$, $b_2$), [$a_3$, $b_3$) [$a_4$, $b_4$)) | $a_i$ $b_i$ ∈ Z, or [$a_i$, $b_i$) ∈ $N_{co}$ (Z), 1 ≤ i ≤ 4} be the collection of half closed open natural row intervals we add two half closed open row intervals componentwise in A; i.e. if x = ([3, 1), [0, -2), [7, 3), [5, 5)) and y = (2, [3, -1), [-2, 5) [-7, 1)) are in A, x + y = ([3, 1) + [2, 2), [0, -2) + [3, -1), [7, 3) + [-2, 5) + [5, 5) + [-7, 1)) = ([5, 3) [3, -3) [5, 8) [-2, 6)) ∈ A. Thus (A, +) is a semigroup under addition (in fact A is group under addition) called the semigroup of natural closed open row interval matrices.

*Example 4.7*: Let S = {([$a_1$, $b_1$], [$a_2$, $b_2$], …, [$a_{20}$, $b_{20}$) | $b_i$, $a_i$ ∈ Q; 1 ≤ i ≤ 20} be the natural row matrix closed intervals semigroup under addition.

*Example 4.8*: Let

$$P = \left\{ \begin{bmatrix} [a_1, b_1] \\ [a_2, b_2] \\ \vdots \\ [a_7, b_7] \end{bmatrix} \right\}$$

$b_i$, $a_i$ ∈ $Z_7$, 1 ≤ i ≤ 7; [$a_i$, $b_i$] ∈ $N_c$ ($Z_7$)} be a natural column closed interval semigroup under addition.

*Example 4.9*: Let

$$W = \left\{ \begin{bmatrix} [a_1, b_1] \\ [a_2, b_2] \\ \vdots \\ [a_{12}, b_{12}] \end{bmatrix} \right\}$$

[$a_i$, $b_i$] ∈ $N_c$ (R)} be the semigroup under addition.

*Example 4.10*: Let

$$P = \left\{ \begin{bmatrix} [a_1, b_1] \\ [a_2, b_2] \\ \vdots \\ [a_{11}, b_{11}] \end{bmatrix} \right\}$$

$[a_i, b_i] \in N_c(Z_{45})\}$ be a finite commutative semigroup under addition.

*Example 4.11*: Consider $S = \{([a_i, b_i]) \mid$ be a $3 \times 2$ closed natural interval matrix from the collection $N_c(Z_4)\}$. S is a semigroup under addition. If $A = ([a_i, b_i])$ and $B = ([c_i, d_i])$ are in S, $A + B = [(a_i, b_i)] + ([c_i, d_i]) = ([a_i + c_i \pmod 4), (b_i + d_i) \bmod 4])$ is in S. Further $A + B = B + A$, so S is a commutative semigroup of finite order. We get semigroup if $N_c(Z_4)$ is replaced by $N_o(Z_4)$ or $N_{oc}(Z_4)$ or $N_{co}(Z_4)$. All semigroups got are commutative and are of same order.

One can use them according to the situations and apply them. Thus we can get both finite and infinite semigroups built using $N_c(Z)$ or $N_c(Z_n)$ or $N_c(R)$ or $N_c(Q)$ (closed intervals can be replaced by open intervals, half open-closed intervals and half closed open intervals).

Thus using these natural class of intervals we have more choices to our solutions and so that we can choose the appropriate one. We see one can define product in case of natural class of intervals of matrices $m \times n$ when $m = 1$ and $n$ any positive integer or when $m = n$ a positive integer.

Now we give the class of natural semigroup matrices.

**DEFINITION 4.4**: *Let $S = \{([a_1, b_1], ..., [a_n, b_n]) \mid a_i, b_i \in N_c(Z); 1 \leq i \leq n\}$ be a natural class of row interval matrices. We define product on S as follows.*
*If $A = ([a_1, b_1], ..., [a_n, b_n])$ and $B = ([c_1, d_1], ..., [c_n, d_n])$ are in S. Then $A \times B = ([a_1, b_1], ..., [a_n, b_n]) \times ([c_1, d_1], ..., [c_n, d_n]) = ([a_1, b_1] \times [c_1, d_1]), ..., ([a_n, b_n] \times [c_n, d_n]) = ([a_1 c_1, b_1 d_1], ..., [a_n c_n, b_n d_n])$ is in S. $(S, \times)$ is a semigroup of natural class of closed intervals matrices. $N_c(Z)$ can be replaced by $N_{oc}(Z)$ or* $N_{co}(Z)$ or $N_o(Z)$ and the definition continues to be true. Also if Z is replaced by Q or R or $Z_n$ the definition is true.

We will now describe this by a few examples.

**Example 4.12**: Let $S = \{((a_1, b_1), \ldots, (a_4, b_4)) \mid (a_i, b_i) \in N_{co}(Q); 1 \leq i \leq 4\}$ be a semigroup of row interval matrices under multiplication.

If $A = ([0, 3), [7, 2), [5, 1), 4)$ and $B = (7, [-2, 1), [0, -7), [-4, 2))$ are in S. Then $AB = ([0, 3) \times [7, 7), [7, 2) \times [-2, 1), [5, 1) \times [0, -7), [4, 4) \times [4, 2)) = ([0, 21), [-14, 2), [0, -7), [-16, 8)) \in S$. Clearly $(S, \times)$ is a commutative semigroup of infinite order.

**Example 4.13**: Let $W = \{((a_1, b_1), (a_2, b_2), \ldots, (a_7, b_7)) \mid (a_i, b_i) \in N_o(Z_{12}); 1 \leq i \leq 7\}$ be a semigroup of row natural open interval matrices under multiplication '×'. $(W, \times)$ is a semigroup of finite order and is commutative.

**Example 4.14**: Let $T = \{([a_1, b_1), [a_2, b_2)) \mid [a_i, b_i) \in N_{co}(Z_7)\}$ be a semigroup under multiplication. T is a commutative semigroup of finite order.

Now in case of multiplicative semigroups with row interval matrices we can define zero divisors, units, idempotents and nilpotents.

Also we can define substructures like ideals and subsemigroup. All these can be done as a matter of routine and we expect the interested reader to do this exercise.

However we will illustrate all these situations by some examples.

**Example 4.15**: Let $S = \{([a_1, b_1), [a_2, b_2), \ldots, [a_6, b_6)) \mid [a_i, b_i) \in N_{co}(Z)\}; 1 \leq i \leq 6\}$ be the semigroup under multiplication. Consider $V = \{([0, a_1), [0, a_2), \ldots, [0, a_6)) \mid a_i \in Z ; 1 \leq i \leq 6\} \subseteq S$ V is a subsemigroup of S. In fact V is an ideal of S. Consider $W = \{([a_1, b_1), [a_2, b_2), \ldots, [a_6, b_6)) \mid [a_i, b_i) \in N_{co}(2Z)\}; 1 \leq i \leq 6\} \subseteq S$ W is a subsemigroup of S which is also an ideal of S. Take $M = \{([a_1, b_1), [a_2, b_2), \ldots, [a_6, b_6)) \mid [a_i, b_i] \in$ $N_{co}(nZ)\}$ ; $1 \le i \le \infty$, $1 \le i \le 6\} \subseteq S$ is a subsemigroup as well as ideal of S. Infact S has infinite number of ideals and subsemigroup.

*Example 4.16*: Let G = $\{([a_1, b_1), [a_2, b_2), [a_3, b_3)) \mid [a_i, b_i) \in N_{co}(Q)\}$ ; $1 \le i \le 3\}$ be the semigroup under component wise multiplication. Clearly G is an infinite semigroup which is commutative.

Take W = $\{([a_1, b_1), [a_2, b_2), [a_3, b_3)) \mid [a_i, b_i) \in N_{oc}(Z)\}$ ; $1 \le i \le 3\} \subseteq G$. W is only subsemigroup of G and is not an ideal of G. Thus all subsemigroup of a semigroup in general are not ideals of G.

Consider P = $\{([a_1, b_1), 0, 0) \mid [a_1, b_1) \in N_{co}(Q)\} \subseteq G$; P is an ideal of G. Thus the semigroup. If V = $\{([a_1, b_1), 0, 0)$ such that $[a_1, b_1) \in N_{co}(3Z)\} \subseteq G$, V is only a subsemigroup and is not an ideal of G.

*Example 4.17*: Let T = $\{([a_1, b_1], [a_2, b_2], \ldots, [a_9, b_9]) [a_i, b_i] \in N_c(Z_2)\}$ ; $1 \le i \le 9\}$ be the semigroup under multiplication. T is a commutative semigroup of finite order. Consider W = $\{([a_1, b_1], [a_2, b_2] 0, \ldots, 0) \mid [a_i, b_i] \in N(Z_2)$ ; $1 \le i \le 2\} \subseteq T$, W is a subsemigroup of T as well as an ideal of T. Thus has finitely many ideals and subsemigroups.

Now we proceed onto give examples of zero divisors units idempotents and S-zero divisors, S-units and S-idempotents in semigroups under multiplication.

*Example 4.18*: Let W = $\{([a_1, b_1], [a_2, b_2], [a_3, b_3]) \mid [a_i, b_i] \in N_c(Q)\}$ ; $1 \le i \le 3\}$ be a semigroup under multiplication. Clearly W is of infinite order and has both subsemigroups and ideals in it.

Consider H = $\{([a_1, b_1], 0, [a_2, b_2]) \mid [a_i, b_i] \in N_c(Q)\}$ ; $1 \le i \le 2\} \subseteq W$, H is a subsemigroup as well as an ideal of W.

*Example 4.19*: Let M= $\{([b_1, a_1], [b_2, a_2], \ldots, [b_{12}, a_{12}]) \mid [b_i, a_i] \in N_c(Z_{12})\}$ ; $1 \le i \le 12\}$ is a semigroup under multiplication $(0, 0, 0, \ldots, [b_{11}, a_{11}], [b_{12}, a_{12}])$ in M is a zero divisor. Also if X =

$\{(0\ 0\ 0, \ldots, [b_{11}, a_{11}]\ [b_{12}, a_{12}]\}$ where $b_{11}, a_{11}, b_{12}, a_{12} \in N_c$ $(Z_{12})\} \subseteq M$ is a zero divisor. $y = ([0,4]\ [0,4] \ldots [0,4])$ and $x = ([0,3]\ [0,3] \ldots [0,3])$ in M. $xy = yx = (0\ 0\ \ldots\ 0)$.

Thus M has several zero divisors. Consider $x = \{([0,4]\ [0,4] \ldots [0,4]) \in M$, then $x^2 = x$.

***Example 4.20:*** Let $S = \{([a_1, b_1], [a_2, b_2], [a_3, b_3], [a_4, b_4])) \mid [a_i, b_i] \in N_{co}(R); 1 \le i \le 4\}$ be a commutative semigroup of natural intervals under multiplication. S has zero divisors, ideals, subsemigroups which are not ideals but has no idempotents but has units.

Next we proceed onto analyse the notion of Smarandache semigroups built using these natural class of intervals.

We define a semigroup to be a Smarandache semigroup if it has a proper subset which is a group [ ].

We will illustrate this situation by some examples.

***Example 4.21:*** Let $W = \{([a_1, b_1], [a_2, b_2], [a_3, b_3], [a_4, b_4]) \mid [a_i, b_i] \in N_c(Z_{19}); 1 \le i \le 4\}$ be a semigroup. Consider $H = \{([a_1, b_1], [a_2, b_2], [a_3, b_3], [a_4, b_4]) \mid [a_i, b_i] \in N_c(1, 18) \subseteq N_c(Z_{19}); 1 \le i \le 4\} \subseteq W$; H is a group under multiplication. So W is a Smarandache semigroup.

***Example 4.22:*** Let $M = \{([a_1, b_1], [a_2, b_2]) \mid [a_i, b_i] \in N_c(Z_{20}); 1 \le i \le 2\}$ be a semigroup of finite order. Consider $P = \{([1, 19], [1, 19]), (1, 1), (19, 19), ([19, 1], [19,1]), (1, [1, 19]), (19, [1, 19]), ([1, 19], 1), ([1,19], 19), ([19, 1], 1), ([19, 1], 19), (19, [19, 1]), (1, [19, 1]), ([1, 19], [19, 1]), ([19, 1], [1, 19])\ (1,19), (19,1)\} \subseteq M$. P is a group with $(1,1) = 1$ as its multiplicative identity so M is a S-semigroup.

***Example 4.23:*** Let $G = \{((a_1, b_1], (a_2, b_2], (a_3, b_3], (a_4, b_4], (a_5,b_5], (a_6, b_6]) \mid (a_i, b_i] \in N_{oc}(240); 1 \le i \le 6\}$ be a semigroup of finite order. Consider $H = \{((a_1, b_1], (a_2, b_2], (a_3, b_3], \ldots, (a_6, b_6]) \mid (a_i, b_i] \in N_{oc}(1, 239); 1 \le i \le 6\} \subseteq G$ is a group of order $6^4$. Hence G is a S-semigroup.

In view of this we have the following theorem.

**THEOREM 4.1:** *Let $G = \{([a_1, b_1], \ldots, [a_m, b_m]) / [a_i, b_i] \in N_c(Z_n); (n < \infty); 1 \leq i \leq m\}$ be a semigroup. G is a Smarandache semigroup.*

We just give the hint of the proof.

Consider $H = \{([a_1, b_1], \ldots, [a_m, b_m]) \mid [a_i, b_i] \in N_c(1, n-1); 1 \leq i \leq m\} \subseteq G$; H is a group. This is true if $N_c(Z_n)$ is replaced by $N_o(Z_n)$ or $N_{oc}(Z_n)$ or $N_{co}(Z_n)$ ($n < \infty$).

Thus we have a large class of semigroups built using natural class of intervals which are Smarandache semigroup.

We can as in case of finite semigroups study all the related properties like S-Lagrange theorem, S-Sylow theorem, S-Cauchy element and so on. For more refer [17].

*Example 4.24:* Let $V = \{([a_1, b_1], \ldots, [a_5, b_5]) / [a_i, b_i] \in N_c(Z); 1 \leq i \leq 5\}$ be a semigroup.
Consider $T = \{([a_1, b_1], [a_2, b_2], \ldots, [a_5, b_5]) / [a_i, b_i] \in N_c(1, -1), 1 \leq i \leq 5\} \subseteq V$; T is a group. So V is a Smarandache group.

*Example 4.25:* Let $V = \{([a_1, b_1], [a_2, b_2], [a_3, b_3]) / [a_i, b_i] \in N_c(Q); 1 \leq i \leq 3\}$ be a semigroup. Take $W = \{([a_1, b_1], [a_2, b_2], [a_3, b_3]) / [a_i, b_i] \in N_c(1, -1); 1 \leq i \leq 3\} \subseteq V$. W is a group so V is a S-semigroup.

$$W = \{([1, 1] [1, 1] [1, 1]), ([-1, -1] [-1, -1] [-1, -1]),$$
$$([1, 1] [1, 1] [-1, -1]), ([-1, -1] [1, 1] [1, 1]),$$
$$([1, 1] [-1, -1] [1, 1]), ([-1, -1] [-1, -1] [1, 1]),$$
$$([-1, -1] [1, 1] [-1, -1]), ([1, 1] [-1, -1] [-1, -1]),$$
$$([-1, 1] [-1, 1] [-1, 1]), ([1, -1] [1, -1] [1, -1]), \ldots\}$$

we have W is a finite order.

In view of this we have the following theorem.

**THEOREM 4.2:** *Let $S = \{([a_1, b_1], ..., [a_n, b_n]) / [a_i, b_i] \in N_c(Z)$ (or $(N_c(Q)$ or $(N_c(R))$; $1 \leq i \leq n\}$ be a semigroup. S is a Smarandache semigroup.*

Proof is direct and hence left as an exercise to the reader. We can place $N_c(Z)$ by $N_o(Z)$, $N_{oc}(Z)$ or $N_{co}(Z)$ still the conclusion of the theorem holds good. Further if Z is replaced by R or Q then also the theorem is true.

Now having seen we have a class of S-semigroup we now proceed onto prove that there exists semigroups built using the class of natural intervals which are not Smarandache semigroups.

*Example 4.26:* Let $A = \{([a_1, b_2], [a_3, a_4], ..., [a_{n-1}, a_n]) / [a_i, a_{i+1}] \in N_c(2Z); 1 \leq i \leq n-1\}$ be a semigroup under componentwise multiplication. Clearly A is not a Smarandache semigroup.

*Example 4.27:* Let $P = \{([a_1, b_1], [a_2, b_2), [a_3, b_3)) / [a_i, b_i) \in N_{co}(5Z); 1 \leq i \leq 3\}$ be a semigroup. P is not a Smarandache semigroup under componentwise multiplication.

*Example 4.28:* Let $V = \{((a_1, b_1), (a_2, b_2), ..., (a_6, b_6)) / [a_i, b_i] \in N_o(7Z); 1 \leq i \leq 6\}$ be a semigroup under multiplication. Clearly V is not a Smarandache semigroup.

*Example 4.29:* Let $T = \{((a_1, b_1], (a_2, b_2], (a_3, b_3], (a_4, b_4]) / (a_i, b_i] \in N_{oc}(4Z); 1 \leq i \leq 4\}$ be a semigroup under multiplication. T is not a Smarandache semigroup.

In view of these examples we have the following theorem.

**THEOREM 4.3:** *Let $V = \{((a_1, b_1], (a_2, b_2], ..., (a_n, b_n]) / (a_i, b_i] \in N_{oc}(tZ); 1 \leq t < \infty; 1 \leq i \leq n\}$ be semigroups under multiplication (as t varies in $(1, \infty)$ we get an infinite class of semigroups. V is not a Smarandache semigroup.*

The above theorem is true even if $N_{oc}(tZ)$ is replaced by $N_o(tZ)$ or $N_c(tZ)$ or $N_{co}(tZ)$, $1 < t < \infty$.

Now having seen examples of both Smarandache semigroups and semigroups which are not Smarandache we now proceed onto study about natural class of column intervals.

In the first place on the set of natural class of column intervals we cannot define multiplication. However we can study them as they form a Smarandache semigroup under addition.

*Example 4.30:* Let

$$P = \left\{ \begin{bmatrix} [a_1, b_1] \\ [a_2, b_2] \\ \vdots \\ [a_{10}, b_{10}] \end{bmatrix} \right\},$$

where $[a_i, b_i] \in N_c(Z_7)$; $1 \le i \le 10\}$ be a semigroup. In fact P is a group under addition and has subsemigroups which are groups under addition.

So P is a S-semigroup.
Take

$$M = \left\{ \begin{bmatrix} [a_1, b_1] \\ [a_2, b_2] \\ 0 \\ \vdots \\ 0 \end{bmatrix} \right\}$$

$[a_1, b_1]$, $[a_2, b_2] \in N_c(Z_7)\} \subseteq P$; M is a group, hence P is a S-semigroup.

*Example 4.31:* Let

$$S = \left\{ \begin{bmatrix} [a_1, b_1) \\ [a_2, b_2) \\ [a_3, b_3) \\ [a_4, b_4) \\ [a_5, b_5) \\ [a_6, b_6) \end{bmatrix} \right.$$

$[a_i, b_i) \in N_{co}(Z); 1 \le i \le 6\}$ be a semigroup under '+'.

Clearly S is a S-semigroup.

We cannot define the operation of multiplication on column matrix intervals.

Likewise we cannot define product on m × n (m ≠ n) interval matrices. Only when m = n we can define product.

We will illustrate this situation by some examples.

Choose

$$A = \begin{bmatrix} [9,10] & [0,3] & [-2,1] \\ [1,2] & [0,0] & [2,-1] \\ [1,1] & [3,1] & [-1,-4] \end{bmatrix}$$

and

$$B = \begin{bmatrix} [0,3] & [0,0] & [-2,-1] \\ [1,-1] & [3,1] & [-3,-9] \\ [8,0] & [5,7] & [1,-5] \end{bmatrix}$$

be two 3 × 3 natural intervals. Now we define the product AB as follows:

$$A \times B = \begin{bmatrix} [9,10] & [0,3] & [-2,1] \\ [1,2] & [0,0] & [2,-1] \\ [1,1] & [3,1] & [-1,-4] \end{bmatrix} \times$$

$$\begin{bmatrix} [0,3] & [0,0] & [-2,-1] \\ [1,-1] & [3,1] & [-3,-9] \\ [8,0] & [5,7] & [1,-5] \end{bmatrix}$$

$$= \begin{bmatrix} [-16,27] & [-10,10] & [-20,-42] \\ [16,6] & [10,-7] & [0,3] \\ [-5,2] & [4,-27] & [-12,10] \end{bmatrix}.$$

*Example 4.32:* Let

$$X = \left\{ \begin{bmatrix} [a_1,b_1] & [a_3,b_3] \\ [a_2,b_2] & [a_4,b_4] \end{bmatrix} \;\middle|\; [a_i, b_i] \in N_c(Z_{12}); 1 \le i \le 4 \right\}$$

be a semigroup under matrix multiplication.

*Example 4.33:* Let W = {6 × 6 natural interval matrices with entries from $N_o(Z)$} be a semigroup under matrix multiplication.

*Example 4.34:* Let V = {8 × 8 interval matrices with entries from $N_{oc}(Q)$} be a semigroup under matrix multiplication.

These semigroups have zero divisors, units and idempotents. Also ideals can be defined in this case.

Further these semigroups are non commutative and have identity given by

$$I_n = \begin{bmatrix} [1,1] & 0 & \cdots & 0 \\ 0 & [1,1] & \cdots & 0 \\ \vdots & \vdots & & \vdots \\ 0 & 0 & \cdots & [1,1] \end{bmatrix}_{n \times n}.$$

Only [1, 1] = 1 on the main diagonal.

The interval collection can be open intervals or closed intervals or half-open-closed interval or half closed-open interval.

Thus we can have S-semigroups also in case of semigroups of n × n interval matrices.

When the semigroup is built using $N_c(Z_n)$ or $N_o(Z_n)$ or $N_{oc}(Z_n)$ or $N_{co}(Z_n)$, $n < \infty$; we see they are of finite order. However in all other cases they are of infinite order.

Substructures definition is a matter of routine and interested reader can study define and give examples of them.

**Chapter Five**

# POLYNOMIAL INTERVALS

In this chapter we build semigroups of polynomials with coefficients from the natural class of intervals and study the properties about them.

Let $p(x) = \sum_{i=0}^{\infty} a_i x^i$ where $a_i \in N_c(Z)$ (or $N_o(Z)$ or $N_{oc}(Z)$ or $N_{co}(Z)$), is defined as the polynomial with coefficients from the natural class of intervals in the variable x. We can add two polynomials by adding the coefficients.

These polynomials will be known as the natural interval polynomials. The polynomials can have coefficients also from $N_o(Z_n)$ or $N_{oc}(Z_n)$ or $N_{co}(Z_n)$ or $N_c(Z_n)$ or $N_o(R)$ or $N_o(Q)$ and so on.

We will illustrate by some examples how multiplication or addition is carried out.

Let $p(x) = (0, 5) x^8 + (-7, -9) x^5 + (8, 0) x^3 + (-3, 2) x^2 + (2, -4) x + (8, 3)$ and $q(x) = (6, 3) x^4 + (-3, 2)x^3 + (7, 1) x^2 + (-3, 2)$

be two natural open interval polynomials. Their sum $p(x) + q(x)$
$= (0, 5) x^8 + (-7, -9) x^5 + (6, 3) x^4 + (5, 2) x^3 + (4, 3) x^2 + (2, -4)x + (5, 5)$.

Now we just show how product can also be defined. $p(x) = [3, 4)x^5 + [2, -3)x^3 + [-1, 2) x^2 + [8, 1)$ and $q(x) = [2, -4) x^4 + [-3, 0) x^3 + [0, -7) x + [8, 0)$ be polynomials with coefficients from $N_{co}(Z)$.

$$p(x) \cdot q(x) = [6, -16)x^9 + [4, +12)x^7 + [-2, -8) x^6 +$$
$$[16, -4) x^4 + [-9, 0)x^8 + [-6, 0) x^6 + [3, 0) x^5 + [-24, 0)x^3 + [0, -28) x^6 + [0, 21)x^4 + [0, -14) x^3 + [0, -7) x + [24, 0) x^5 + [16, 0) x^3 + [-8, 0) x^2 + [64, 0)$$
$$= [6, -16) x^9 + [-9, 0)x^8 + [4, 12)x^7 + [-8, -36) x^6 + [27, 0)x^5 +$$
$$[16, 17) x^4 + [-8, -14) x^3 + [-8, 0)x^2 + [0, -7) x + [64, 0).$$

Polynomial semigroups can be defined with respect to addition or with respect to multiplication.

We will illustrate these situations by some examples.

**Example 5.1:** Let $S = \left\{ \sum_{i=0}^{\infty} [a_i, b_i]x^i \mid [a_i, b_i] \in N_c(Q) \right\}$ be a

semigroup under addition.
Note S is also a semigroup under multiplication.

**Example 5.2:** Let $T = \left\{ \sum_{i=0}^{8} (a_i, b_i]x^i \mid (a_i, b_i] \in N_{oc}(Z) \right\}$ be a

semigroup under addition. T is commutative infinite semigroup.
We see T is not a semigroup under multiplication.

**Example 5.3:** Let $W = \left\{ \sum_{i=0}^{20} (a_i, b_i)x^{20} \mid (a_i, b_i) \in N_o(Z_{12}) \right\}$ be

the semigroup under addition W is a finite commutative semigroup with identity.

***Example 5.4:*** Let $M = \left\{ \sum_{i=0}^{\infty} [a_i, b_i]x^i \mid [a_i, b_i] \in N_c(Z_6) \right\}$ be a semigroup under multiplication.
1. Prove M has zero divisors.
2. Find ideals in M.
3. Can M have subsemigroups which are not ideals?
4. Can M have idempotents?
5. Is M a S-semigroup?
6. Does M contain S-ideals?
7. Can M have S-semigroups which are not S-ideals?

By answering these questions the reader becomes familiar with properties of these semigroups.

***Example 5.5:*** Let $M = \left\{ \sum_{i=0}^{8} (a_i, b_i]x^i \mid (a_i, b_i] \in N_{oc}(Z_2) \right\}$ be a semigroup under addition.
1. Find the order of M.
2. Can M have S-ideals?
3. Is M a S-semigroup?
4. Find subsemigroups which are not ideals in M.
5. Can M have S-Cauchy elements?
6. Can M have atleast one subsemigroup H such that o(H) / o(M)?

By studying these questions the reader can in general understand the properties of polynomial semigroup built using intervals from modulo integers.

***Example 5.6:*** Let $P = \left\{ \sum_{i=0}^{\infty} [a_i, b_i]x^i \mid [a_i, b_i] \in N_c(Z_3) \right\}$ be a semigroup under multiplication.

Clearly T = {[1, 1], [1, 2], [2, 1], [2, 2]} $\subseteq$ P is a group so P is a S-semigroup.

In view of this example we give few more examples.

***Example 5.7:*** Let $S = \left\{ \sum_{i=0}^{\infty} [a_i, b_i) x^i \,\middle|\, [a_i, b_i) \in N_{co}(Z_8) \right\}$ be a semigroup under multiplication. Take $P = \{[1, 1), [7, 7), [7, 1), [1, 7)\} \subseteq S$; P is a group.
Hence S is a S-semigroup.

$$I = \left\{ \sum_{i=0}^{\infty} [a_i, b_i) x^i \,\middle|\, [a_i, b_i) \in N_{co}(0, 2, 4, 6) \subseteq N_{co}(Z_8) \right\} \subseteq S$$

is an ideal of S. Clearly I is not an S-ideal of S.

***Example 5.8:*** Let $V = \left\{ \sum_{i=0}^{\infty} (a_i, b_i] x^i \,\middle/\, (a_i, b_i] \in N_{oc}(Z_{12}) \right\}$ be a semigroup. V is a S-semigroup for take $G = \{(1, 1], (1, 11], (11, 1], (11, 11]\} \subseteq V$ is a group under multiplication.

Consider $T = \left\{ \sum_{i=0}^{\infty} (a_i, b_i] x^i \,\middle|\, (a_i, b_i] \in N_{oc}(0, 6) \subseteq N_{oc}(Z_{12}) \right\} \subseteq V$; T is an ideal of V which is not a S-ideal.

Consider $M = \left\{ \sum_{i=0}^{\infty} (a_i, b_i] x^{2i} \,\middle|\, (a_i, b_i] \in N_{oc}(0, 6) \subseteq N_{oc}(Z_{12}) \right\} \subseteq V$; M is a subsemigroup of V and not an ideal of V. However M is not a S-subsemigroup.

Consider $H = \left\{ \sum_{i=0}^{\infty} (a_i, b_i] x^{2i} \,\middle|\, (a_i, b_i] \in N_{oc}(Z_{12}) \right\} \subseteq V$; H is a S-subsemigroup of V and is not an ideal of V.

In view of these we have the following theorem.

**THEOREM 5.1:** *Let* $S = \left\{ \sum_{i=0}^{\infty} [a_i, b_i] x^i \,\middle|\, [a_i, b_i] \in N_c(Z_n); n < \infty \right\}$ *be a semigroup under multiplication*
  i) *S is a S-semigroup.*
  ii) *S has S-subsemigroups.*
  iii) *S has subsemigroups which are not S-subsemigroups.*

iv) S has ideals which are not S-ideals.
v) S has subsemigroups which are not S-ideals.

***Example 5.9***: Let $V = \left\{ \sum_{i=0}^{9} (a_i, b_i)x^i \,\Big|\, (a_i, b_i) \in N_o(Z_5); 0 \leq i \leq 9; x^{10} = 1 \right\}$ be a semigroup under multiplication. V has subsemigroups, V is commutative and is of finite order.

***Example 5.10***: Let $M = \left\{ \sum_{i=0}^{\infty} [a_i, b_i]x^i \,\Big|\, [a_i, b_i] \in N_c(Z_{80}) \right\}$ be a semigroup under product. Take $T = \left\{ \sum_{i=0}^{\infty} [a_i, b_i]x^i \text{ where } [a_i, b_i] \in N_c \{0, 2, \ldots, 78\} \subseteq N_c(Z_{80}) \right\} \subseteq M$; T is a subsemigroup as well as ideal of M. Consider $W = \left\{ \sum_{i=0}^{\infty} [a_i, b_i]x^i \,\Big|\, [a_i, b_i] \in N_c(1, 79) \right\} \subseteq M$, W is only a subsemigroup and is not an ideal of M.

***Example 5.11***: Let $T = \left\{ \sum_{i=0}^{\infty} [a_i, b_i)x^i \text{ where } [a_i, b_i) \in N_{co}(Z_{120}) \right\}$ be a semigroup under product. $M = \left\{ \sum_{i=0}^{\infty} [a_i, b_i)x^i \,\Big|\, [a_i, b_i) \in N_{co}\{0, 2, \ldots, 118\} \subseteq N_{co}(Z_{120}) \right\} \subseteq T$ is a subsemigroup as well as ideal of T.

Take $V = \left\{ \sum_{i=0}^{\infty} [a_i, b_i)x^i \,\Big|\, [a_i, b_i) \in N_{co}\{1, 119\} \subseteq N_{co}(Z_{120}) \right\} \subseteq T$; V is a subsemigroup and not an ideal of T.

In view of this we have the following theorem.

**THEOREM 5.2**: *Let* $S = \left\{ \sum_{i=0}^{\infty} [a_i, b_i] x^i \mid [a_i, b_i] \in N_c(Z) \right\}$ *be a semigroup under multiplication.*

    i)   *S is a S-semigroup*
    ii)  *S has ideals*
    iii) *S has subsemigroups which are not ideals*
    iv) *S has zero divisors*
    v)  *S has no nontrivial idempotents.*

The proof is direct and hence left as an exercise to the reader.

*Example 5.12*: Let $S = \left\{ \sum_{i=0}^{\infty} [a_i, b_i) x^i \mid [a_i, b_i) \in N_{co}(Q) \right\}$ be a semigroup under multiplication. Consider P = {[1, 1), [-1, -1), [1, -1), [-1, 1)} ⊆ S. P is a group given by the following table under multiplication.

| × | [1, 1) | [-1, -1) | [1, -1) | [-1, 1) |
|---|---|---|---|---|
| [1, 1) | [1, 1] | [-1, -1) | [1, -1) | [-1, 1) |
| [-1, -1) | [-1, -1) | [1, 1) | [-1, 1) | [1, -1) |
| [1, -1) | [1, -1) | [-1, 1) | [1, 1) | [-1, -1) |
| [-1, 1) | [-1, 1) | [1, -1) | [-1, -1) | [1, 1) |

Hence S is a S-semigroup. Consider $H = \left\{ \sum_{i=0}^{\infty} [a_i, b_i) x^i \mid [a_i, b_i) \in N_{co}(Z) \right\} \subseteq S$. H is only a subsemigorup and is not an ideal of H. S has zero divisors for take $p(x) = \sum_{i=0}^{\infty} [0, a_i) x^i$ and $q(x) = \sum_{i=0}^{\infty} [a_i, 0) x^i$ in S. Clearly p(x) q(x) = [0, 0). Take $I = \left\{ \sum_{i=0}^{\infty} [0, a_i) x^i \mid [0, a_i) \in N_{co}(Q) \right\} \subseteq S$ is an ideal of S.

Clearly $I = \left\{ \sum_{i=0}^{\infty} [0, a_i) x^i \right\}$ is of infinite order where $[0, a_i) \in N_{co}$ (Q). Consider $W = \left\{ \sum_{i=0}^{\infty} [0, a_i) x^i \mid [0, a_i) \in N_{co} (Z) \subseteq N_{co} (Q) \right\}$ $\subseteq S$; W is only a subsemigroup and is not an ideal of S.

S has units for take $[3, \frac{1}{2})$ in S; $[1/3, 2)$ in S is such that $[1/3, 2), [3, \frac{1}{2}) = [1, 1) = 1$ the multiplicative identity of S.

Clearly $[1, 0)$ and $[0, 1)$ are nontrivial idempotents in S for $[0, 1) [0, 1) = [0, 1)$ and $[1, 0) [1, 0) = [1, 0)$.

***Example 5.13:*** Let $S = \left\{ \sum_{i=0}^{7} [a_i, b_i] x^i \mid [a_i, b_i] \in N_c (Z); x^8 = 1 \right\}$ be a semigroup under multiplication. S is a S-semigroup. S has ideals and subsemigroups. Further S has zero divisors.

***Example 5.14:*** Let $S = \left\{ \sum_{i=0}^{3} [a_i, b_i) x^i \mid [a_i, b_i) \in N_{co} (R) \text{ where } x^4 = 1 \right\}$ be a semigroup of infinite order, S has zero divisors. S is a S-semigroup and S has subrings which are not ideals of S.

***Example 5.15:*** Let $S = \left\{ \sum_{i=0}^{\infty} [a_i, b_i] x^i \mid [a_i, b_i] \in N_c (Z^+ \cup \{0\}) \right\}$ be a semigroup under multiplication. S has zero divisors and S is not a S-semigroup.

***Example 5.16:*** Let $T = \left\{ \sum_{i=0}^{5} [a_i, b_i) x^i \mid [a_i, b_i) \in N_{co} (R^+ \cup \{0\}) \text{ with } x^6 = 1 \right\}$ be a semigroup under multiplication, S is not a S-semigroup. $P = \left\{ \sum_{i=0}^{5} [a_i, 0) x^i \mid x^6 = 1 \right\} \subseteq T$ is an ideal of S. $A = W = \left\{ \sum_{i=0}^{5} [a_i, b_i) x^i \mid [a_i, b_i) \in N_{co} (Z^+ \cup \{0\}) \text{ with } x^5 = 1 \right\} \subseteq T$ is a subsemigroup of S which is not an ideal of S.

However S has zero divisors and units but no idempotents other than [1, 0) and [0, 1). Clearly S has no nilpotent elements in it. T is a commutative semigroup of infinite order.

Study of S-weakly Lagrange semigroup, S-Lagrange semigroup, S-p-Sylow semigroup etc. can be carried out in case of semigroup, $T = \left\{ \sum_{i=0}^{n} [a_i, b_i] x^i \mid n < \infty; [a_i, b_i] \in N_c (Z_m) \text{ with } x^n = 1; m < \infty \right\}$ as a matter of routine. For all these semigroups are of finite order which has ideals, subsemigroups and are S-subsemigroup with zero divisors and units. However idempotents and nilpotents exist if m is a composite number.

**Example 5.17:** Let $T = \left\{ \sum_{i=0}^{4} [a_i, b_i) x^i \mid x^5 = 1, [a_i, b_i) \in N_{co} (Z_4) \right\}$ be a semigroup. T is a S-semigroup for W = {[1, 1), [1, 3), [3, 1), [3, 3)} ⊆ T is a group under multiplication. Consider $p(x) = \sum_{i=0}^{4} [a_i, 0) x^i$ and $q(x) = \sum_{i=0}^{4} [0, a_i) x^i$ in T p(x) . q(x) = 0. So T has zero divisors. Take $G = \left\{ \sum_{i=0}^{4} [a_i, 0) x^i \mid [a_i, 0) \in N_{co} (Z_4); x^5 = 1 \right\} \subseteq T$. G is an ideal of G. Now $M = \left\{ \sum_{i=0}^{2} [a_i, b_i) x^{2i} \mid [a_i, b_i) \in N_{co} (Z_4) \right\} \subseteq T$ is only subsemigroup and is not an ideal of T. In fact M is a S-subsemigroup of T.

**Example 5.18:** Let $G = \left\{ \sum_{i=0}^{2} [a_i, b_i] x^i \mid [a_i, b_i] \in N_c (Z_2) \text{ with } x^3 = 1 \right\}$ be a semigroup of finite order.
G = {0, 1, [0, 1], [1, 0], x, [0, 1] x, [1, 0] x, $x^2$, [0, 1] $x^2$, [1, 0] $x^2$, 1 + x, 1 + [0, 1] x, 1 + [1, 0] x 1 + $x^2$, 1 + [0, 1]$x^2$, 1 + [1, 0] $x^2$ [1, 0] + x [1, 0] + [1, 0]x, [1, 0] + $x^2$ [1, 0] + [0, 1] x [1, 0] + [1, 0]$x^2$, [1, 0] + [0, 1]$x^2$, [0, 1] + x, [0, 1] + $x^2$, [0, 1] + [1, 0]x,

$[0, 1] + [1, 0] x^2$, $[0, 1] + [0, 1] x$, $[0, 1] + [1, 0]x^2$, $x + x^2$, $x + [0, 1]x^2$, $x + [1, 0]x^2$, $[1, 0] x + x^2$, $[1, 0] x + [1, 0] x^2$, $[1, 0] x + [0, 1] x^2$, $[0, 1] x + x^2$ $[0, 1] x + [0, 1] x^2$, $[0, 1] x + [1, 0]x^2$, $1 + [0, 1] x + x^2$, $1 + [0, 1] + [0, 1]x^2$, $1 + [0, 1] + [1, 0]x^2$, $1 + x + x^2$, $1 + x + [0, 1]x^2$, $1 + x + [1, 0]x^2$, $1 + [1, 0] x + x^2$, $1 + [1, 0] x + [1, 0]x^2$, $1 + [1, 0] x + [0, 1] x^2$, $[0, 1] + x + x^2$, $[0, 1] + [0, 1] x + [0, 1]x^2$, $[0, 1] + x + [0, 1]x^2$, $[0, 1] + [1, 0] x + [0, 1]x^2$, $[0, 1] + [0, 1]x + x^2$, $[0, 1] + [0, 1] x + [1, 0] x^2$, $[0, 1] + [1, 0] x + [1, 0] x^2$, $[0, 1] + [1, 0] x + x^2$, $[0, 1] + x + [1, 0] x^2$ $[1, 0] + x + x^2$, $[1, 0] + [1, 0] x + x^2$ $[1, 0] + x + [1, 0] x^2$, $[1, 0] + [1, 0] x + [1, 0] x^2$, $[1, 0] + [1, 0] x + [0, 1]x^2$, $[1, 0] + [0, 1] x + [1, 0] x^2$, $[1, 0] + [0, 1] x + [0, 1]x^2$, $[1, 0] + [0, 1] x + x^2$, $[1, 0] + x + [0, 1] x^2$ }

Clearly G has 64 elements in it. G has zero divisors.
G is a S-semigroup as $\{1, x, x^2\} = H$ is a group in G.
$T = \{0, 1, x, x^2\} \subseteq G$ is a S-subsemigroup.
We see under addition the interval polynomials are groups.
However under multiplication they do not enjoy the group structure.

**Example 5.19:** Let $G = \left\{ \sum_{i=0}^{\infty} [a_i, b_i) x^i \mid [a_i, b_i) \in N_{co}(Z) \right\}$ be a infinite abelian group under addition. Clearly $Z[x] = \left\{ \sum_{i=0}^{\infty} a_i x^i \mid a_i \in Z \right\} \subseteq G$ is a subgroup of G.

$Z^+[x] = \left\{ \sum_{i=0}^{\infty} a_i x^i \mid a_i \in Z^+ \right\}$ is a semigroup in G. Thus G is a Smarandache special definite group. We can study subgroups in these groups.

Next we proceed onto study matrix interval groups.

*Example 5.20:* Let P = {([$a_1$, $b_1$], [$a_2$, $b_2$], …, [$a_8$, $b_8$]) | [$a_i$, $b_i$] ∈ $N_c$ ($Z_7$)}, P is a group under component wise addition. Clearly P is commutative group of finite order. P has subgroups.

*Example 5.21:* Let M = {(($a_1$, $b_1$), ($a_2$, $b_2$) ($a_3$, $b_3$)) where ($a_i$, $b_i$) ∈ $N_o$ (Z); 1 ≤ i ≤ 3} be a group under addition of infinite order, M has many subgroups.

*Example 5.22:* Let M = {([$a_1$, $b_1$], [$a_2$, $b_2$], …, [$a_5$, $b_5$)) | [$a_i$, $b_i$] ∈ $N_{co}$ ($Z_{40}$) ; 1 ≤ i ≤ 5} be an abelian finite group. Clearly M has subgroups.

*Example 5.23:* Let T= {(($a_1$, $b_1$], ($a_2$, $b_2$], …, ($a_{10}$, $b_{10}$]) | ($a_i$, $b_i$] ∈ $N_{oc}$ (R) ; 1 ≤ i ≤ 10}. T is an additive abelian group of infinite order. T has infinitely many subgroups.

We need these concepts to built natural interval Smarandache vector spaces of a special type in later chapters.

**Chapter Six**

# NEW TYPES OF RINGS OF NATURAL INTERVALS

In this chapter we introduce the notion of rings and Smarandache rings using interval polynomials and interval matrices. We describe some essential properties and illustrate them by examples.

**DEFINITION 6.1**: *Let $S = \{([a_1, b_1], …, [a_n, b_n]) / [a_i, b_i] \in N_c(R) ; 1 \leq i \leq n\}$ be a commutative abelian group under addition and a commutative semigroup under multiplication. Thus $(S, +, \times)$ is ring with unit $(1, 1, …, 1)$ known as the natural class of matrix interval ring.*

*If we replace $N_c(R)$ by $N_o(R)$ or $N_{oc}(R)$ or $N_{co}(R)$ or $N_c(Z)$ or $N_c(Q)$ or $N_c(Z_n)$ $(n < \infty)$ or by closed open interval or open intervals or open closed intervals we see S is a ring.*

We will first illustrate this situation by some examples.

*Example 6.1:* Let $M = \{((a_1, b_1], (a_2, b_2], (a_3, b_3], (a_4, b_4]) \mid (a_i, b_i] \in N_{oc}(R); 1 \leq i \leq 4\}$ be a ring. This ring is commutative and is of infinite order.

*Example 6.2:* Let $R = \{([a_1, b_1), [a_2, b_2), \ldots, [a_8, b_8)) \mid [a_i, b_i) \in N_{co}(Z_5); 1 \leq i \leq 8\}$ be a commutative ring of finite order.

*Example 6.3:* Let $T = \{((a_1, b_1), (a_2, b_2), (a_3, b_3), (a_4, b_4), (a_5, b_5)) \mid (a_i, b_i) \in N_c(Z); 1 \leq i \leq 5\}$ be a ring of infinite order which is commutative with unit.

*Example 6.4:* Let $W = \{((a_1, b_1], (a_2, b_2], \ldots, (a_6, b_6]) \mid (a_i, b_i] \in N_{oc}(Z_{12}); 1 \leq i \leq 6\}$ be a finite commutative ring with unit.

The substructures like subrings and ideals can be taken as a work which is simple and direct.

We will illustrate these concepts by some examples.

*Example 6.5:* Let $S = \{([a_1, b_1], [a_2, b_2], \ldots, [a_6, b_6]) \mid [a_i, b_i] \in N_c(Z); 1 \leq i \leq 6\}$ be a ring. Consider $T = \{([a_1, b_1], [a_2, b_2], \ldots, [a_6, b_6]) \mid [a_i, b_i] \in N_c(3Z); 1 \leq i \leq 6\} \subseteq S$; T is a subring of S. T does not contain the identity. T is infact an ideal of S. Take $W = \{([a_1, b_1], [a_2, b_2], 0, 0, 0, 0) \mid [a_1, b_2], [a_2, b_2] \in N_c(Z)\} \subseteq S$ W too is a subring as well as ideal of S.

*Example 6.6:* Let $R = \{([a_1, b_1), [a_2, b_2), [a_3, b_3)) \mid [a_i, b_i) \in N_{co}(Z_9); 1 \leq i \leq 3\}$ be a ring. $T = \{(0, [a_1, b_1), 0) \mid [a_1, b_1) \in N_{co}(Z_9)\} \subseteq R$ is a subring as well as ideal of R.

$M = \{([0, a_1), [0, a_2), [0, a_3)) \mid a_i \in Z_9; 1 \leq i \leq 3\} \subseteq R$ is a subring as well as ideal of R.

Thus R has both ideals and subrings.

Consider $V = \{([a_1, b_1), [a_2, b_2), [a_3, b_3))$ where $[a_1, b_1) \in N_{co}\{0, 3, 6\} \subseteq N_{co}(Z_9)\} \subseteq R$ is an subring as well as a ideal of R.

*Example 6.7:* Let $R = \{([a_1, b_1], [a_2, b_2])$ where $[a_i, b_i] \in N_c(R); 1 \leq i \leq 2\}$ be a ring. Consider $H = \{([a_1, b_1], [a_2, b_2]) \mid [a_i, b_i] \in N_c(Z); 1 \leq i \leq 2\} \subseteq R$; H is a subring of R and is not an ideal of

R. However R has ideals. Consider T = {([$a_1$, $b_1$], 0) | [$a_1$, $b_1$] ∈ $N_c$ (R)} ⊆ R; T is a subring as well as an ideal of R.

Thus in these rings also every subring in general is not an ideal, however every ideal is a subring.

Now we can define the notion of zero divisors, units and idempotents in these rings which is a matter of routine and is left as an exercise to the reader.

We will illustrate them with examples.

***Example 6.8:*** Let R = {(([$a_1$, $b_1$), [$a_2$, $b_2$) , …, [$a_{10}$, $b_{10}$)) | [$a_i$, $b_i$) ∈ $N_{co}$ (Z) ; 1 ≤ i ≤ 10} be a ring. R has zero divisors for take x = ([0, $a_1$), …, [0, $a_{10}$)) and y = ([$a_1$, 0), ..., [$a_{10}$, 0)) in R.

Clearly x·y = 0. In fact R has several other zero divisors. However R has no units or idempotents except idempotents of the form x = (1, 0, …, 1) or y = (1 1 1, 0 … 0) or z = (0 0 1 1 1 0 0 1 1 0) and so on.

***Example 6.9:*** Let P = {(($a_1$, $b_1$), ($a_2$, $b_2$), ($a_3$, $b_3$)) | ($a_i$, $b_i$) ∈ $N_o$ (R); 1 ≤ i ≤ 3} be a ring. R has zero divisors and units but has no idempotents. x = (0, ($a_1$, $b_1$), 0) and y = (($a_1$, $b_1$), 0, 0) in P is such that x·y = (0, 0, 0). Take x = ((1/2, 8), (7/3, 4), (2, 1/8)) and y = ((2, 1/8), (3/7, ¼), (1/2, 8)) in P we see xy = (1 1 1) is the unit in P. Consider x = ((0, $a_1$), (0, $a_2$), ($a_3$, 0)) and y = (($a_1$, 0), ($a_2$, 0), (0, $a_3$)) in P. xy = (0, 0, 0).

***Example 6.10:*** Let P = {([$a_1$, $b_1$], [$a_2$, $b_2$] , …, [$a_{12}$, $b_{12}$]) | [$a_i$, $b_i$] ∈ $N_c$ (R) ; 1 ≤ i ≤ 12} be a ring. P is a S-ring as I = {(a, 0, …, 0) | a ∈ R ⊆ $N_c$(R)} is a field. However P has zero divisors, units but has no idempoents. P also has ideals which are S-ideals. P has subrings which are S-subrings.

***Example 6.11:*** Let M = {(($a_1$, $b_1$], ($a_2$, $b_2$] , …, ($a_6$, $b_6$]) | ($a_i$, $b_i$] ∈ $N_{oc}$ ($Z_7$) ; 1 ≤ i ≤ 6} be a ring. M is a S-ring. For T = {(0, a, 0, 0, 0, 0) | a ∈ $Z_7$ ⊆ $N_{oc}$ ($Z_7$)} ⊆ M is a field. Hence M is a S-ring, M has zero divisors and units but has no proper idempotents.

It has idempotents of the form (1 0 0 1 1 0), (0 1 1 0 1 1), (1 1 1 1 0 1) and so on. M has both S-ideals and S-subrings.

***Example 6.12:*** Let W = {(($a_1$, $b_1$), ($a_2$, $b_2$), ($a_3$, $b_3$)) | ($a_i$, $b_i$) ∈ $N_o$(Z)} be a ring. Clearly W is not a S-ring. But has zero divisors no units or idempotents. However W has ideals and subrings. W is of infinite order.

Several subrings of infinite order which are ideals can be obtained. Now the property of quotient rings can be studied as in case of rings $N_c$(Z) or $N_o$(R) and so on. The reader is left with this task.

Now we will give only one example before we proceed to study square matrix rings using natural class of intervals.

***Example 6.13:*** M = {([a, b], [c, d]) | [a, b], [c, d] ∈ $N_c$($Z_2$)} be a ring.

M = {([1, 0], [1, 0]), (1, 1), (0, 0), (0, 1), (1, 0), ([0, 1], 0), (0, [0, 1]), ([0, 1], [0, 1]), ([1, 0], 0), (0, [1, 0]), ([1, 0], [0, 1]), ([0, 1], [1, 0]), (1, [0, 1]), ([0, 1], 1), ([1, 0], 1), (1, [1, 0])}

is a ring of order 16.

M is a S-ring as T = {[0, 0], [1, 1]} ⊆ M is a field isomorphic to $Z_2$.

Consider V = {0, (0, [0, 1]), (0, 1), (0, [1, 0])} ⊆ M is an ideal of M. Consider
M/V = {V, (1, 1) + V, (1, 0) + V, ([1, 0], 0) + V, ([1, 0], [1, 0]) + V, ([1, 0], 1) + V, (1, [0, 1]) + V, ([1, 0], [0, 1]) + V, (1, [1, 0]) + V, ([0, 1], 1) + V, ([0, 1], [1, 0]) + V, ([0, 1], [0, 1]) + V, ([1, 0], [1, 0] + V}.

M/V has thirteen elements. Clearly M/V has no zero divisors. Thus M/V is a semifield with 13 elements and is of characteristic two.

So using this method we can get several semifields of different orders with varying characteristics. This answers an open problem [15, 18, 19] of existence of semifields of different characteristic other than zero. Now these quotient rings using the natural class of row matrices yields semifields.

Now we proceed onto recall that a n × n interval matrix with natural intervals defines a ring of both finite or infinite order.

We will illustrate this by some examples.

*Example 6.14:* Let

$$M = \left\{ \begin{bmatrix} [a_1b_1] & [a_3b_3] \\ [a_2b_2] & [a_4b_4] \end{bmatrix} \middle| [a_i, b_i] \in N_c(Z); 1 \le i \le 4 \right\},$$

M is a ring with interval matrix addition and multiplication. Clearly M is an infinite ring which is non commutative.

*Example 6.15:* Let V = {All 3 × 3 matrices with intervals from $N_{oc}(Z_{10})$}. V is a ring of finite order which is non commutative. V has zero divisors, idempotents and units;

$$\begin{pmatrix} (1,1] & 0 & 0 \\ 0 & (1,1] & 0 \\ 0 & 0 & (1,1] \end{pmatrix} = I_{3 \times 3}$$

is the identity element of V.

*Example 6.16:* Let T = {All 5 × 5 interval matrices with elements from $N_{co}(Q)$}, T is a ring infact a S-ring of infinite order and is non commutative.

*Example 6.17:* Let W = {All 10 × 10 interval matrices with elements form $N_o(Z_2)$}, W is a finite non commutative ring which is a S-ring.

**THEOREM 6.1:** *Let M = {All n × n interval matrices with intervals from $N_{oc}(Z_p)$ or $N_o(R)$ or $N_c(Q)$ or $N_c(Z_p)$ or $N_{co}(Z_p)$ or $N_o(Z_p)$ or $N_c(R)$ or $N_o(Q)$ and so on} be a ring. M is a S-ring (p- a prime).*

*Proof:* Obvious from the fact that

$$P = \left\{ \begin{bmatrix} a & 0 & \ldots & 0 \\ 0 & 0 & & 0 \\ \vdots & \vdots & & \vdots \\ 0 & 0 & \ldots & 0 \end{bmatrix} \right|$$

$a \in N_c(Z_p)$ or $N_o(Q)$ or $N(R)$, p a prime$\} \subseteq M$ is field of characteristic p, $Z_p$ or field Q or field R. Hence the claim.

However certain rings with interval entries from $N_c(Z_n)$ when $Z_n$ is a S-ring will be a S-ring.

**THEOREM 6.2**: *Let $V = \{$all $m \times m$ interval matrices with entries form $N_c(Z_n)$ or $N_o(Z_n)$ or $N_{oc}(Z_n)$ or $N_{co}(Z_n)$ where $Z_n$ is a S-ring$\}$, then V is a S-ring.*

We just give the hint of the proof.

Take $P = \{a_1, \ldots, a_r \,/\, a_i \in T \subseteq Z_n$ where T is a field in $Z_n$ and $1 \leq i \leq r < n-1\}$.

Consider

$$W = \left\{ \begin{pmatrix} a & 0 & \ldots & 0 \\ 0 & 0 & & 0 \\ \vdots & \vdots & & \vdots \\ 0 & 0 & \ldots & 0 \end{pmatrix} \right.$$

interval matrices with intervals from $N_c(P \subseteq Z_n)$ or $N_{co}(P \subseteq Z_n)$ or $N_{oc}(P \subseteq Z_n)$ or $N_o(P \subseteq Z_n)\} \subseteq V$; W is a field. Thus V is a S-ring.

We will illustrate this situation by some examples.

*Example 6.18*: Let V = {All $7 \times 7$ interval matrices with entries from $N_c(Z_3)$}, V is a ring. V is a S-ring for take

$$P = \left\{ \begin{pmatrix} 0 & 0 & \cdots & 0 \\ 0 & 0 & & 0 \\ \vdots & \vdots & & \vdots \\ 0 & 0 & \cdots & 0 \end{pmatrix}_{7 \times 7}, \begin{pmatrix} [1,1] & 0 & \cdots & 0 \\ 0 & 0 & & 0 \\ \vdots & \vdots & & \vdots \\ 0 & 0 & \cdots & 0 \end{pmatrix}_{7 \times 7} \right.$$

$$\left. \begin{pmatrix} [2,2] & 0 & \cdots & 0 \\ 0 & 0 & & 0 \\ \vdots & \vdots & & \vdots \\ 0 & 0 & \cdots & 0 \end{pmatrix}_{7 \times 7} \right\} \subseteq V,$$

P is a field. Hence V is a S-ring.

***Example 6.19:*** Let M = {All 3 × 3 interval matrices with entries from $N_o(Z_6)$} be a ring.
Consider

$$P = \left\{ \begin{pmatrix} [0,0] & 0 & 0 \\ 0 & 0 & 0 \\ 0 & 0 & 0 \end{pmatrix}, \begin{pmatrix} [3,3] & 0 & 0 \\ 0 & 0 & 0 \\ 0 & 0 & 0 \end{pmatrix} \right\} \subseteq M.$$

P is a field isomorphic with $Z_2$ as $3 + 3 \equiv 6 \equiv 0 \pmod{6}$ and $3 \cdot 3 \equiv 3 \pmod 6$. Thus M is a S-ring.

***Example 6.20:*** Let R = {all 5 × 5 interval matrices with intervals from $N_{oc}(Z_{10})$} be a ring.
Consider

$$T = \left\{ \begin{pmatrix} 0 & 0 & 0 & 0 & 0 \\ 0 & 0 & 0 & 0 & 0 \\ 0 & 0 & 0 & 0 & 0 \\ 0 & 0 & 0 & 0 & 0 \\ 0 & 0 & 0 & 0 & 0 \end{pmatrix}_{5 \times 5}, \begin{pmatrix} 5 & 0 & 0 & 0 & 0 \\ 0 & 0 & 0 & 0 & 0 \\ 0 & 0 & 0 & 0 & 0 \\ 0 & 0 & 0 & 0 & 0 \\ 0 & 0 & 0 & 0 & 0 \end{pmatrix} \right\} \subseteq R,$$

T is a field and hence R is a S-ring.

In view of this we have the following theorem.

**THEOREM 6.3**: *Let $T = \{$all $n \times n$ interval matrices with intervals from $N_c(Z_{2p})$ or $N_o(Z_{2p})$ or $N_{oc}(Z_{2p})$ or $N_{co}(Z_{2p})$, $p$ a prime$\}$ be a ring; T is a S-ring.*

*Hint:* Consider

$$M = \left\{\begin{pmatrix} 0 & 0 & \cdots & 0 \\ 0 & 0 & & 0 \\ \vdots & \vdots & & \vdots \\ 0 & 0 & \cdots & 0 \end{pmatrix}, \begin{pmatrix} p & 0 & \cdots & 0 \\ 0 & 0 & & 0 \\ \vdots & \vdots & & \vdots \\ 0 & 0 & \cdots & 0 \end{pmatrix}\right\} \subseteq T$$

is a field.
Hence the claim.

*Example 6.21:* Let

$$R = \left\{\begin{bmatrix} [a_1 b_1] & [a_2 b_2] \\ [a_3 b_3] & [a_4 b_4] \end{bmatrix}\right\}$$

where $[a_i, b_i] \in N_c(Z_{12})\}$ be a ring.
Consider

$$W = \left\{\begin{pmatrix} 0 & 0 \\ 0 & 0 \end{pmatrix}, \begin{pmatrix} 4 & 0 \\ 0 & 0 \end{pmatrix}, \begin{pmatrix} 8 & 0 \\ 0 & 0 \end{pmatrix}\right\} \subseteq R.$$

W is a field with 4 as its unit.
Hence R is a S-ring.

*Example 6.22:* Let $M = \{$all $3 \times 3$ interval matrices with entries from $N_c(Z_{30})\}$ be a ring.
Consider

$$P = \left\{ \begin{pmatrix} 0 & 0 & 0 \\ 0 & 0 & 0 \\ 0 & 0 & 0 \end{pmatrix}, \begin{pmatrix} 10 & 0 & 0 \\ 0 & 0 & 0 \\ 0 & 0 & 0 \end{pmatrix}, \begin{pmatrix} 20 & 0 & 0 \\ 0 & 0 & 0 \\ 0 & 0 & 0 \end{pmatrix} \right\} \subseteq M$$

is a field so M is a S-ring.

*Example 6.23:* Let M = {all 4 × 4 interval matrices with entries from $N_{oc}(Z_{40})$} be a ring.

$$T = \left\{ \begin{pmatrix} 0 & 0 & 0 & 0 \\ 0 & 0 & 0 & 0 \\ 0 & 0 & 0 & 0 \\ 0 & 0 & 0 & 0 \end{pmatrix}, \begin{pmatrix} 8 & 0 & 0 & 0 \\ 0 & 0 & 0 & 0 \\ 0 & 0 & 0 & 0 \\ 0 & 0 & 0 & 0 \end{pmatrix}, \begin{pmatrix} 16 & 0 & 0 & 0 \\ 0 & 0 & 0 & 0 \\ 0 & 0 & 0 & 0 \\ 0 & 0 & 0 & 0 \end{pmatrix}, \right.$$

$$\left. \begin{pmatrix} 24 & 0 & 0 & 0 \\ 0 & 0 & 0 & 0 \\ 0 & 0 & 0 & 0 \\ 0 & 0 & 0 & 0 \end{pmatrix}, \begin{pmatrix} 32 & 0 & 0 & 0 \\ 0 & 0 & 0 & 0 \\ 0 & 0 & 0 & 0 \\ 0 & 0 & 0 & 0 \end{pmatrix} \right\} \subseteq M;$$

is a field. So M is a S-ring.

*Example 6.24:* Let R = {All 10 × 10 interval matrices with entries from $N_c(Z_{60})$} be a ring.

Consider M = $\left\{ \begin{bmatrix} a & 0 & \ldots & 0 \\ 0 & 0 & \ldots & 0 \\ \vdots & \vdots & & \vdots \\ 0 & 0 & \ldots & 0 \end{bmatrix} \right.$ a = [a, a] ∈ {0, 12, 24, 36,

48} $\subseteq Z_{60}$} $\subseteq$ R; M is a field isomorphic to $Z_5$ with 36 acting as identity. So R is a S-ring.

It is left for the reader to discuss about zero divisors, units, idempotents and nilpotents in these ring. Also interested reader can study S-units, S-zero divisors and S-idempotents of these rings. Further these rings by using $Z_p$ we can get finite division which are not fields.

This task is also left for the reader to describe and illustrate with examples.

Now we proceed onto study polynomial rings using these natural interval coefficients.

**DEFINITION 6.2**: *Let $R = \left\{ \sum_{i=0}^{\infty} [a_i, b_i] x^i \mid [a_i, b_i] \in N_c(Z) \text{ (or } (N_c(R) \text{ or } N_c(Q) \text{ or } N_c(Z_n))\right\}$; R under polynomial addition and multiplication is a ring infact a commutative ring with unit of infinite order.*

We will give examples of them.

***Example 6.25:*** Let $R = \left\{ \sum_{i=0}^{\infty} (a_i, b_i] x^i \mid (a_i, b_i] \in N_{oc}(Z_2) \right\}$ be a ring. Clearly R is of infinite order and R is a S-ring.

***Example 6.26:*** Let $W = \left\{ \sum_{i=0}^{\infty} [a_i, b_i) x^i \mid [a_i, b_i) \in N_{co}(Z_{20}) \right\}$ be a ring, W is of infinite order and is commutative. W is a S-ring as $T = \left\{ \sum_{i=0}^{\infty} a_i x^i \mid a_i \in \{0, 4, 8, 12, 16\} \subseteq Z_{20} \right\} \subseteq W$ is a subring.

Infact $G = \{0, 4, 8, 12, 16\} \subseteq Z_{20}$ is a field isomorphic to $Z_5$. So W is a S-ring.

***Example 6.27:*** Let $M = \left\{ \sum_{i=0}^{\infty} [a_i, b_i) x^i \mid [a_i, b_i) \in N_{co}(Z) \right\}$ be a ring. M is not a S-ring. M has subrings and ideals.

*Example 6.28:* Let $R = \left\{ \sum_{i=0}^{\infty} [a_i, b_i] x^i \mid [a_i, b_i] \in N_c(Q) \right\}$ be a ring. R is a S-ring.

*Example 6.29:* Let $P = \left\{ \sum_{i=0}^{\infty} (a_i, b_i) x^i \mid (a_i, b_i) \in N_o(Z_{15}) \right\}$ is a ring.

Now we can find rings of finite order.

*Example 6.30:* Let $M = \left\{ \sum_{i=0}^{3} [a_i, b_i) x^i \mid [a_i, b_i) \in N_{co}(Z_3); 0 \leq i \leq 3; x^4 = 1 \right\}$ be a ring of finite order. M has subrings and M is a S-ring. M = {0, 1, 2, x, 2x, $2x^2$, $x^2$, $x^3$, $2x^3$, [0, 1)x, [0, 1)2x [0, 1)$x^2$, [0, 1)$2x^2$, [0, 1)$x^3$, [0, 1)$2x^3$, …, [0, 1) $x^3$ + [1, 2)$x^2$ + [2, 1)x + [1, 0)}.

Take $I = \left\{ \sum_{i=0}^{3} [0, a_i) x^i \mid a_i \in Z_3 \right\} \subseteq M$. I is an ideal of M. We can find M/I.

*Example 6.31:* Let $W = \left\{ \sum_{i=0}^{5} [a_i, b_i] x^i \mid [a_i, b_i] \in N_c(Q); x^6 = 1 \right\}$ be a ring; W is a S-ring. W has ideals and subrings.

Consider $T = \left\{ \sum_{i=0}^{5} [a_i, b_i] x^i \mid [a_i, b_i] \in N_c(Z); x^6 = 1 \right\} \subseteq W$; T is only a subring and is not an ideal of W. Infact T is not a S-subring of infinite order.

*Example 6.32:* Let $W = \left\{ \sum_{i=0}^{2} [a_i, b_i) x^i \mid [a_i, b_i) \in N_{co}(Z_2); x^3 = 1 \right\}$ be a ring, W is a S-ring. W has ideals and subrings. W has units and zero divisors.

*Example 6.33:* Let $R = \left\{ \sum_{i=0}^{8} (a_i, b_i) x^i \mid (a_i, b_i) \in N_o(Z_{10}); 0 \leq i \leq 8; x^9 = 1 \right\}$ be a ring, R is a S-ring. R has ideals and subrings. R is of finite order. R has idempotents, units and zero divisors.

*Example 6.34:* Let $N = \left\{ \sum_{i=0}^{5} (a_i, b_i] x^i \mid (a_i, b_i] \in N_{oc}(Z_8); x^6 = 1 \right\}$ be a ring. N has zero divisors units and idempotents. N is of finite order. N has subrings and ideals.

*Example 6.35:* Let $P = \left\{ \sum_{i=0}^{3} (a_i, b_i) x^i \mid (a_i, b_i) \in N_o(Z_9); 0 \leq i \leq 3; x^4 = 1 \right\}$ be a ring, P is finite order has nontrivial units and zero divisors.

Several properties associated with rings can be studied in this case of interval coefficient polynomial rings also.

Now we now proceed onto define Smarandache vector spaces of type I using the natural class of intervals.

Chapter Seven

# VECTOR SPACES USING NATURAL INTERVALS

In this chapter we introduce the notion of vector spaces using natural intervals and Smarandache vector spaces of type II and modified Smarandache vector space of type II which we call as quasi module Smarandache vector space of type II and illustrate them by examples.

Also we study the special properties associated with them.

**DEFINITION 7.1**: *Let V be an additive abelian group of natural class of intervals. F be a field. If V is a vector space over F we call V a natural class of interval vector space over F.*

We will first illustrate this situation by some examples.

*Example 7.1:* Let $V = \{N_c(Q)\}$ is an abelian group with respect to addition. Q be a field. V is a natural class of interval vector space over Q.

The basis for V is $\{[0, 1], [1, 0]\}$ over Q; for any $[a, b] = a [1, 0] + b [0, 1] = [a, b]$.

So dimension of V is two.

**Example 7.2**: Let $V = \{N_o(R)\}$ be a vector space over Q. V is of infinite dimension over Q.

Just like usual vector spaces V can be of dimension finite or infinite.

**Example 7.3**: Let $V = (N_{oc}(Z_7))$ be a vector space over $Z_7$, $\{(0, 1], (1, 0]\} \subseteq V$ is a basis of V over $Z_7$. Thus dimension of V over $Z_7$ is two.

**Example 7.4**: Let $M = (N_{co}(Z_{11}))$ be a vector space over $Z_{11}$ of dimension two over $Z_{11}$.

**Example 7.5**: Let $V = \{N_c(R)\}$ be a vector space over R.
Dimension of V over R is two given by $\{[0, 1], [1, 0]\}$.
If R is replaced by Q then dimension of V over Q is infinite.
We see if $V = \{(N_c(Q)\}$ be an abelian group V is not a vector space over R, V is a vector space only over Q.

**Example 7.6**: Let $V = (N_{oc}(R))$, V is a vector space over Q of infinite dimension.

**Example 7.7**: Let $V = \{([a_1, a_2], [b_1, b_2], [c_1, c_2]) \mid [a_1, a_2], [b_1, b_2]$ and $[c_1, c_2] \in N_c(Q)\}$ be a vector space over Q.
Now $B = \{([1, 0], 0, 0), ([0, 1], 0, 0), (0, [1, 0], 0), (0, [0, 1], 0), (0, 0, [0, 1]), (0, 0, [1, 0])\} \subseteq V$ is a basis of V over Q. Take

$$x = ([a_1, b_1], [a_2, b_2] [a_3, b_3])$$
$$= a_1 ([1, 0], 0, 0) + b_1 ([0, 1], 0, 0) + a_2 (0, [1, 0], 0) + b_2 (0, [0, 1], 0) + a_3 (0, 0, [1, 0]) + b_3 (0, 0, [0, 1])$$
$$= ([a_1, 0], 0, 0) + ([0, b_1], 0, 0) + (0, [a_2, 0], 0) + (0, [0, b_2], 0) + (0, 0, [a_3, 0]) + (0, 0, [0, b_3])$$
$$= ([a_1, b_1], 0, 0) + (0, [a_2, b_2], 0) + (0, 0, [a_3, b_3])$$
$$= ([a_1, b_1], [a_2, b_2] [a_3, b_3]).$$

Thus V is of dimension six over Q.

***Example 7.8:*** Let $V = \left\{ \begin{bmatrix} [a_1,b_1] \\ [a_2,b_2] \\ [a_3,b_3] \\ [a_4,b_4] \end{bmatrix} \middle| [a_i, b_i] \in N_c(Z_7) ; 1 \leq i \leq 4 \right\}$

be a vector space over $Z_7$. V has only finite number of elements. V is generated by eight elements given by

$$B = \left\{ \begin{bmatrix} [1,0] \\ 0 \\ 0 \\ 0 \end{bmatrix}, \begin{bmatrix} [0,1] \\ 0 \\ 0 \\ 0 \end{bmatrix}, \begin{bmatrix} 0 \\ [1,0] \\ 0 \\ 0 \end{bmatrix}, \begin{bmatrix} 0 \\ [0,1] \\ 0 \\ 0 \end{bmatrix}, \right.$$

$$\left. \begin{bmatrix} 0 \\ 0 \\ [1,0] \\ 0 \end{bmatrix}, \begin{bmatrix} 0 \\ 0 \\ [0,1] \\ 0 \end{bmatrix}, \begin{bmatrix} 0 \\ 0 \\ 0 \\ [1,0] \end{bmatrix}, \begin{bmatrix} 0 \\ 0 \\ 0 \\ [0,1] \end{bmatrix} \right\} \subseteq V.$$

***Example 7.9:*** Let $V = \{((a_1, a_2), (b_1, b_2), \ldots, (a_n, b_n)) \mid (a_i, b_i) \in N_o(Z_{11}); 1 \leq i \leq n\}$ be a vector space over $Z_{11}$. Dimension of V over $Z_{11}$ is 2n.

However V has only finitely many elements. The basis B of V over $Z_{11}$ is given by $B = \{((1, 0), 0, \ldots, 0), ((0, 1) 0, 0, \ldots, 0), (0, (1, 0), 0, \ldots, 0), (0, (0, 1), 0, \ldots, 0), \ldots, (0, 0, \ldots, 0, (1, 0)), (0, 0, \ldots, 0, (0, 1))\} \subseteq V$ is a basis having exactly 2n elements in it.

***Example 7.10:*** Let $V = \left\{ \begin{bmatrix} [a_1,b_1] \\ [a_2,b_2] \\ \vdots \\ [a_{12},b_{12}] \end{bmatrix} \middle| [a_i, b_i] \in N_c(R) ; 1 \leq i \leq 12 \right\}$ be a vector space over Q.

Dimension of V over Q is infinite.

**Example 7.11:** Let $V = \left\{ \begin{pmatrix} (a_1,b_1] \\ (a_2,b_2] \\ \vdots \\ (a_7,b_7] \end{pmatrix} \middle| (a_i, b_i] \in N_{oc}(R) ; 1 \leq i \leq 7 \right\}$ be a vector space over R.

$$B = \left\{ \begin{bmatrix} (1,0] \\ 0 \\ \vdots \\ 0 \end{bmatrix}, \begin{bmatrix} (0,1] \\ 0 \\ \vdots \\ 0 \end{bmatrix}, \begin{bmatrix} 0 \\ (1,0] \\ \vdots \\ 0 \end{bmatrix}, \begin{bmatrix} 0 \\ (0,1] \\ \vdots \\ 0 \end{bmatrix}, \ldots, \begin{bmatrix} 0 \\ \vdots \\ 0 \\ (1,0] \end{bmatrix}, \begin{bmatrix} 0 \\ \vdots \\ 0 \\ (0,1] \end{bmatrix} \right\} \subseteq V$$

is a basis of V and dimension of V over R is 14.

Thus as in case of usual vector spaces we see in case of natural class of intervals, the dimension depends on the field over which they are defined.

**Example 7.12:** Let

$$V = \left\{ \begin{bmatrix} [a_1,b_1] & [a_4,b_4] & [a_7,b_7] & [a_{10},b_{10}] \\ [a_2,b_2] & [a_5,b_5] & [a_8,b_8] & [a_{11},b_{11}] \\ [a_3,b_3] & [a_6,b_6] & [a_9,b_9] & [a_{12},b_{12}] \end{bmatrix} \right\}$$

$[a_i, b_i] \in N_c(Z_3) ; 1 \leq i \leq 12\}$ be a vector space over the field $Z_3$. Number of elements in V is finite. The dimension of V over $Z_3$ is 24.

$$B = \left\{ \begin{bmatrix} [1,0] & 0 & 0 & 0 \\ 0 & 0 & 0 & 0 \\ 0 & 0 & 0 & 0 \end{bmatrix}, \begin{bmatrix} [0,1] & 0 & 0 & 0 \\ 0 & 0 & 0 & 0 \\ 0 & 0 & 0 & 0 \end{bmatrix}, \ldots, \right.$$

$$\left. \begin{bmatrix} 0 & 0 & 0 & 0 \\ 0 & 0 & 0 & 0 \\ 0 & 0 & 0 & [1,0] \end{bmatrix} \begin{bmatrix} 0 & 0 & 0 & 0 \\ 0 & 0 & 0 & 0 \\ 0 & 0 & 0 & [1,0] \end{bmatrix} \right\} \subseteq V$$

is a basis of V over $Z_3$.

***Example 7.13:*** Let V = {3 × 3 interval matrices with entries from $N_{oc}(Z_5)$} be the vector space over the field $Z_5$. Clearly number of elements in V is finite. Dimension of V over $Z_5$ is $2 \times 3^2 = 18$. The basis

$$B = \left\{ \begin{bmatrix} (1,0] & 0 & 0 \\ 0 & 0 & 0 \\ 0 & 0 & 0 \end{bmatrix}, \begin{bmatrix} (0,1] & 0 & 0 \\ 0 & 0 & 0 \\ 0 & 0 & 0 \end{bmatrix}, \right.$$

$$\begin{bmatrix} 0 & (0,1] & 0 \\ 0 & 0 & 0 \\ 0 & 0 & 0 \end{bmatrix}, \begin{bmatrix} 0 & (1,0] & 0 \\ 0 & 0 & 0 \\ 0 & 0 & 0 \end{bmatrix}, \ldots,$$

$$\left. \begin{bmatrix} 0 & 0 & 0 \\ 0 & 0 & 0 \\ 0 & 0 & (0,1] \end{bmatrix}, \begin{bmatrix} 0 & 0 & 0 \\ 0 & 0 & 0 \\ 0 & 0 & (1,0] \end{bmatrix} \right\} \subseteq V$$

over $Z_5$.

***Example 7.14:*** Let V = {7 × 2 interval matrices with intervals from $N_{oc}(Q)$} be a vector space over field Q. V is of dimension 28 over Q.

***Example 7.15:*** Let V be a collection of 2 × 3 interval matrices with intervals from $N_c(R)$ over the field R. Dimension of V over R is 12. The basis B of V over R is given by

$$\left\{ \begin{bmatrix} [1,0] & 0 & 0 \\ 0 & 0 & 0 \end{bmatrix}, \begin{bmatrix} [0,1] & 0 & 0 \\ 0 & 0 & 0 \end{bmatrix}, \begin{bmatrix} 0 & [0,1] & 0 \\ 0 & 0 & 0 \end{bmatrix}, \begin{bmatrix} 0 & [1,0] & 0 \\ 0 & 0 & 0 \end{bmatrix}, \right.$$

$$\begin{bmatrix} 0 & 0 & [1,0] \\ 0 & 0 & 0 \end{bmatrix}, \begin{bmatrix} 0 & 0 & [0,1] \\ 0 & 0 & 0 \end{bmatrix}, \begin{bmatrix} 0 & 0 & 0 \\ [0,1] & 0 & 0 \end{bmatrix}, \begin{bmatrix} 0 & 0 & 0 \\ [1,0] & 0 & 0 \end{bmatrix},$$

$$\left. \begin{bmatrix} 0 & 0 & 0 \\ 0 & [1,0] & 0 \end{bmatrix}, \begin{bmatrix} 0 & 0 & 0 \\ 0 & [0,1] & 0 \end{bmatrix}, \begin{bmatrix} 0 & 0 & 0 \\ 0 & 0 & [1,0] \end{bmatrix}, \begin{bmatrix} 0 & 0 & 0 \\ 0 & 0 & [0,1] \end{bmatrix} \right\}$$

$\subseteq$ V is a basis of V over R with 12 elements in it.

However if R in example 7.15 is replaced by Q then dimension of V over Q is infinite. Thus V will become an infinite dimensional vector space over Q.

Now if in a vector space V a special type of product can be defined then V becomes a linear algebra.

***Example 7.16:*** Let V = {$N_c(R)$} be a vector space over the field R. V is a linear algebra over R.

***Example 7.17:*** Let V = {$N_o(Q)$} be a vector space over the field Q. V is a linear algebra over the field Q.

***Example 7.18:*** Let V = {$N_{oc}(Z_7)$} be a vector space over the field $Z_7$. V is a linear algebra over $Z_7$.

***Example 7.19:*** Let V = $\left\{ \begin{bmatrix} [a_1, b_1] \\ [a_2, b_2] \\ [a_3, b_3] \\ [a_4, b_4] \end{bmatrix} \middle| [a_i, b_i] \in N_c(R) ; 1 \le i \le 4 \right\}$

be a vector space over the field R. V is not a linear algebra over R.

In view of all these we have the following theorem.

**THEOREM 7.1**: *Let V be a linear algebra over a field F then V is a vector space over F. However if V is a vector space over F, V in general is not linear algebra over F.*

The later part is evident from example 7.19.

***Example 7.20***: Let $V = \{([a_1, b_1], [a_2, b_2), [a_3, b_3)) \mid [a_i, b_i) \in N_{co}(Z_{19}) ; 1 \leq i \leq 3\}$ be a vector space over the field $F = Z_{19}$. V is a linear algebra over $Z_{19}$.

***Example 7.21***: Let $V = \{3 \times 3$ interval matrices with intervals from $N_o(R)\}$ be a vector space over the field R. V is a linear algebra over R.

***Example 7.22***: Let V be a collection of $10 \times 5$ interval matrices with intervals from $N_{co}(Z_{43})$. V is only a vector space over the field $Z_{43}$ and is never a linear algebra over V.

***Example 7.23***: Let $V = \{3 \times 7$ interval matrices with intervals from $N_c(Z_{13})\}$ be a vector space over the field $Z_{13}$. V is not a linear algebra over $Z_{13}$.

***Example 7.24***: Let V be a collection of $n \times n$ interval matrices with intervals from $N_{oc}(Z_{29})$ over the field $Z_{29}$. V is a linear algebra over $Z_{29}$. V has a basis of $2 \times n^2$ elements and dimension of V is finite over $Z_{29}$. However the number of elements in V is also finite.

Now we can define as in case of vector space linear transformation of two interval vector spaces provided they are defined over the same field F.

***Example 7.25***: Let $V = N_c(Z_{11})$ and $W = N_{co}(Z_{11})$ be two vector spaces defined over the field $F = Z_{11}$. $T : V \to W$ defined by $T([a_i, b_i]) = [a_i, b_i)$ is a linear transformation of V to W.

***Example 7.26***: Let $V = N_o(R)$ and $W = \{([a_1, b_1], [a_2, b_2]) \mid [a_i, b_i] \in N_c(R); i = 1, 2\}$ be two vector spaces defined over the field R. Let $T : V \to W$ be a map such that $T((a_1, b_1)) = ([a_1, b_1], [a_1, b_1])$. T is a linear transformation of V to W.

***Example 7.27:*** Let

$$V = \left\{ \begin{bmatrix} [a_1,b_1] \\ [a_2,b_2] \\ \vdots \\ [a_6,b_6] \end{bmatrix} \right.$$

where $[a_i, b_i] \in N_c(Z_{11}); 1 \leq i \leq 6\}$ be a vector space over the field $F = Z_{11}$. Consider

$$W = \left\{ \begin{bmatrix} [a_1,b_1] & [a_4,b_4] \\ [a_2,b_2] & [a_5,b_5] \\ [a_3,b_3] & [a_6,b_6] \end{bmatrix} \right.$$

where $[a_i, b_i] \in N_c(Z_{11}); 1 \leq i \leq 6\}$ be a vector space over $Z_{11}$. Let $T : V \to W$ be a map such that

$$T\left( \begin{bmatrix} [a_1,b_1] \\ [a_2,b_2] \\ \vdots \\ [a_6,b_6] \end{bmatrix} \right) = \begin{bmatrix} [a_1,b_1] & [a_4,b_4] \\ [a_2,b_2] & [a_5,b_5] \\ [a_3,b_3] & [a_6,b_6] \end{bmatrix}.$$

T is a linear transformation of V to W.

***Example 7.28:*** Let $V = \{([a_1, b_1), [a_2, b_2) \ldots, [a_6, b_6)) \mid [a_i, b_i) \in N_{co}(Z_3); 1 \leq i \leq 6\}$ be a vector space over the field $Z_3$. $T : V \to V$ is the linear operator on V.

$T = ([a_1, b_1), [a_2, b_2), \ldots, [a_6, b_6)) = ([a_6, b_6), [a_5, b_5), \ldots, [a_2, b_2), [a_1, b_1))$ is a linear operator on V.

Now we can define as in case of linear transformation of usual vectors the kernel in case of linear transformation of natural class of interval vector spaces.

So if $T : V \to W$ is a linear transformation of the two vector spaces built using natural class of intervals over the same field F then kernel $T = \ker T = \{x \in V \mid T(x) = (0)\} \subseteq V$.

We will just give some examples first and then prove ker T is also a subspace of V.

**Example 7.29:** $V = \{((a_1, b_1), (a_2, b_2), (a_3, b_3)) \mid (a_i, b_i) \in N_o(Z_5); 1 \le i \le 3\}$ and $W = \{([a_1, b_1], [a_2, b_2]) \text{ where } [a_i, b_i] \in N_c(Z_5); 1 \le i \le 2\}$ be two vector spaces over the field $Z_5$.

Let $T : V \to W$

$T(((a_1, b_1), (a_2, b_2), (a_3, b_3))) = ([a_1, b_1], [a_3, b_3])$. Now kernel $T = \{T((a_1, b_1), (a_2, b_2), (a_3, b_3)) = ([0, 0], [0, 0]) = (0, 0)\}$

$= T(0, (a_2, b_2), 0) = (0, 0)$.

It is easily verified ker T is a proper subspace of V.

**Example 7.30:** Let $P = \left\{ \sum_{i=0}^{\infty} [a_i\ b_i]\ x^i \mid [a_i, b_i] \in N_c(Q) \right\}$ be a vector space over the field Q. Clearly P is of infinite dimension over Q.

Next we proceed onto define the notion of Smarandache double interval vector space of type II and doubly Smarandache interval vector space of type II. Even if we do not mention the term double from the context one can easily understand the situation.

Let V be an additive abelian group built using the natural class of intervals, F be a S-ring of natural class of intervals. If V is a vector space over the field P contained in the S-ring or V is a module over the S-ring then we call V to be a doubly Smarandache interval vector space of type II. The terms doubly and interval can be left out for from the context one can easily understand them.

We will illustrate this by examples.

**Example 7.31:** Let $V = \{N_c(Q) \times N_c(Q) \times N_c(Q)\} = \{([a_1, b_1], [a_2, b_2], [a_3, b_3]) \mid [a_i, b_i] \in N_c(Q); 1 \le i \le 3\}$ be a doubly Smarandache vector interval space over the S-ring $F = N_c(Q)$.

It is both interesting and important to note that V contains a normal or usual vector space over a field. Consider $P = (Q \times Q \times Q) = \{(a, b, c)\} \subseteq \{N_c(Q) \times N_c(Q) \times N_c(Q)\}$ is a vector space over Q of dimension 3. We call these vector spaces as inherited vector space of the S-vector space over the S-ring.

***Example 7.32***: Let $V = \left\{ \begin{bmatrix} [a_1,b_1) & [a_3,b_3) \\ [a_2,b_2) & [a_4,b_4) \end{bmatrix} \right.$ where $[a_i, b_i) \in$ $N_{co}(Z_{12})$, $1 \leq i \leq 4\}$ be a S-vector space over the S-ring $N_{co}(Z_{12})$. Take $P = \left\{ \begin{bmatrix} a & b \\ c & d \end{bmatrix} \right.$ where $a, b, c, d \in \{0, 4, 8\} \subseteq Z_{12}\} \subseteq V$ is a vector space over the field $F = \{0, 4, 8\} \subseteq N_{co}(Z_{12})$. Thus P is the inherited vector space of V.

***Example 7.33***: Let

$$P = \left\{ \begin{bmatrix} [a_1,b_1] & [a_5,b_6] \\ [a_2,b_2] & [a_6,b_6] \\ [a_3,b_3] & [a_7,b_7] \\ [a_4,b_4] & [a_8,b_8] \end{bmatrix} \right| [a_i, b_i] \in N_c(R), 1 \leq i \leq 8 \right\}$$

be a S-vector space over the S-ring $N_c(R)$.

$$V = \left\{ \begin{bmatrix} a_1 & a_5 \\ a_2 & a_6 \\ a_3 & a_7 \\ a_4 & a_8 \end{bmatrix} \right| a_i \in R; 1 \leq i \leq 8 \right\} \subseteq P$$

is an inherited vector space of P over the field $R \subseteq N_c(R)$. Also

$$W = \left\{ \begin{bmatrix} a_1 & a_5 \\ a_2 & a_6 \\ a_3 & a_7 \\ a_4 & a_8 \end{bmatrix} \right| a_i \in Q; 1 \leq i \leq 8 \right\} \subseteq P$$

is an inherited vector space of P over the field $Q \subseteq N_c(R)$.

In view of this we have the following theorem.

**THEOREM 7.2**: *Let V be a doubly Smarandache interval vector space over the interval S-ring S. Then V has a proper subset P, $P \subseteq V$; such that P is the inherited vector subspace of V over the field F; $F \subseteq S$.*

Proof is direct from the definition and properties of S-rings.

We will give some more examples before we proceed to define other properties.

*Example 7.34*: Let $V = \{((a_1, b_1), (a_2, b_2), (a_3, b_3), (a_4, b_4)) \mid (a_i, b_i) \in N_o(Z_{13}); 1 \leq i \leq 4\}$ be an interval Smarandache double vector space over the S-ring $N_o(Z_{13})$. Take $P = \{(a_1, a_2, a_3, a_4) \mid a_i \in Z_{13}; 1 \leq i \leq 4\} \subseteq V$ is the inherited vector space of V over the field $Z_{13}$.

*Example 7.35*: Let V = {all 5 × 5 interval matrices with entries from $N_{co}(Z_{43})$} be a double Smarandache interval vector space over the S-ring $N_{co}(Z_{43})$.

Consider M = {all 5 × 5 matrices with entries from $Z_{43}$} $\subseteq$ V; is the inherited vector space over the field $Z_{43} \subseteq N_{co}(Z_{43})$ of V.

**THEOREM 7.3**: *Let V be any double Smarandache interval vector space over $N_c(Z_p)$ (or $N_o(Z_p)$ or $N_{oc}(Z_p)$ or $N_{co}(Z_p)$; p a prime) over the S-ring $N_c(Z_p)$ (or $N_o(Z_p)$ so on respectively).*

*Then $M \subseteq V$ (where M is the inherited algebraic structure of V) is an inherited vector subspace over $Z_p$.*

The proof follows from the fact every set $N_c(Z_p)$ (or $N_o(Z_p)$ or $N_{co}(Z_p)$ or $N_{oc}(Z_p)$) contains $Z_p$. The same holds good for V hence the claim.

*Example 7.36*: Let $V = \{([a_1, b_1), [a_2, b_2), [a_3, b_3)) \mid [a_i, b_i) \in N_{co}(Z_{19}); 1 \leq i \leq 3\}$ be a interval vector space over the interval S-ring $N_{co}(Z_{19})$.

Consider $M = \{(a_1, a_2, a_3) \mid a_i \in Z_{19}; 1 \leq i \leq 3\} \subseteq V$; M is a vector subspace over the field $Z_{19}$.

**Example 7.37**: Let $V = \left\{ \begin{bmatrix} [a_1, b_1] \\ [a_2, b_2] \\ [a_3, b_3] \\ [a_4, b_4] \end{bmatrix} \middle| [a_i, b_i] \in N_c(Z_{15}), 1 \le i \le 4 \right\}$ be a double interval vector space over the S-ring $N_c(Z_{15})$. V is a S-vector space for take $P = \left\{ \begin{bmatrix} a_1 \\ a_2 \\ a_3 \\ a_4 \end{bmatrix} \middle| a_i \in \{0, 5, 10\} \subseteq Z_{15} \right\}$ is a vector subspace over the field $F = \{0, 5, 10\} \subseteq Z_{15}$.

Clearly V is not a linear algebra.

**Example 7.38**: Let

$$V = \left\{ \begin{pmatrix} (a_1, b_1) & (a_3, b_3) & (a_5, b_5) & (a_7, b_7) \\ (a_2, b_2) & (a_4, b_4) & (a_6, b_6) & (a_8, b_8) \end{pmatrix} \right\}$$

$1 \le i \le 8$, $(a_i, b_i) \in N_c(Z_{21})\}$ be a doubly interval vector space over the S-ring $N_c(Z_{21})$.

Consider

$$W = \left\{ \begin{pmatrix} a_1 & a_2 & a_3 & a_4 \\ a_5 & a_6 & a_7 & a_8 \end{pmatrix} \middle| a_i \in \{0, 7, 14\}, 1 \le i \le 8 \right\} \subseteq V$$

is a vector space over the field $S = \{0, 7, 14\} \subseteq Z_{21}$.

**Example 7.39**: Let

$$M = \left\{ \begin{pmatrix} [a_1, b_1) & [a_3, b_3) \\ [a_2, b_2) & [a_4, b_4) \end{pmatrix} \right\}$$

where $[a_i, b_i) \in N_{co}(Z_{35})$; $1 \le i \le 4\}$ be a doubly interval linear algebra over the S-ring $N_{co}(Z_{35})$. Consider $H = \left\{ \begin{pmatrix} a & b \\ c & d \end{pmatrix} \middle| a, b, c, d \in \{0, 7, 14, 21, 28\} \subseteq Z_{35} \right\} \subseteq M$ is a usual linear algebra

over the field $F = \{0, 7, 14, 21, 28\} \subseteq Z_{35}$. Thus M is a doubly S-interval linear algebra.

Now having seen examples of them we can define linear transformation provided they are defined over the same S-ring.

We will illustrate this situation as the definition is direct.

*Example 7.40*: Let

$$V = \left\{ \begin{pmatrix} [a_1,b_1] & [a_3,b_3] & [a_5,b_5] \\ [a_2,b_2] & [a_4,b_4] & [a_6,b_6] \end{pmatrix} \right.$$

$[a_i, b_i] \in N_c (Z_{15})$; $1 \le i \le 6\}$ be a double interval vector space over the S-ring $N_c(Z_{15})$.

Consider $W = \{([a_1, b_1], [a_2, b_2], \ldots, [a_6, b_6])$ where $[a_i, b_i] \in N_c (Z_{15})$; $1 \le i \le 6\}$ be a double interval vector space over $N_c(Z_{15})$.

Now $V_1 = \left\{ \begin{pmatrix} a_1 & a_2 & a_3 \\ a_4 & a_5 & a_6 \end{pmatrix} \right.$ $a_i \in \{0, 5, 10\} \subseteq Z_{15}\}$; $1 \le i \le 6\} \subseteq W$ is a vector space over the field, $F = \{0, 5, 10\} \subseteq Z_{15}$ and $W_1 = \{(a_1, a_2, \ldots, a_6) \mid a_i \in \{0, 5, 10\} \subseteq Z_{15}$; $1 \le i \le 6\}$ a vector space over the field $F = \{0, 5, 10\} \subseteq Z_{15}$. Define $T : V_1 \to W_1$ a linear transformation from $V_1$ to $W_1$. T is defined as the linear transformation of V to W.

Study in this direction is interesting and the reader is expected to develop and work in this direction. On same lines linear operator of interval S-linear algebras and S-vector spaces are defined.

We will however illustrate these by some examples.

*Example 7.41*: Let $V = \{6 \times 6$ natural interval matrices with intervals of the form $(a, b] \in N_{oc} (Z_{10})\}$ be a doubly interval vector space over the S-ring $N_{oc} (Z_{10})$.

Suppose $W = \{6 \times 6$ matrices with entries from $\{0, 5\} \subseteq Z_{10}\} \subseteq V$ be a vector subspace over the field, $F = \{0, 5\} \subseteq Z_{10}$. Define $T : V \to V$ such that $T (M) =$ upper triangular $6 \times 6$ matrices with entries from $\{0, 5\}$. T is a S-linear operator on V.

Thus we can for any double Smarandache vector space V over the S-ring define the notion of normal operator, diagonalizable normal operator and so on. To this end we will proceed onto define some properties in case of double Smarandache vector space.

We cannot define inner product space on usual interval vector spaces for double S-interval vector spaces as they are defined over S-rings.

We can easily with simple appropriate modifications derive the properties related with interval vector spaces and double Smarandache interval vector spaces.

We will illustrate the results in case of these two types of interval vector spaces.

Let $V = \{([a_1, b_1], [a_2, b_2]) \mid [a_i, b_i] \in N_c(Q)\}$ be a vector space over the field Q.

Now $\{([1, 0], 0), ([0, 1], 0), (0, [1, 0]), (0, [0, 1])\}$ is a basis of V over Q.

For
$([3, 2], [5, 7]) = 3([1, 0], 0) + 2([0, 1], 0) + 5(0, [1, 0]) + 7(0, [0, 1])$
$\phantom{([3, 2], [5, 7])} = ([3, 0], 0) + ([0, 2] + 0) + (0, [5, 0]) + (0, [0, 7])$
$\phantom{([3, 2], [5, 7])} = ([3, 2], 0) + (0, [5, 7]) = ([3, 2], [5, 7])$.

Thus dimension of V over Q is four.

We can as in case of usual vector spaces find the collection of linear transformations (linear operators) and study them. This task is simple when the interval vector spaces or linear algebras are defined over Q or R or $Z_p$ (p a prime). All results follow simply without any complications.

The little problem that arises when these interval linear algebras or interval vector spaces are over interval S-rings viz. $N_c(Q)$ or $N_c(R)$ or $N_c(Z_p)$ or $N_c(Z_n)$.

This problem is also over come very easily by considering the appropriate Smarandache vector spaces and obtaining the analogous results [15, 18-9].

Already polynomial with interval coefficients have been introduced so forming characteristic interval polynomial equations are a matter of routine. Also all techniques used in case of polynomials with real coefficients can be imitated.

Further every interval polynomial

$p(x) = [a_0, b_0] + [a_1, b_1]x + \ldots + [a_n, b_n] x^n$ with $[a_i, b_i] \in N_c(R)$ be a realized as $[a_0 + a_1 x + \ldots + a_n x^n, b_0 + b_1 x + \ldots + b_n x^n]$ where $a_i, b_i \in R$; $0 \le i \le n$.

So solving the usual equations and then putting back into intervals makes the problem easy as we have in our collection of $N_c(R)$ degenerate interval R, increasing interval [a, b]; a < b and decreasing interval [a, b], a > b. So all algebraic operations on R can be easily extended to $N_c(R)$ (or $N_o(R)$ or $N_{oc}(R)$ or $N_{co}(R)$) and solutions got in a very easy way.

Next to make the study of algebraic structures complete in the following chapter we introduce the notion of natural class of interval semirings and natural class of interval near ring.

**DEFINITION 7.2**: *Let S be the collection of natural intervals from $Z^+ \cup \{0\}$ or $R^+ \cup \{0\}$ or $Q^+ \cup \{0\}$. S under usual addition and multiplication is an interval semiring.*

*That is if $S = \{N_c (R^+ \cup \{0\})$ or $N_c (Q^+ \cup \{0\})$ or $N_c (Z^+ \cup \{0\})\}$ then S is an interval semiring or a natural class of interval semiring.*

*Even if closed intervals are replaced by open intervals or closed-open intervals or closed-open intervals S will continue to be a semiring.*

We will illustrate this situation by some examples.

**Example 7.42**: Let $S = \{N_{oc} (Z^+ \cup \{0\})\}$ be the natural class of open-closed intervals. S is a semiring $S = \{(a, b] \mid a, b \in Z^+ \cup \{0\}\}$.

**Example 7.43**: Let $P = \{(a, b) \mid a, b \in Q^+ \cup \{0\}\} = N_o (Q^+ \cup \{0\})$ be a natural class of interval semiring.

**Example 7.44**: Let $W = \{N_c (R^+ \cup \{0\})\}$ be the interval semirings.

All the interval semirings given in the above examples are of infinite order.

We see all these interval semirings are strict semirings for a + b = 0 in them imply a = 0 or b = 0.

Further all these interval semirings are of infinite order and commutative. Also these semirings are not semifields as they contain infinite number of interval zero divisors.

For if (0, 3) and (19, 0) are elements the interval semiring S = $\{N_o(Z^+ \cup \{0\})\}$ we see (0, 3) (19, 0) = (0, 0) = 0 is a zero divisor in S. Infact S has infinitely many zero divisors for if x = (0, n) and y = (m, 0) where m, n $\in Z^+ \cup \{0\}$ then x · y = (0, 0).

Thus S has infinitely many zero divisors.

But all these semirings are S-semirings for they contain interval semifields.

Consider S = $\{N_c (R^+ \cup \{0\})\}$ be an interval semiring. Take V = $\{[a, a] \mid a \in R^+ \cup \{0\}\} \subseteq$ S. V is an interval semifield where each interval is a degenerate interval and V = $R^+ \cup \{0\}$. Thus S contains a field so S is a S-semiring.

If we take R = $N_o (Z^+ \cup \{0\})$ to be an interval semiring, clearly R is a S-semiring as T = $\{(a, a) \mid a \in Z^+ \cup \{0\}\} \subseteq$ R is a degenerate interval semifield.

Similarly P = $\{N_{oc} (Q^+ \cup \{0\})\}$ is a S-semiring as W = $\{(a, a] \mid a \in Q^+ \cup \{0\}\} \subseteq$ P is a semifield of P. Thus we have the following theorem.

**THEOREM 7.**4: *Let S = $\{N_c (Z^+ \cup \{0\})$ (or $N_c (R^+ \cup \{0\}$ or $N_c (Q^+ \cup \{0\})\}$, S is a S-semiring*.

Proof is direct and hence is left as an exercise to the reader.

Further if the closed interval is replaced by open intervals or closed-open intervals or open-closed intervals the results continue to be true.

**Chapter Eight**

# ALGEBRAIC STRUCTURES USING FUZZY NATURAL CLASS OF INTERVALS

In this chapter we introduce two types of fuzzy algebraic structures using natural class fuzzy intervals and other using maps. We will discuss both types and give examples of them.

**DEFINITION 8.1**: *Let $N_c([0, 1]) = \{$collection of all closed fuzzy intervals $[a, b]$ where $a, b \in [0, 1]$ $(a < b, b < a$ or $a = b\}$ $(N_o([0, 1])$ will be the natural class of open fuzzy intervals, $N_{oc}([0, 1])$ fuzzy natural class of half open-closed intervals and $N_{co}([0, 1])$ the fuzzy natural class of half closed-open intervals).*

We will give examples of them.

*Example 8.1*: Let V = {[0.5, 0.2], 0.9, [0, 0.3], [0, 8, 0.1], [0.7, 0]} be the set of closed natural fuzzy intervals, we see some of them increasing, some decreasing and 0.9 = [0.9, 0.9] is a degenerate fuzzy interval to a fuzzy number.

*Example 8.2*: Let S = {(0.3, 0.1), (0, 0.7), (0.5, 0) (1, 0.3), (0.9, 1), 0.8, 1, (1 0)} be the class of natural fuzzy open intervals. Some of them are increasing, some decreasing and some degenerate.

Likewise we can have fuzzy half open-closed natural class of intervals and half closed-open natural class of intervals.

Consider P = {[0, 0.3), [1, 0.2), [0.8, 0.8), [0.2, 1), [0.9, 0.4), [0.6, 0.9)} P is the half closed-open natural class of fuzzy intervals.

T = {(0.8, 0.2], (0, 0.9], (1, 0.3], (0.2, 0.1], (0.5, 0.5], (0.3, 0.3]} is the half open-closed fuzzy intervals.

Now we proceed onto define operations on them. Let x = [0.1, 0.9] and y = [0.6, 0.2]. The product xy = [0.06, 0.18], this product is defined as the natural product and always under this product $N_c([0, 1])$ is a closed set.

Similarly $N_o([0, 1])$, $N_{co}([0, 1])$ and $N_{oc}([0, 1])$ are closed under product.

Now consider $N_o([0, 1])$ for any x = (0.3, 0.7) and y = (1, 0.4) in $N_o([0, 1])$ define min (x, y) = min {(0.3, 0.7), (1, 0.4)} = (min {0.3, 1}, min {0.7, 0.4}) = (0.3, 0.4).

We say $N_o([0, 1])$ under min operation is closed.

Now we can define max operation on $N_o([0, 1])$. For x = (0.7, 0.3) and y = (0.5, 0.8) be in $N_o([0, 1])$ max {x, y} = {max {0.7, 0.5}, max {0.3, 0.8}) = (0.7, 0.8) is in $N_o([0, 1])$.

Hence $N_o([0, 1])$ is also closed under max operation.

We will define now the natural fuzzy semigroups of intervals from $N_o([0, 1])$ or $N_{oc}([0, 1]$ or $N_{co}([0, 1])$ or $N_c([0, 1])$.

**DEFINITION 8.2**: *Let S = $N_o$ ([0, 1]) be the collection of natural fuzzy intervals. $N_o([0, 1])$ under the min operation is a commutative fuzzy natural interval semigroup of infinite order.*

***Example 8.3***: Let $S = \{N_{oc}([0, 1]), \min\}$ be a fuzzy interval semigroup of infinite order and commutative.

Infact '0' acts as the identity for min $\{(0, 0], (a, b]\} = (0, 0] = 0$.

Now if in the definition 8.2 min operation is replaced by max operation then $N_o([0, 1])$ is a commutative fuzzy semigroup with identity. Infact $1 = (1, 1)$ acts as the identity for max $\{(1, 1), (a, b)\} = (1, 1) = 1$.

Now if in the definition 8.2 min is replaced by the usual product × (multiplication) we get the fuzzy interval semigroup of infinite order which is commutative.

Clearly fuzzy semigroup under product is distinct different from the fuzzy semigroup under max or min operation for consider $x = [0.3, 0.7]$ and $y = [0.5, 0.2]$ in $N_c([0, 1])$. Then $x \cdot y = [0.3, 0.7] [0.5, 0.2] = [0.15, 0.14]$ but min $(x, y) = [0.3, 0.2]$ and max $(x, y) = [0.5, 0.7]$. Thus we see all the three operations are distinctly different.

Now having seen the examples of natural fuzzy semigroup we define fuzzy semigroups using the semigroups $N_o(R)$ or $N_c(Q)$ or $N_{co}(Z)$.

We will just give the definition and illustrate them by some examples.

**DEFINITION 8.3**: *Let $S = \{N_{oc}(R)\}$ be a semigroup under product the map $\eta : S \to [0, 1]$ with $\eta(x, y) \geq \min(\eta(x), \eta(y))$, $(S, \eta)$ is the fuzzy natural class of interval semigroups.*

Interested reader can give examples of them.

Thus we have two ways of defining fuzzy interval semigroups.

The development of these concepts can be considered as a matter of routine. For every algebraic structure V using intervals (V, η) the corresponding fuzzy structure can be defined.

We can also define the natural class of neutrosophic intervals. Just in $N_c(Z)$ or $N_c(R)$ or $N_c(Q)$ replace Z or R or Q by ZI or RI or QI then we obtain the natural class of pure neutrosophic intervals and all operations done in case of Z or Q or R can be verbatim carried out in case of ZI or QI or RI.

Now if we replace in $N_c(R)$ or $N_o(R)$ or $N_{oc}(R)$ or $N_{co}(R)$ R by $\langle Z \cup I \rangle$ or $\langle R \cup I \rangle$ or $\rangle Q \cup I \rangle$ ($\langle \ \rangle$ means generates the elements of Q and I or R and I or Z and I and any elements is of the form a + bI; a, b are reals and I the neutrosophic number such that $I^2 = I$). Then also we can have algebraic structures built using them [9-11].

# ALGEBRAIC STRUCTURES USING NEUTROSOPHIC INTERVALS

In this chapter we define different types of algebraic structures using the natural class of pure neutrosophic intervals and neutrosophic intervals [13].

$N_c(Z\,I) = \{[a, b] \mid a = xI \text{ and } b = yI;\ x, y \in Z\} = \{[xI, yI] \mid xI, yI \in ZI\}$. This will be known as pure neutrosophic integer closed intervals.

$N_o(QI) = \{(aI, bI) \mid aI, bI \in QI\}$ will be known as the natural class of pure neutrosophic rational open intervals. $N_{oc}(RI) = \{(aI, bI] \mid aI, bI \in RI\}$ will denote the pure neutrosophic open-closed intervals of reals.

$N_c(Z_nI) = \{[aI, bI] \mid aI, bI \in Z_nI\}$ denotes the natural class of pure neutrosophic closed modulo integer intervals.

Now $N_o(\langle Z \cup I \rangle) = \{(a + bI, c + dI) \mid a, b, c, d \in Z\}$ denotes the natural class of integer neutrosophic open intervals.

Likewise $N_{oc}(\langle Q \cup I \rangle)$ denotes the natural class of open closed neutrosophic rational intervals and $N_c(\langle Z_n \cup I \rangle)$ denotes the natural class of closed modulo integer intervals which are neutrosophic.

Now we proceed onto define algebraic structures using them.

**DEFINITION 9.1**: *Let $S = \{N_c (ZI), +\}$ be the natural class of neutrosophic interval semigroups of closed intervals under addition of infinite order.*

$N_c(ZI)$ can be replaced by $N_o(ZI)$, $N_{oc}(ZI)$ or $N_{co}(ZI)$. Further ZI can be replaced by QI or RI or $Z_nI$ and the definition holds good.

*Example 9.1*: Let $S = \{N_c (ZI), +\}$ be a semigroup of infinite order which pure neutrosophic interval semigroup.

*Example 9.2*: Let $M = \{N_{oc} (Z_5I)\}$ be the pure neutrosophic open-closed interval semigroup of finite order.

*Example 9.3*: Let $T = \{N_c (QI)\}$, T is a pure neutrosophic interval semigroup of infinite order under multiplication. Infact T has infinitely many zero divisors and units but has no neutrosophic idempotents except of the form [I, I], [I, 0], [0, I].

*Example 9.4*: Let $W = \{N_o (Z_{12} I)\}$ be the pure neutrosophic interval semigroup under multiplication of finite order.
    W has neutrosophic units, neutrosophic zero divisors, neutrosophic idempotents and neutrosophic nilpotents.
   [0, 4I] is such that $[0, 4I]^2 = [0, 4I]$ is an idempotents.
   [0,3] [0,4I] = 0 is a zero divisor. $[6I, 0]^2 = 0$ is a nilpotent element of W.
   $[0, 11I]^2 = [0, I]$ is only a semiunit.
   [1 II, I] × [1 II, I] = I is the unit [I, 11I] [1, 11] = I is the unit in W.
   Now we give examples of neutrosophic semigroups which are not pure neutrosophic.

***Example 9.5***: Let $M = \{N_c (<Z \cup I>)\}$ be the neutrosophic interval semigroup under multiplication. M is of infinite order and M has zero divisors.

Clearly $P = \{N_c (ZI)\} \subseteq M$, P is a subsemigroup of M which is pure neutrosophic. Also $T = \{N_c (Z)\} \subseteq M$ is a subsemigroup of M which is not neutrosophic. Infact M as many subsemigroups which are not neutrosophic.

Also M has many pure neutrosophic subsemigroups. Further $M_n = \{N_c (<nZ \cup I>)\} \subseteq M$ is a neutrosophic subsemigroup; for $n = 1, 2, \ldots$.

***Example 9.6***: Let $P = \{N_{oc} (<Q \cup I>)\}$ be the semigroup under addition. P has subsemigroups. P is a neutrosophic interval semigroup. $T = \{N_{oc} (Q)\} \subseteq P$ is an interval subsemigroup which is not neutrosophic.

Likewise $W = \{N_{oc} (QI)\} \subseteq P$ is a pure neutrosophic subsemigroup.

***Example 9.7***: Let $M = \{N_o (<Z_{15} \cup I>)\}$ be a neutrosophic semigroup under multiplication. M has subsemigroups which are pure neutrosophic neutrosophic as well as which are not neutrosophic.

Clearly M is of finite order.

Now we can define ideals in neutrosophic semigroups which is considered as a matter of routine. So we just give examples of them.

***Example 9.8***: Let $W = \{N_c (<Z \cup I>)\}$ be a neutrosophic interval semigroup under multiplication.

$T = \{N_c (ZI)\} \subseteq W$ is an ideal of W. However $V = N_c (Z) \subseteq W$ is only a subsemigroup of W and is not an ideal of W. Further $M = \{N_c (<nZ \cup I>)\} \subseteq W$ is a neutrosophic interval subsemigroup of W which is also an ideal of W.

However $N_c (nZ) \subseteq W$ are only neutrosophic subsemigroups which are not ideals ($2 \leq n < \infty$).

Now we just give an example of a Smarandache neutrosophic natural interval semigroup.

**Example 9.9**: Let $M = \{N_c (<Z \cup I>)\}$ be a semigroup under multiplication of neutrosophic intervals. Consider $P = \{[1, 1], [-1, -1] [1, -1], [-1, 1]\} \subseteq M$, M is a Smarandache neutrosophic interval semigroup as P is a group given by the following table.

| × | [1,1] | [-1,-1] | [1,-1] | [-1,1] |
|---|---|---|---|---|
| [1,1] | [1,1] | [-1,-1] | [1,-1] | [-1,1] |
| [-1,-1] | [-1,-1] | [1,1] | [-1,1] | [1,-1] |
| [1,-1] | [1,-1] | [-1,1] | [1,1] | [-1,-1] |
| [-1,1] | [-1,1] | [1,-1] | [-1,-1] | [1,1] |

**Example 9.10**: Let $T = \{N_o (<Q \cup I>)\}$ be a neutrosophic interval semigroup under multiplication. T is a S-neutrosophic interval semigroup.

All neutrosophic interval semigroups in general are not S-neutrosophic interval semigroup.

This is illustrated by examples.

**Example 9.11**: Consider $S = \{N_{oc} (3Z \cup I)\}$ be a neutrosophic interval semigroup under multiplication. S is not a S-semigroup for S has no proper subset which is a group under multiplication.

**Example 9.12**: Let $R = \{N_{co} (<5Z \cup I>)\}$ be a neutrosophic interval semigroup. Clearly R is not a S-semigroup.

**Example 9.13**: Let $M = \{N_c (<Z_5 \cup I>)\}$ be a semigroup of neutrosophic intervals under multiplication. $H = N_c (<1, 4, I>) \subseteq M$ is a group (neutrosophic) so M is a S-semigroup. Clearly M is of finite order.

We can define neutrosophic interval semigroup of polynomials and matrices as in case of interval semigroups. Here only neutrosophic intervals take the place of intervals.

We will only illustrate this situation by some examples.

**Example 9.14**: Let V = {3 × 3 neutrosophic intervals with entries from $N_c$ (ZI)} V is a neutrosophic interval semigroup under multiplication of infinite order. Clearly V is non commutative.

**Example 9.15**: Let M = {2 × 7 pure neutrosophic intervals with entries from $N_o$(QI)}. M is a pure neutrosophic semigroup of infinite order under addition.

**Example 9.16**: Let R = {5 × 3 pure neutrosophic interval matrices with entries from $N_c(Z_{12}I)$} be a pure neutrosophic semigroup. R is of finite order. As the operation in R is addition modulo 12, R is commutative.

**Example 9.17**: Let P = {2 × 2 interval pure neutrosophic matrices with entries from $N_o(Z_6I)$} under multiplication be a semigroup. P is non commutative and is of finite order. P has zero divisors units and idempotents in it.

**Example 9.18**: Let M = {8 ×8 interval neutrosophic matrices with entries from $N_c(<Z \cup I>)$} be a semigroup under multiplication. M has zero divisors and is of infinite order and M is a S-semigroup.

**Example 9.19**: Let S = {All 1 × 18 interval neutrosophic row matrices with entries from $N_{oc}(<Q \cup I>)$} be a semigroup under multiplication. S has zero divisors and units. S is a S-semigroup.

**Example 9.20**: Let T = $\left\{ \begin{bmatrix} a_1 \\ a_2 \\ a_3 \\ a_4 \\ a_5 \\ a_6 \end{bmatrix} \middle| a_i \in N_c (<Q \cup I>); 1 \leq i \leq 6 \right\}$ be

an additive semigroup of neutrosophic interval matrices with entries from $N_c$ (<Q∪I>). T is a S-semigroup, T has subsemigroups and S-subsemigroups. T is of infinite order and is commutative.

*Example 9.21*: Let P = {all 5 × 5 neutrosophic interval matrices with entries from $N_{oc}$ (<$Z_{10}$∪I>)} be a semigroup under multiplication P is of finite order. P is non-commutative and has zero divisors and idempotents. P is a S-semigroup. Now P has pure neutrosophic subsemigroup as well as subsemigroups which are not neutrosophic.

T = {all 5 × 5 interval matrices with entries from $N_{oc}$ ($Z_{10}$)} ⊆ P is a subsemigroup of P.

W = {all 5 × 5 interval matrices with entries from $N_{oc}$($Z_{10}$I)} ⊆ P is a pure neutrosophic interval subsemigroup.

Now we can define neutrosophic interval groups and pure neutrosophic interval groups the definition is a matter of routine.

We will illustrate this by some examples.

*Example 9.22*: Let V = {all 3 × 3 neutrosophic intervals from $N_o$(ZI)} be a neutrosophic interval group under addition.

*Example 9.23*: Let P = {$N_c$(QI)} be an interval neutrosophic group under addition.

*Example 9.24*: W = {$N_{oc}$(<Z∪I>)} be an interval neutrosophic group under addition.

*Example 9.25*: M = {All 3 × 8 neutrosophic interval matrices with intervals $N_c$(<Q∪I>)} be a neutrosophic interval group under addition.

*Example 9.26*: Let F = {$N_c$(QI \ {0})} be an interval neutrosophic S-semigroup.

*Example 9.27*: Let M = {all 3 × 9 neutrosophic intervals with entries from $N_o$(<R∪I>)} be a neutrosophic interval group under addition, M is of infinite order.

*Example 9.28*: Let W = {1 × 9 neutrosophic intervals with entries from $N_{oc}$ (<Q∪I>)}. W under addition is a group of infinite order.

*Example 9.29*: M = {9 × 1 neutrosophic interval matrices with entries from $N_{oc}$ (<Z∪I>)} be a neutrosophic interval group under addition.

We can define subgroups and Smarandache neutrosophic special definite groups as in case of usual interval groups.

This task is also left as an exercise to the reader.

*Example 9.30*: Let V = {$N_c(Z_{12}I)$} be a group under addition. V is of finite order and has subgroups.

*Example 9.31*: Let T = {all 3 × 4 neutrosophic interval matrices with entries from $N_{oc}$ (<$Z_{15}$∪I>)} be a neutrosophic interval group under addition. T is of finite order. T has subgroups T is abelian.

*Example 9.32*: Let M = {$N_{oc}$ (ZI)} be a pure neutrosophic interval group under addition.

Take T = {N($Z^+I$ ∪ {0})} ⊆ M, T is a pure neutrosophic Smarandache special definite group.

Next we proceed on to define the notion of pure neutrosophic semirings and neutrosophic semirings of infinite order.

We can also define neutrosophic interval matrix semirings and neutrosophic interval polynomial semirings of infinite order. We give the properties related with them. However most of the properties can be proved by the interested reader.

**DEFINITION 9.2**: *Let S = {$N_c(Z^+I$ ∪ {0})} be the set of natural class of closed neutrosophic intervals; (S, + ×) is a pure neutrosophic semiring of integers.*

Infact S is a strict semiring with zero divisors.

So $Z^+I$ ∪ {0} is not a semifield.

Clearly $Z^+I$ can be replaced by $Q^+I$ or $R^+I$, still the results continue to be true.

We will just give examples of semirings built using the natural class of neutrosophic intervals (closed, open or open-closed or closed-open).

***Example 9.33***: Let $S = \{N_o(5Z^+I \cup \{0\})\}$ S be a pure neutrosophic semiring of open intervals.
Clearly S is commutative but has no multiplicative unit. S is of infinite order, S has zero divisors.

***Example 9.34***: Let $P = \{N_{oc}(Q^+I \cup \{0\})\}$ be a neutrosophic semiring of rationals built using open closed intervals of infinite order. P is commutative.

***Example 9.35***: Let $T = \{N_o(R^+I \cup \{0\}), +, \times\}$ be a neutrosophic real semiring of infinite order.
Now having seen some direct examples we now proceed onto give other examples of pure neutrosophic semirings.

***Example 9.36***: Let $S = \left\{\sum_{i=0}^{\infty} [a_i\ b_i]\ x^i \mid [a_i, b_i] \in N_c(Q^+I \cup \{0\})\right\}$ be a pure neutrosophic polynomial semiring of infinite order.

***Example 9.37***: Let $R = \left\{\sum_{i=0}^{\infty} [a_i\ b_i]\ x^i \mid [a_i, b_i] \in N_c(R^+I \cup \{0\})\right\}$ be a pure neutrosophic polynomial semiring of infinite order.

***Example 9.38***: Let $M = \left\{\sum_{i=0}^{\infty} [a_i\ b_i)\ x^i \mid [a_i, b_i) \in N_{co}(Q^+I \cup \{0\})\right\}$ be the pure neutrosophic polynomial ring of infinite order which has zero divisors so is not a semifield.
Now as a matter of routine interested reader can define semiring of polynomial matrices. However we give examples of them.

***Example 9.39***: Let $S = \left\{ \begin{bmatrix} a_1 & a_2 \\ a_4 & a_3 \end{bmatrix} \middle| a_i \in N_{co}(Z^+I \cup \{0\}), 1 \le i \le 4 \right\}$ be a matrix pure neutrosophic semiring. Clearly S is not a semifield.

***Example 9.40***: Let T = {all 3 × 3 pure neutrosophic interval matrices with entries from $N_{co}(Z^+I \cup \{0\})$}. T is a semiring with usual matrix addition and multiplication; T is not a semifield. T has zero divisors. T is non commutative and of infinite order.

***Example 9.41***: Let S = {all 5 × 5 upper triangular matrices with intervals from $N_o(Q^+I \cup \{0\})$}, S is a semiring of infinite order and is noncommutative.

Now we will define semiring which are neutrosophic but not pure neutrosophic. Consider S = {[$a_i$, $b_i$] | [$a_i$, $b_i$] $\in N_c(Z^+I \cup \{0\})$}; S under usual addition and multiplication is a semiring called the neutrosophic interval semiring. We can replace $<Z^+ \cup I> \cup \{0\}$ by $<Q^+ \cup I> \cup \{0\}$ or $<R^+ \cup I> \cup \{0\}$.

We will only illustrate these situations by some examples.

***Example 9.42***: Let $S = (N_{oc} (<R^+ \cup I>) \cup \{0\}) = \{(a + bI, c + dI] | a, b, c, d, \in (R^+ \cup \{0\}\}$ S is a neutrosophic semiring under usual addition and multiplication.

We will show how the addition and multiplication are carried out on S.

Take x = (3 – 4I, 2 + 5I] and y = (2 + 5I, 3 + I] in S. Now x · y = (3 – 4I × 2 + 5I, 2 + 5I × 3 + I] = (6 – 8I + 15I - 20I, 6 + 15I +2I + 5I] = (6 – 13I, 6 + 22I].

***Example 9.43***: Let $S = (N_o (<Z^+ \cup I> \cup 0))$, S is a neutrosophic semiring of natural class of intervals. That is S = {(a + bI, c + dI) | a, b, c, d $\in Z^+ \cup \{0\}$} is neutrosophic semiring with + and × as the operations.

***Example 9.44***: Let $M = \{N_o (<Q^+ \cup I> \cup \{0\})\}$ be a neutrosophic interval semiring.

Now having seen natural class of interval neutrosophic semiring we can proceed onto give examples of substructures in them. However the definition of them is direct and hence is left as an exercise to the reader.

***Example 9.45***: Let $S = \{N_{oc} (<R^+ \cup I> \cup \{0\})\}$ be a neutrosophic semiring. Consider $P = \{N_{oc} (<Q^+ \cup I> \cup \{0\})\} \subseteq S$; P is a neutrosphic subsemiring. Clearly the neutrosophic interval semiring has zero divisors but is a strict neutrosophic interval semiring.

Infact S has several neutrosophic interval subsemirings.

***Example 9.46***: Let $S = \{N_c (<Z^+I \cup \{0\}>)\}$ be a pure neutrosophic semiring. Consider $P = \{N_c (3Z^+I \cup \{0\})\} \subseteq S$; P is a neutrosophic subsemiring of S. Infact S has infinitely many neutrosophic subsemirings.

***Example 9.47***: Let $S = \{N_{co} (<Q^+ \cup I> \cup \{0\})\}$ be a neutrosophic semiring. We see $P = \{N_{co} (Q^+ \cup \{0\})\} \subseteq S$ is also an interval semiring but P is not a neutrosophic semiring. Thus we call subsemirings in S which are not neutrosophic subsemirings as pseudo neutrosophic interval subsemirings.

The following theorems are direct and is left as an exercise to the reader.

**THEOREM 9.1** : *Let $S = \{N_c (Z^+I \cup \{0\})\}$ be a semiring ($N_c(Z^+I \cup \{0\})$ can be replaced by $N_c(Q^+I \cup \{0\})$ or $N_c(R^+I \cup \{0\})$ also closed intervals by open intervals or half open-closed intervals or half closed-open intervals). S is a pure neutrosophic interval semiring and has no pseudo neutrosophic interval subsemirings.*

**THEOREM 9.2**: *Let $S = \{N_c (<R^+ \cup I> \cup \{0\})\}$ be a neutrosophic interval semiring. S has pseudo neutrosphic interval subsemirings.*

Proof is direct.

Clearly $(N_c(<R^+ \cup I> \cup \{0\})$ can be replaced by $N_c(<Q^+ \cup I> \cup \{0\})$ or $N_c(<Z^+ \cup I> \cup \{0\})$ and still the conclusion of the theorem holds good. Further closed interval is replaced by open or open-closed or closed-open intervals still the conclusion of the theorem holds good.

*Example 9.48*: S = $\{N_o((<Q^+ \cup I>) \cup \{0\})\}$ be an interval neutrosophic semiring.

Clearly $\{Q^+ \cup \{0\}\}$ = P corresponds to the degenerate intervals of the form (a, a) where a $\in Q^+ \cup \{0\}$ is a strict subsemiring which is semifield and is not a neutrosophic semifield. Thus S is a pseudo neutrosophic S-ring.

Also T = $\{Q^+I \cup \{0\}\} \subseteq$ S; and $Q^+I$ is a neutrosophic semifield and a pure neutrosophic subsemiring of S. So S is a neutrosophic S-ring.

Inview if this we have the following theorem.

**THEOREM 9.3**: *Let S = $\{N_c(<Z^+ \cup I> \cup \{0\})$ be a neutrosophic interval semiring.*
        *1. S has S-ideals.*
        *2. S has zero divisors.*

If $Z^+$ is replaced by $Q^+$ or $R^+$ still the conclusions of the theorem hold good.

*Example 9.49*: Let S = $\{N_c (Z^+I \cup (0))\}$ be a pure neutrosophic interval semiring. S contains the degenerate pure neutrosophic set P = $\{[a, a] \mid a \in Z^+I \cup \{0\}\}$ and P is a neutrosophic semifield.

Thus S is a S-ring. Hence we can say all interval semirings built using $Z^+ \cup \{0\}$ or $Q^+ \cup \{0\}$ or $R^+ \cup \{0\}$ are S-semirings.

However we have a class of interval semirings which are not S-semirings.

We will illustrate this situation by an example.

*Example 9.50*: Let S = $\{N_c (3Z^+I \cup \{0\})\}$ be an interval neutrosophic semiring. S is not a S-semiring.

*Example 9.51*: Let T = $\{N_o (12Z^+I \cup \{0\})\}$ be an interval neutrosophic semiring.
      T is not a S-semiring.

**THEOREM 9.4**: *Let $S = \{N_c (nZ^+I \cup \{0\}) \mid 2 \leq n \in Z^+\}$ be a neutrosophic interval semiring. S is not a S-semiring.*

Proof is direct (if closed interval is replaced by open interval or open-closed interval or closed-open interval then also the conclusions of the theorem holds good).

***Example 9.52***: Let $S = \{N_o (<5Z^+ \cup I> \cup \{0\})\}$ be a neutrosophic interval semiring. S is not a S-semiring.

***Example 9.53***: Let $S = \{N_{co} (<6Z^+ \cup I> \cup \{0\})\}$ be a neutrosophic interval semiring S is not a S-semiring.
    In view of these two examples we proceed onto state the following theorem. The proof of which is direct.

**THEOREM 9.5**: *Let $T = \{N_o (<nZ^+ \cup I> \cup \{0\})\}$ be a neutrosophic interval semiring, $2 \leq n < \infty$. T is not a S-semiring.*

Thus we have an infinite class of neutrosophic interval semirings and pure neutrosophic interval semirings which are not S-semirings. But all these classes of neutrosophic interval semirings be them S-semirings or not have zero divisors but have no idempotent or nilpotent.
    Now we will give other type of interval semirings and neutrosophic interval semirings.

***Example 9.54***: Let $S = \{([a_1, b_1] [a_2, b_2] [a_3, b_3]) \mid [a_i, b_i] \in N_c (Z^+ \cup \{0\})\}$ be a interval semiring under usual addition and multiplication.

***Example 9.55***: Let $S = \left\{ \begin{bmatrix} [a_1,b_1] & [a_2,b_2] \\ [a_3,b_3] & [a_4,b_4] \end{bmatrix} \right.$ where $[a_i, b_i] \in$ $N_c(Q^+ \cup \{0\})\}$ be a natural class of interval semiring under usual matrix addition and multiplication.
    If we replace $N_c (Q^+ \cup \{0\})$ by $N_c (Q^+ I \cup \{0\})$ or $N_c (<Q^+ \cup I> \cup \{0\})$ we get the neutrosophic interval semirings.

**Example 9.56**: Let W = {7 × 7 interval matrices with entries from $N_{oc}(Z_{92})$} be the matrix interval semiring of finite order. W is a S-semiring.

**Example 9.57**: Let M = {all n × n interval matrices from $N_c(Z^+ \cup \{0\})$} be a interval matrix semiring.

Clearly M is of infinite order commutative, have zero divisors and is a S-semiring.

If we replace in these examples the intervals by neutrosophic intervals and pure neutrosophic intervals the results and examples are true. Interested reader can develop these concepts and give examples of them which is a matter of routine. Now we proceed onto give examples of interval polynomial semirings.

**Example 9.58**: Let $S = \left\{ \sum_{i=0}^{\infty} [a_i, b_i] x^i \mid [a_i, b_i] \in N_c(Z^+ \cup \{0\}) \right\}$ be the polynomial interval semiring of infinite order.

**Example 9.59**: Let $R = \left\{ \sum_{i=0}^{\infty} (a_i, b_i] x^i \mid (a_i, b_i] \in N_{oc}(Q^+ \cup \{0\}) \right\}$ is an interval polynomial semiring which is a S-semiring.

**Example 9.60**: Let $S = \left\{ \sum_{i=0}^{\infty} (a_i, b_i) x^i \mid (a_i, b_i) \in N_o(Z_n) \right\}$ is an interval polynomial semiring.

Now in these semirings if $Z^+ \cup \{0\}$ is replaced by $Z^+I \cup \{0\}$ or $\langle Z^+ \cup I \rangle \cup \{0\}$ or $Q^+I \cup \{0\}$ or $R^+I \cup \{0\}$ or $\langle R^+ \cup I \rangle \cup \{0\}$ we get the interval neutrosophic polynomial semiring.

**Example 9.61**: Let $T = \left\{ \sum_{i=0}^{\infty} [a_i, b_i] x^i \mid [a_i, b_i] \in N_c(Z^+I \cup \{0\}) \right\}$ is a neutrosophic interval polynomial semiring of infinite order which is commutative.

***Example 9.62***: Let $R = \left\{ \sum_{i=0}^{\infty} (a_i, b_i) x^i \mid (a_i, b_i) \in N_o(<Z_n \cup I>) \right\}$ be a neutrosophic interval polynomial semiring R is a S-semiring.

All polynomial semirings built using intervals $N_o(<Z^+ \cup I> \cup \{0\})$ $N_c(<Q^+ \cup I> \cup \{0\})$ and $N_{oc}(<R^+ \cup I> \cup \{0\})$ are all S-semirings.

As in case of polynomial interval semirings we can in case of polynomial interval neutrosophic semirings also derive all related properties as a matter of routine.

We can have all substructures and special elements like units, zero divisors, idempotents and their Smarandache analogue.

Neutrosophic interval semirings also follow all the properties with appropriate simple modifications [13].

Now we proceed onto define neutrosophic interval rings using $N_c(Z)$ or $N_{oc}(R)$ or $N_o(Q)$ and so on.

**DEFINITION 9.3** : *Let $S = \{N_c(ZI)\}$ be the collection of all closed intervals.*
  1. *S is a group under addition (S is a commutative group).*
  2. *S is a semigroup under multiplication.*
  3. *The operation on S distributes that is $a(b + c) = ab + ac$ and $(a + b)c = ac + bc$ for all a, b, c in $N_c(Z)$.*

*Thus S is a ring defined as the pure neutrosophic interval ring of integers.*

If ZI is replaced by $<Z \cup I>$ we just get the neutrosophic ring of integers.

If ZI is replaced by QI (or $<Q \cup I>$) in definition 9.3 we get the pure neutrosophic interval ring of rationals (or neutrosophic ring of rationals).

If ZI is replaced by RI (or $<R \cup I>$) in definition 9.3 we get the pure neutrosophic interval ring of reals (or neutrosophic ring of reals).

If ZI is replaced by $Z_nI$ (or $<Z_n \cup I>$) in the definition we get the pure neutrosophic ring of modulo integers (or neutrosophic ring of modulo integers).

We see the closed interval in the definition can be replaced by open intervals, open-closed intervals or closed-open intervals and the definition continues to hold good.

We now proceed onto give some examples.

**Example 9.63**: Let $R = \{N_{oc}(ZI)\} = \{(a_i, b_i] \mid a_i, b_i \in ZI\}$ be a pure neutrosophic ring which is commutative and is of infinite order.

**Example 9.64**: Let $W = \{N_c(QI)\} = \{[a_i, b_i] / a_i, b_i \in QI\}$ be the pure neutrosophic commutative ring of infinite order.

**Example 9.65**: Let $M = \{N_{co}(RI)\} = \{[a_i, b_i) / a_i, b_i \in RI\}$ be the pure neutrosophic commutative ring of infinite order.

**Example 9.66**: Let $P = \{N_o(Z_{20}I)\} = \{(a, b) \mid a, b, \in Z_{20}I\}$ be a pure neutrosophic commutative ring of finite order.

**Example 9.67**: Let $R = \{N_c(<R \cup I>)\} = \{[a, b] \mid a, b \in <R \cup I>\}$ be the neutrosophic ring of reals of closed neutrosophic intervals.

We will just show how addition and multiplication are performed in R.
Consider $x = [5 + 2I, -7 + 5I]$ and $y = [-3 + 8I, -I]$ in R.
$x + y = [2 + 10I, -7 + 4I]$ is in R.
$x \cdot y = [5 + 2I, -7 + 5I] \times [-3 + 8I, -I]$
$= [-15 - 6I + 40I + 16I, 7I - 5I]$
$= [-15 + 50I, 2I]$ is in R.

**Example 9.68**: Let $R = \{N_{oc}(<Z_7 \cup I>)\}$ be a ring of finite order. R is a neutrosophic interval ring.

**Example 9.69**: Let $M = \{N_{oc}(<Z_{240} \cup I>)\}$ be a commutative neutrosophic interval ring with unit of finite order.

**Example 9.70**: Let $W = \{N_{oc}(<3Z^+ \cup I>)\}$ be a neutrosophic interval ring of finite order and has no unit.

***Example 9.71***: Let M = {N$_o$ (PI) | P = {0, 2, 4, 6, 8, 10, 12, 14, 16, 18} ⊆ Z$_{20}$} be a finite neutrosophic interval ring with no unit and M is commutative.

We can define substructures and special elements as in case of ring.

We will illustrate these by examples and leave the task of defining them to the reader as it is a matter of routine.

***Example 9.72***: Let R = {N$_c$ (ZI)} be the neutrosophic interval ring, T = {N$_c$ (8ZI)} ⊆ R be a subring as well as an ideal of R. Clearly T is not a maximal ideal.

Consider M = {N$_c$ (3ZI)} ⊆ R is a subring as well as an ideal of M but is a maximal ideal of R.

In fact R has infinitely many maximal ideals and subrings.

***Example 9.73***: Let W = {N$_o$ (<R∪I>)} be a neutrosophic interval ring. Take M = {(0,a) | a ∈ <R∪I>} ⊆ W is again a subring as well as ideal of W.

Infact T = {(a, 0) / a ∈ RI ⊆ < R ∪ I>} ⊆ W is again a subring as well as ideal of W.

But V = {(0, a) | a ∈ N$_o$ (R) ⊆ N$_o$ (<R∪I>)} is only a subring and not an ideal of V.

Thus we have subrings which are not ideals of W in case of neutrosophic rings also.

We cannot in general define neutrosophic interval fields.

We have only S-neutrosophic interval rings that too not all interval rings are S-neutrosophic rings. We see in general neutrosophic ring N$_c$(<nZ$^+$ ∪ I> ∪ {0}) are not S-neutrosophic for n = 1, 2, … .

Further we do not have interval fields or neutrosophic fields using these natural class of intervals. We can use only S-interval rings. Likewise we do not have interval semifield or neutrosophic interval semifields we only use S-semirings in building the S-interval semivector spaces and S-neutrosophic interval semivector spaces. Also use S-rings to study S-linear algebra and S-neutrosophic linear algebra.

Now the concept of substructures zero divisors, idempotents etc. can be derived as in case of usual rings [16]. The definition of neutrosophic interval vector space can be defined on three types of S-rings viz., and over fields and neutrosophic fields.

We will just briefly give the definition of these concepts.

Let $N_c(ZI)$ (or $N_c(<Z \cup I>)$ or $N_o(QI)$ or $N_{co}(RI)$ or $N_{oc}(<Q \cup I>)$ or $N_{co}(<R \cup I>)$ be S-neutrosophic interval S-rings.

Q or R or $Z_p$ be usual fields. QI or RI or $Z_pI$ be pure neutrosophic fields (p a prime).

**DEFINITION 9.4**: *Let V be a interval neutrosophic abelian group under addition with entries from $N_o(ZI)$ (or $N_o(<Z \cup I>)$ or $N_c(QI)$ or $N_{oc}(<Q \cup I>)$ or $N_{co}(RI)$ or $N_{oc}(<R \cup I>)$); here they can be open intervals or closed intervals or half open-closed or half closed open intervals). Suppose F is a real field Q or R then we define V to be a neutrosophic interval vector space of type I over F if for all $v \in V$ and $a \in F$, va and av $\in V$ and*

1. *$(v_1 + v_2) a = v_1 a + v_2 a$*
2. *$(a + b) v_1 = a v_1 + b v_1$*
3. *$0 \cdot v = 0$*
4. *$1 \cdot v = v$ for all $v, v_1, v_2$ in V and $a, b \in F$.*

We will give an example of this.

***Example 9.74***: Let $V = N_o(QI)$ be an additive group of neutrosophic intervals. Let $F = Q$ be the field. V is a neutrosophic interval vector space over Q.

***Example 9.75***: Let $V = N_{oc}(<Q \cup I>)$ be a neutrosophic interval vector space over Q. Infact V is a linear algebra.

We see V is not a neutrosophic interval vector space over R.

***Example 9.76***: Let $V = N_c(RI)$ be neutrosophic interval vector space over Q (or over R).

Clearly V is a pure neutrosophic interval linear algebra over Q (or over R).

***Example 9.77***: Let $S = \{N_{oc}(<R \cup I>)\}$ be a neutrosophic vector space over R of type I.

Now finding subspaces, basis and linear transformations in case of neutrosophic interval vector spaces of type I is a matter of routine and is left as an exercise to the reader.

Now we proceed onto give one example of a finite neutrosophic interval vector spaces of type I.

**Example 9.78**: Let $V = \{N_c(Z_5I)\}$ be a neutrosophic interval vector space of type I over the field $Z_5$.

We now define type II neutrosophic interval vector spaces.

**DEFINITION 9.5**: *Let $V = N_c(QI)$ (or $N_c(<Q\cup I>))$ be the collection of all neutrosophic intervals. Consider the neutrosophic field QI. V is a vector space over QI; defined as the neutrosophic interval vector space of type II.*

We can replace QI by RI or $<R\cup I>$ and the result will be true.
We will give examples of them.

**Example 9.79**: Let $V = \{N_{oc}(QI)\}$ be a neutrosophic interval vector space of over the neutrosophic field $QI = F$ of type II.
Clearly dimension of V over F is 2. The basis of V over F is $\{(I, 0], (0, I]\}$. For any $x = (aI, bI] = a(I, 0] + b(0, I]$.

**Example 9.80**: Let $V = \{N_c(<R\cup I>)\}$ be a neutrosophic interval vector space over the neutrosophic field RI of type II. V is also a neutrosophic interval linear algebra of type II.

**Example 9.81**: Let $V = \{(a_1, a_2, \ldots, a_8) \mid a_i \in N_c(QI); 1 \leq i \leq 8\}$ be a neutrosophic interval vector space over the field Q of type I.

**Example 9.82**: Let
$$W = \left\{ \begin{bmatrix} a_1 \\ a_2 \\ a_3 \\ a_4 \end{bmatrix} \middle| a_i \in N_o(RI); 1 \leq i \leq 4 \right\}$$

be a neutrosophic interval vector space over Q of type I. W is infinite dimensional and W is not a neutrosophic interval linear algebra over Q.

*Example 9.83*: Let M = {all 3 × 5 neutrosophic interval matrices with entries from $N_o$ (QI)} be a neutrosophic interval vector space over the field Q of type I. Clearly M is not a neutrosophic interval linear algebra over Q.

*Note*: If we define V or W or M over the neutrosophic field QI we would get neutrosophic interval vector space of type II.

*Example 9.84*: Let

$$V = \left\{ \sum_{i=0}^{\infty} (a_i, b_i) x^i \mid (a_i, b_i) \in N_o (RI) \right\}$$

be a neutrosophic interval vector space of type I over Q (or R). (If V is defined on the field QI (or RI) then V is a neutrosophic interval vector space of type II).
   Clearly this V is a neutrosophic interval linear algebra over Q (or R) of type I (or type II when defined over QI or RI).

*Example 9.85*: Let V = {9 × 9 neutrosophic interval matrices with entries from $N_c(<Q \cup I>)$} be the neutrosophic interval vector space over the field Q of type I (or over the field QI of type II). Clearly V is a neutrosophic interval linear algebra of type I over Q (or type II over QI).
   We see as in case of usual interval vector spaces in case of neutrosophic interval vector spaces of type I and II also we see every neutrosophic linear algebra is a neutrosophic vector space and not conversely.

Next we proceed onto define type III and type IV neutrosophic interval vector spaces.

**DEFINITION 9.6**: *Let $V = \{N_c (QI)\}$ (or $N_o (QI)$ or $N_{oc} (RI)$ or $N_{co} (<Q \cup I>)$ or $N_o (<R \cup I>)$ and so on) be a Smarandache neutrosophic interval vector space over the interval S-ring $N_c(Z)$ (or $N_c(Q)$ or $N_c(R)$). We call V a S-neutrosophic interval vector space of type III over the interval S-ring $N_c(Q)$ (or $N_c(Z)$). It is to be noted if V contains closed intervals then the interval S-ring must also be closed intervals.*

It cannot be defined over open intervals or half open-closed intervals or half closed-open intervals.
   Likewise if V contains open intervals the interval S-ring must also be a open interval S-ring of type III.

Further if V is the collection of half open-closed neutrosophic intervals then V can only be defined over the half open-closed interval S-ring. We will give examples of them.

***Example 9.86***: Let $T = \{N_{oc} (QI)\}$ be the S-neutrosophic interval vector space over the interval S-ring $F = N_{oc}(Q)$ (or $F = N_{oc} (Z)$).

***Example 9.87***: Let $V = N_o(RI)\}$ be a S-neutrosophic interval vector space over the interval S-ring $F = N_o(Q)$.

***Example 9.88***: Let $W = \{N_{co} (QI)\}$ be a S-neutrosophic interval vector space over the interval S-ring $F = N_{co} (Q)$.

***Example 9.89***: Let $P = \{$all $3 \times 2$ neutrosophic interval matrices with entries from $N_{co} (<R \cup I>)$ be a S-neutrosophic vector space of type III over the interval S-ring $F = N_{co}(Q)$.

***Example 9.90***: Let

$$V = \left\{ \sum_{i=0}^{\infty} (a_i, b_i) x^i \mid (a_i, b_i] \in N_{oc} (<Q \cup I>) \right\}$$

be a S-neutrosophic interval linear algebra of type III over the interval S-ring $F = N_{oc}(Q)$.

*Example 9.91*: Let S = {all 10 × 10 neutrosophic interval matrices with entries from $N_c(<Q \cup I>)$} be a S-neutrosophic interval linear algebra of type III over the interval S-ring $W = N_c(Q)$.

Now having seen the definition and some examples all the notions related to vector spaces can be derived as a matter of routine without any difficulty.

Now we proceed onto define type IV S-neutrosophic interval vector spaces defined over the neutrosophic interval S-rings.

**DEFINITION 9.7**: *Let $V = N_c(<R \cup I>)$ be an additive abelian group of neutrosophic real intervals. V is a Smarandache neutrosophic interval vector space of type IV over the Smarandache neutrosophic interval S-ring $F = N_c(QI)$ (or $N_c(<Q \cup I>)$ or $N_c(<R \cup I>)$ or $N_c(RI)$) if the vector space postulates are true or equivalently we can say in the definition of type III vector spaces replace the interval S-ring by neutrosophic interval S-ring.*

We will illustrate this situation by some examples.

*Example 9.92*: Let M = {$(a_1, a_2, a_3, a_4) | a_i \in N_{co}(<Q \cup I>)$} be a S-neutrosophic interval vector space of type IV over the neutrosophic interval S-ring $T = N_{co}(QI)$.

*Example 9.93*: Let B = {all 5 × 2 neutrosophic interval matrices with entries from $N_{oc}(QI)$} be a S-neutrosophic interval vector space of type IV over the S neutrosophic interval S-ring $F = N_{oc}(QI)$. Clearly B is not a S-neutrosophic interval linear algebra over F of type IV.

*Example 9.94*: Let M = {8 × 8 neutrosophic interval matrices with entries from $N_c(<R \cup I>)$} be a S-neutrosophic interval linear algebra of type IV over the neutrosophic interval S-ring F = $N_c(RI)$.

*Example 9.95*: Let P = {2 × 9 neutrosophic interval matrices with entries from $N_{oc}$ (QI)} be a S-neutrosophic interval vector space over the neutrosophic interval S-ring F = $N_{oc}$ (QI). Clearly P is not a S-neutrosophic interval linear algebra of type IV. Other related properties for these structures can also be derived as a matter of routine.

Now we proceed onto define four types of neutrosophic interval semivector spaces.

**DEFINITION 9.8**: *Let V be a semigroup under addition with zero built using the natural class of neutrosophic – intervals. If V is a semivector space over a real semifield. We define V to be a neutrosophic semivector space of intervals of type I.*

We will illustrate this situation by some examples.

*Example 9.96*: Let

$$V = \left\{ \begin{bmatrix} a_1 \\ a_2 \\ a_3 \\ a_4 \end{bmatrix} \middle| a_i \in N_o (Z^+I \cup \{0\}); 1 \leq i \leq 4 \right\}$$

be a semivector space of neutrosophic intervals or neutrosophic interval semivector space over the semifield $S = Z^+ \cup \{0\}$ of type I.

*Example 9.97*: Let V = $N_c$ ($Q^+I \cup \{0\}$) be a semivector space of neutrosophic intervals of type I over the semifield $S = Z^+ \cup \{0\}$.

*Example 9.98*: Let M = {all 2 × 2 neutrosophic interval matrices with entries from $N_{oc}$ ($<Z^+ \cup I> \cup \{0\}$)} be a semivector space of neutrosophic interval matrices over the semifield $S = Z^+ \cup \{0\}$ of type I. Infact M is a semilinear algebra of type I.

*Example 9.99*: Let

$$V = \left\{ \sum_{i=0}^{\infty} [a_i, b_i) x^i \text{ where } [a_i, b_i) \in N_{co}(Q^+I \cup \{0\}) \right\}$$

be a neutrosophic interval polynomial semivector space of type I over the semifield $F = Q^+ \cup \{0\}$. Infact V is a linear algebra of type I over F.

*Example 9.100*: Let V = {3 × 7 neutrosophic interval matrices with entries from $(a_i, b_i] \in N_{oc}(\langle R^+ \cup I \rangle \cup \{0\})$} be a neutrosophic interval polynomial semivector space of type I over the semifield $S = Q^+ \cup \{0\}$. Clearly V is not a semilinear algebra over S of type I.

We can as in case of semivector spaces define substructures, transformations and so on for type I semivector spaces. This task is easy and it can be taken as a matter of routine.

Now if in the definition replace the real semifield by neutrosophic semifield we define these neutrosophic interval semivector as type II semivector space.

We will illustrate this situation by some examples.

*Example 9.101*: Let $V = \{N_{oc}(Z^+I \cup \{0\})\}$ be a neutrosophic interval semivector space of type II over the neutrosophic semifield $F = Z^+I \cup \{0\}$.

Clearly V is not a neutrosophic interval semivector space of type II over the neutrosophic semifield $QI \cup \{0\}$.

*Example 9.102*: Let M = {2 × 2 neutrosophic interval matrices with entries from $N_{oc}(Q^+I \cup \{0\})$} be a neutrosophic interval semivector space of type II over the neutrosophic semifield $S = Q^+I \cup \{0\}$.

*Example 9.103*: Let T = {all 3 × 2 neutrosophic interval matrices with entries from $N_c(Z_{19}I)$} be a neutrosophic interval semivector space over the neutrosophic semifield $S = Z_{19}I$. Clearly T is not a neutrosophic interval semi linear algebra.

**Example 9.104**: Let P = {all 6 × 6 neutrosophic interval matrices with entries from $N_o (\langle R^+ \cup I \rangle \cup \{0\})$} be the interval neutrosophic semivector space of type II over the neutrosophic semifield $S = R^+I \cup \{0\}$.

The notion of defining and analyzing the properties of basis, substructures and transformations or operations is a matter of routine.

Next if in the definition the semifield is replaced by the natural class of real intervals we get the S-neutrosophic interval semivector space of type III.

We will illustrate this situation by some examples.

**Example 9.105**: Let $S = \{N_{oc} (Z^+I \cup \{0\})\}$ be a Smarandache neutrosophic interval semivector space of type III over the real interval S-semiring $F = N_{oc} (Z^+ \cup \{0\})$.

**Example 9.106**: Let $V = \{N_c (R^+I \cup \{0\})\}$ be the Smarandache neutrosophic interval semivector space of type III over the real interval S-semiring $F = N_c (Z^+ \cup \{0\})$.

**Example 9.107**: Let $V = \{N_{oc} (Q^+I \cup \{0\})\}$ be the Smarandache neutrosophic interval semivector space of type III over the interval semiring $F = N_{oc} (Q^+ \cup \{0\})$.

Now all properties related with these S-neutrosophic interval semivector spaces can be derived in case of type III S-semivector spaces also with appropriate modifications without any difficulty.

Next in the definition if we replace the semifield by the neutrosophic interval S-semiring then we get type IV interval neutrosophic S-semivector spaces or S-interval neutrosophic semivector spaces.

We will illustrate this situation by some examples.

***Example 9.108***: Let $V = N_c (Z^+I \cup \{0\})$ be a Smarandache neutrosophic interval semivector space of type IV over the Smarandache neutrosophic semiring $F = N_c(Z^+I \cup \{0\})$.

***Example 9.109***: Let $S = \{N_{oc} (Q^+I \cup \{0\})\}$ be a Smarandache neutrosophic interval semivector space of type IV over the neutrosophic interval S-semiring. $N_{oc}(Z^+I \cup \{0\})$.

***Example 9.110***: Let $W$ = {all $2 \times 2$ neutrosophic interval matrices with entries from $N_c(<R^+ \cup I> \cup \{0\})$} be a S-neutrosophic interval semivector space of type IV over the neutrosophic interval S-semiring $N_c(Z^+I \cup \{0\})$.

All properties can be derived for these type IV neutrosophic interval S-semivector spaces also.

Finally we just indicate how all these results can be carried out if we replace neutrosophic real intervals $N_c(<R \cup I>)$ built using by $N_c (<[0,1] \cup [0,I]>)$ where $N_c (<[0,1] \cup [0,I]>) = \{[a, b]$ where $a = a_1 + a_2I$ and $b = b_1 + b_2I$ with $a_1, a_2, b_1, b_2 \in [0,1]\}$.

These intervals will be known as fuzzy neutrosophic intervals. $N_{oc} (<[0,1] \cup [0,I]>)$, $N_{co}(<[0,1] \cup [0, I]>)$ and $N_o (<[0, 1] \cup [0, I]>)$ can also be defined. All results so far studied and defined for real neutrosophic intervals hold good for fuzzy neutrosophic intervals also.

By default of notations we accept [a,b] or (a, b) or [a, b) or (a, b] as an interval even if a and b are not comparable.

Thus $(0.5 + 0.8I, 0.9 + 0.3I)$ $(0, 0.I)$ $(0.I, 0.03)$, $(0.003I, 0.2)$ are all neutrosophic fuzzy intervals. Thus we at this stage relinquish the comparison of the components of the interval.

**Chapter Ten**

# APPLICATIONS OF THE ALGEBRAIC STRUCTURES BUILT USING NATURAL CLASS OF INTERVALS

These natural class of intervals are basically introduced by the authors so that all classical arithmetic operations of reals can be easily extended to $N_c(Q)$ or $N_c(R)$ or $N_c(Z)$ without any difficulty. We see clearly $N_c(Z) \subset N_c(Q) \subseteq N_c(R)$.

If the closed intervals are replaced by class of open intervals or open-closed intervals or by closed-open intervals, still the containment relation hold good. However we do not mix up open interval $N_o(Z)$ with the closed intervals $N_c(Z)$ or open-closed intervals $N_{oc}(Z)$ or closed-open intervals $N_{co}(Z)$.

These structures will find applications in stiffness matrices and in all places where finite element methods are used.

Further by approximating to an interval in place of a fixed value we have more flexibility in choosing the needed values.

These algebraic structures are quite new and in due course of time they will find applications. Also these structures will be more useful than the usual interval used so far, as we do not use any max-min or other operations only usual classical operations will be used. We have built algebraic structures using these natural class of intervals $N_c(Z)$ or so on. Also we have intervals constructed using neutrosophic intervals, fuzzy intervals and fuzzy neutrosophic intervals. These intervals will also find their applications in due course of time.

**Chapter Eleven**

# SUGGESTED PROBLEMS

In this chapter we suggest around 202 problems some of which are research level some simple some just difficult. Interested reader can solve them for it will improve the understanding of this notions in thic book .

1. Is $S = \{[a, b] \mid [a,b] \in N_c(Z_5)\}$ a semigroup under multiplication?

2. Does a semigroup of natural intervals of order 27 exist?

3. Is $P = \{(a, b) \mid (a, b) \in N_o(Z_2)\}$ a semigroup under addition?

4. Obtain some interesting properties about the semigroup $S = \{(a, b] / (a, b] \in N_{oc}(Z)\}$ under multiplication.

5. Let $G = \{(a, b] \mid (a, b] \in N_{oc}(Z_8)\}$ be a semigroup under multiplication.
    i. Find zero divisors in G.

ii. Can G have S-zero divisors?
iii. Find idempotents in G.
iv. Can G have S-idempotents?
v. Does G have S-units?
vi. Find the order of G.
vii. Can G have ideals?
viii. Can G have subsemigroups which are not ideals of G.?

6. Let $W = \{[a, b] \mid [a, b] \in N_c(Z)\}$ be a group under addition.
   i. Find subgroups in W.
   ii. Give a nontrivial automorphism on W.
   iii. Let $P = \{[a, b] \mid [a, b] \in N_c(3Z)\} \subseteq W$. Find W/P. What could be the algebraic structure on $W/P = \{[a, b] + P \mid [a, b] \in N_c(Z)\}$.
   iv. Is W/P of finite order?

7. Let $S = \{[a, b) \mid [a, b) \in N_{co}(Q)\}$ be the semigroup.
   i. Find ideals in S.
   ii. Find a subsemigroup of S which is not an ideal of S.
   iii. Prove S has zero divisors.
   iv. Can S have nontrivial units?
   v. Prove S cannot have nontrivial idempotents or nilpotents.

8. Let $S = \{[a, b] \mid [a, b] \in N_c(Z_{12})\}$ be a semigroup.
   i. Find the order of S.
   ii. Can S have ideals?
   iii. Can S have idemoptents?
   iv. Give some examples of idempotents and units in S.

9. Let $P = \{(a, b] \mid (a, b] \in N_{oc}(Z_{10})\}$ be a semigroup under multiplication
   i. Find the order of P.
   ii. Does P have ideals?
   iii. Does P have subsemigroups which are not ideals?
   iv. Can P have zero divisors?
   v. Does P contain S-idempotents?
   vi. Can P have S-units?
   vii. Find subsemigroups in P whose order does not divide order of P.

10. Let G = {[a_i, b_i) | [a_i, b_i) ∈ $N_{co}(Z_{40})$} be the group under addition.
   i. Find order of G.
   ii. Find subgroups of G.
   iii. Is Lagrange's theorem true in G?
   iv. Does G have Sylow subgroups?
   v. Find the order of at least 10 elements in G and verify Cauchy theorem for them.
   vi. H be a subgroup of G. Find G/H.

11. Let G = {[a_i, b_i] | [a_i, b_i] ∈ $N_c(Z_{40})$} be a semigroup under multiplication.
   i. What is the order of G?
   ii. Find zero divisors in G.
   iii. Is G a S-semigroup?
   iv. Does G have S-subsemigroups?
   v. Does G satisfy S-weak Lagrange theorem?
   vi. Can G have S-zero divisors?
   vii. Enumerate any other property related with G.
   viii. Does G contain S-Cauchy elements?

12. Let S = {[a_1, a_2], [b_1, b_2]) | [a_1, a_2], [b_1, b_2] ∈ $N_c(Z_6)$} be a semigroup under multiplication.
   i. Find the order of S.
   ii. Can S have zero divisors?
   iii. Give example of idempotents if any in S.
   iv. Find S-ideals if any in S.
   v. Find S-subsemigroups if any in S.
   vi. Does S satisfy S-weakly Lagrange theorem?

13. Let G = {(a, b) | (a, b) ∈ $N_o(Z_5)$} be a group under addition.
   i. Find the order of G.
   ii. Prove Lagrange theorem for finite groups is true in case of G.
   iii. Find atleast three subgroups of G.
   iv. Find Cauchy elements in G.
   v. Does G have p-Sylow subgroup?

vi. Find two subgroups $H_1$, $H_2$ in G $H_1 \neq H_2$ and find $G/H_1$ and $G/H_2$ so that $G/H_1 \cong G/H_2$.

14. Let $V = \{(a, b) \mid (a, b) \in N_o(Z)\}$ be a group under addition;
   i. Find subgroups of V.
   ii. Can G have an element of finite order?
   iii. Is $V \cong Z$ (Z a group under addition of positive and negative integers)?

15. Let $M = \{(a, b] \mid (a, b], \in N_{oc}(Z_6)\}$ be a semigroup under multiplication?
   i. What is the order of M?
   ii. Is M a S-semigroup?
   iii. Find S-ideals of M.

16. Let $T = \{(a, b) \mid (a, b) \in N_o(R)\}$ be a semigroup under multiplication.
   i. Find subsemigroups which are not ideals.
   ii. Is T a S-semigroup?
   iii. Find ideals in T.
   iv. Can T have S-ideals?
   v. Is $\{(1, -1), (-1, 1), (1, 1), (-1, -1)\} \subseteq T$ a group under multiplication?

17. Let $P = \{([a_1, b_1), [a_2, b_2), [a_3, b_3)) \mid [a_i, b_i) \in N_{co}(Z_{12}); 1 \leq i \leq 3\}$ be a semigroup under product.
   i. Find order of P.
   ii. Can P have ideals?
   iii. Is P a S-semigroup?
   iv. Can P have S-subsemigroups?
   v. Can P have S-ideals?
   vi. Does P have S-zero divisors?
   vii. Can P have S-idempotents?
   viii. Does P have S-Cauchy elements?

18. Let P in problem (17) be under addition ;
   i. Is P a group?
   ii. Verify Lagrange's theorem for P, if P is a group.
   iii. Find the order of P.

19. Let $V = \left\{ \begin{bmatrix} [b_1, a_1] \\ [b_2, a_2] \\ [b_3, a_3] \\ [b_4, a_4] \end{bmatrix} \middle| [b_i, a_i] \in N_c(Z_{18}); 1 \leq i \leq 4 \right\}$ be a semigroup under addition.
   i. Find order of V.
   ii. Find subsemigroups of V.

20. Let $S = \left\{ \sum_{i=0}^{5} [a_i, b_i] x^i \middle| x^6 = 1, [a_i, b_i] \in N_c(Z_7) \right\}$ be a semigroup under multiplication.
   i. What is the order of S?
   ii. Find ideals in S.
   iii. Is S a S-semigroup?
   iv. Find some subsemigroups of S which are not S-subsemigroups.

21. Let $S = N_c(Z_{10}) = \{[a, b] \mid [a, b] \in N_c(10)\}$ be the semigroup under multiplication.
   i. Find order of S.
   ii. Is S a S-semigroup?
   iii. Find S-ideals if any in S.
   iv. Find idempotents in S.
   v. Can S have S-semigroups which are not S-ideals?

22. Let $S = \{[a, b] \mid [a, b] \in N_c(Z_{13})\}$ be a semigroup under addition.
   i. Can S have ideals?
   ii. What is the order of S?
   iii. Find subsemigroups of S.
   iv. Is S a S-semigroup?
   v. Find S-subsemigroups if any in S.

23. Let $S = \{(a, b) \mid (a, b) \in N_{oc}(3Z)\}$ be a semigroup.
   i. Can S be a S-semigroup?
   ii. Does S have ideals?

   iii. Find subsemigroups of S which are not ideals of S.

24. Let $M = \{(a, b) \mid (a, b) \in N_o(R)\}$ be a semigroup under multiplication.
   i. Is M a S-semigroup?
   ii. Find ideals in M.
   iii. Find subsemigroups in M which are not ideals.
   iv. Find S-subsemigroups on M which are not S-ideals.
   v. Can M have zero divisors?
   vi. Can M have idempotents?
   vii. Can M have S-units?

25. Let $G = \{(a, b] \mid (a, b] \in N_{oc}(Z_{10})\}$ be a group under addition.
   i. What is the order of G?
   ii. Is G a Smarandache special definite group?
   iii. Find subgroups of G.
   iv. Find for some subgroup H the quotient group G/H.

26. Let $W = \{[a, b] \mid [a, b] \in N_c(Z_{15})\}$ and $V = \{[a, b] \mid [a, b] \in N_c(Z_{20})\}$ be groups under addition.
   i. Find orders of W and V.
   ii. Find a homomorphism in $\phi$ from W to V which has non trivial kernel.
   iii. Prove Cauchy theorem for V and W.
   iv. If on V and W the operation '+' is replaced by × can V and W be groups?; justify your answer.
   v. Does there exist a homomorphism $\phi$ from V and W so that V/ker $\phi \cong$ W?

27. Let $M = \{([a, b], [c, d]) \mid [a, b], [c, d] \in N_c(Z_{20})\}$ be a semigroup under multiplication.
   i. Find order of M.
   ii. Prove M has zero divisors.
   iii. Is M a S-semigroup?
   iv. Can M have S-zero divisors?
   v. Will M be a S-weakly Lagrange semigroup?
   vi. Find ideals in M.

28. Let (G, +) be a S-special definite group of natural class of intervals built using $N_{oc}$ (Z). Determine some interesting properties about G. Is Z embeddable in G?

29. Obtain some interesting properties enjoyed by $(N_c(Z_{12}), \times)$.

30. Enumerate the special properties enjoyed by $(N_{oc}(Z_{10}), +)$ and compare it with $(N_{oc}(Z_{10}), \times)$. What are the algebraic structures that can be associated with them?

31. Find some interesting properties about the rings $\{N_{oc}(R), +, \times\}$.

32. Let $R = \{N_c(Z_{40}), +, \times\}$ be a ring.
    i. Find order of R.
    ii. Find ideals in R.
    iii. Is R a S-ring?
    iv. Find subrings in R which are not ideals.
    v. Does R contain S-zero divisors?
    vi. Find S-idempotents if any in R.

33. Let $S = \{N_{co}(Z_{15}), +, \times\}$ be a ring of finite order.
    i. Find the order of S.
    ii. Find S-ideals if any in S.
    iii. Does S have S-zero divisors?
    iv. Find an ideal I in S, and what is the algebraic structure enjoyed by S/I?
    v. Can S have S-units?
    vi. Is S a S-ring?

34. Let $S = \{N_{co}(Z), +, \times\}$ be a ring. Enumerate all properties associated with S.
    i. Is S a S-ring?
    ii. Can ever S be an integral domain?
    iii. Can S have subrings which are not ideals?
    iv. Can S have S-ideals?
    v. Find an ideal I and find S/I.
    vi. Can S have S-zero divisors?

35. Let $V = \left\{ \begin{bmatrix} [a_1,b_1) \\ [a_2,b_2) \\ \vdots \\ [a_{10},b_{10}) \end{bmatrix} \middle| [a_i, b_i) \in N_{co}(Z_{20}) ; 1 \leq i \leq 10 \right\}$ be a semigroup.
    i. Prove V is finite and find its order.
    ii. Is V a S-semigroup?
    iii. Can V have zero divisors?
    iv. Find S-subsemigroup if any in V.
    v. Can V have S-ideals?
    vi. Does V satisfy S-Lagrange theorem for subgroups?

36. Let $V = \{([a_1, b_1), [a_2, b_2), \ldots, [a_8,b_8)) \mid [a_i, b_i) \in N_{co}(Z_3) ; 1 \leq i \leq 8\}$ be a semigroup under multiplication.
    i. What is the order of V?
    ii. Is V a S-semigroup?
    iii. Find S-zero divisors if any in V.
    iv. Find S-ideals if any in V.
    v. Can V have ideals which are not S-ideals?
    vi. Can V have S-Cauchy elements?
    vii. Prove V cannot be a group.

37. Let $S = \left\{ \begin{bmatrix} [b_1,a_1] & [b_4,a_4] \\ [b_2,a_2] & [b_5,a_5] \\ [b_3,a_3] & [b_6,a_6] \end{bmatrix} \middle| [b_i, a_i] \in N_c(Z_9) ; 1 \leq i \leq 6 \right\}$ be a semigroup under addition.
    i. Find the order of S.
    ii. Can S have S-zero divisors?
    iii. Find zero divisors in S.
    iv. Find S-ideals if any in S.
    v. Find subsemigroups of S.
    vi. Is S a S-semigorup?
    vii. Does S satisfy S-weakly Lagrange subsemigroup condition?
    viii. Can S have idempotents?

38. Let $S = \left\{ \sum_{i=0}^{7} [a_i, b_i) x^i \mid [a_i, b_i) \in N_{co}(Z_{11}) \right\}$ be a semigroup.
    i. Find order of S.
    ii. Does S satisfy S-weakly Lagrange theorem?
    iii. Does S have S-Cauchy elements?
    iv. Is S a S-semigroup?
    v. Can S have S-ideals?
    vi. Can S have subsemigroups which are not S-subsemigroups?
    vii. Can S have zero divisors?

39. Let $S = \left\{ \sum_{i=0}^{7} [a_i, b_i] x^i \mid [a_i, b_i] \in N_c(Z_5), x^8 = 1 \right\}$ be a semigroup under multiplication.
    i. Find the order of S.
    ii. Can S have S-ideals?
    iii. Can S have zero divisors?
    iv. Does S contain subsemigroups which are not ideals?
    v. Can S have S-subsemigroups which are not S-ideals?
    vi. Can S have idempotents?
    vii. Does S satisfy S-Lagrange theorem?

40. Let $G = \left\{ \sum_{i=0}^{12} [a_i, b_i) x^i \mid [a_i, b_i) \in N_{co}(Z_3) \right\}$ be a semigroup under addition.
    i. Find order of G.
    ii. Find subsemigroups of G.
    iii. Is G a S-semigroup?

41. Let $S = \{N_o(Z_5), +, \times\}$ be a ring.
    i. Is S a S-ring?
    ii. Find order of S.
    iii. If $I = \{(0, a) \mid a \in Z_5\} \subseteq S$ be an ideal. Find S/I.
    iv. Is S/I an integral domain or a semifield?
    v. What is the order of S/I?
    vi. Does S have S-zero divisors?
    vii. Find all ideals in S.

viii. Can S have S-ideals?

42. Let R = {$N_{co}$ ($Z_{18}$), ×, +} be a ring.
   i. Is R a S-ring?
   ii. Find S-ideals if any in R.
   iii. Does R have S-units?
   iv. Find order of R.
   v. If I = {[a, 0) | [a, 0) ∈ $N_{co}$ ($Z_{18}$)} ⊆ R is an ideal of R find R/I.
   vi. Is R/I an integral domain or a semifield?
   vii. What is the order of R/I?
   viii. Find S-idempotents of any in R.

43. Let R = {$N_c$ ($Z_{11}$), +, ×} and S = {$N_c$ ($Z_{12}$), +, ×} be rings.
   i. Find the similarities and dissimilarities between S and R.
   ii. Define a ring homomorphism ϕ from R to S with nontrivial kernel.
   iii. If kernel ϕ is M; is M an ideal of R?
   iv. Can R ≅ S?
   v. Are these rings S-rings?
   vi. Can R and S have S-ideals?
   vii. Can R and S have S-zero divisors?

44. Let R = ($N_o$ (3Z), +, ×) be a ring.
   i. Is R a S-ring?
   ii. Can R have S-ideals?
   iii. Can R have S-zero divisors?
   iv. Is I = {(0, a) | a ∈ 3Z} ⊆ R be an ideal of R?
   v. Find R/I.

45. Obtain some nice properties enjoyed by R = ($N_{oc}$ ($Z_p$), +, ×), p a prime.

46. Does there exist a semifield of natural intervals of order 53?

47. Let S = $\left\{ \begin{pmatrix} [a_1, b_1] & [a_3, b_3] \\ [a_2, b_2] & [a_4, b_4] \end{pmatrix} \middle| [a_i, b_i] \in N_c (Z_5) \right\}$ be a ring.
   i. Is S a S-ring?
   ii. Find subrings in S.

    iii. Does S have S-ideals?
    iv. Does S contain S-subrings which are not ideals of S?
    v. Find zero divisors in S?
    vi. Can S have S-zero divisors?
    vii. What is the order of S?
    viii. Find units in S.

48. Let $V = \left\{ \sum_{i=0}^{3} [a_i, b_i] x^i \mid [a_i, b_i] \in N_c(Z_5) \right\}, x^4 = 1, 0 \leq i \leq 3\}$

    be a ring.
    i. Is V a S-ring?
    ii. Can V have S-subrings?
    iii. Can V have S-ideals?
    iv. Find for $I = \left\{ \sum_{i=0}^{3} [0, a_i] x^i \mid x^4 = 1, 0 \leq i \leq 3; a_i \in Z_5 \right\}$, the ideal V/I.
    v. Find order of V.
    vi. Can V have S-zero divisors?
    vii. Find whether V/I is an integral domain or a semifield or a semiring.

49. Let $R = \left\{ \sum_{i=0}^{\infty} (a_i, b_i] x^i \mid (a_i, b_i] \in N_{oc}(Z) \right\}$ be a ring.

    i. Is R a S-ring?
    ii. Can R have S-subring?
    iii. Find subrings of R.
    iv. Prove R has zero divisors.
    v. Find ideals of R.
    vi. If $S = \left\{ \sum_{i=0}^{\infty} (0, a_i] x^i \mid (0, a_i] \in N_{oc}(2Z) \right\} \subseteq R$; prove S is an ideal of R. Is S a S-ideal of R? Justify your answer.

50. Let S = {3 × 5 open natural interval matrices with entries from $N_o(Z_5)$} be a semigroup under addition.
    i. Find order of S.
    ii. Is S a S-semigroup?

iii. Find subsemigroups and S-subsemigroups if any in S.

51. Let S = {5 × 5 closed interval matrices with entries from $N_c (Z_2)$} be a semigroup.
    i. Find order of S.
    ii. Is S commutative?
    iii. Find ideals in S.
    iv. Can S have S-ideals?
    v. Find S-subsemigroups in S which are not ideals.
    vi. Does S have S-zero divisors?
    vii. Find S-idempotents if any in S.

52. Let T = {3 × 3 interval matrices from $N_{oc} (Z_7)$} be a semigroup.
    i. Find order of T.
    ii. Is T a S-semigroup?
    iii. Is every subsemigroup of S a S-subsemigroup?

53. Let V = {$N_c (Z_{14})$} be a ring.
    i. What is the order of V?
    ii. Is V a S-ring?
    iii. Find subrings of V which are not S-subrings.
    iv. Find ideals in V.
    v. What is the order of V/I = {[0,a]+I} where I is the ideal of the form {[0,a] / a ∈ $Z_{14}$}?
    vi. Is V/I an integral domain? Justify!

54. Obtain some interesting properties enjoyed by interval ring constructed using natural class of intervals.

55. Prove or disprove all classical results cannot be true in S = {$N_c (Z_{43})$}, a ring.

56. Find the special properties enjoyed by R = $N_{oc} (Z_{40})$, the ring of finite order.

57. Let S = $N_0 (Z_{45})$ and R = $N_c (Z_{54})$ be two rings. Find a homomorphism which has a nontrivial kernel.

58. Prove we can find semifields of finite order which are not fields in the ring $S = \{N_{oc}(Z_{17})\}$.

59. Obtain some striking properties enjoyed by the ring $R = \{N_c(Z), +, \times\}$.

60. How does vector space $V = \{N_c(Q)\}$ over the field $Q$ differ from the vector space $Q$ over $Q$?

61. Let $V = \{N_c(Z_7)\}$ be a vector space over the field $Z_7$.
    i. What is the dimension of V over $Z_7$?
    ii. Find a basis of V.
    iii. Does V have subspaces over $Z_7$?
    iv. Find a linear operator on V.
    v. If T = {collection of all linear operators on V}, what is the algebraic property enjoyed by T.
    vi. Is T a vector space over the field $Z_7$?

62. Let $V = \{([a_1, a_2], [a_2, b_2]) \mid [a_i, b_i] \in N_c(Z_{19}) ; 1 \leq i \leq 2\}$ be a vector space over the field $F = Z_{19}$.
    i. What is the dimension of V over F?
    ii. Find the number of elements in V.
    iii. Find subvector spaces of V.
    iv. Find a linear operator T on V with nontrivial kernel.
    v. Find a basis of V.

63. Let V = {all $2 \times 2$ open interval matrices with entries from $N_o(Z_{61})$}, be a vector space over $Z_{61}$.
    i. Is V a linear algebra?
    ii. Find basis of V.
    iii. Is V finite dimensional?
    iv. Find the number of elements in V.
    v. Find subspaces of V.
    vi. Find a linear operator T on V so that $T^{-1}$ exists.

64. Let R = {$2 \times 2$ closed interval matrices with entries from $N_c(Z_3)$}, be a ring.
    i. Is R commutative? (Prove your claim)
    ii. Find the order of R.

   iii. Does R have right ideals?
   iv. Can R have left ideals only?
   v. Does R contain ideals?
   vi. Is R a S-ring?
   vii. Does R contain a S-subring which is not an ideal?

65. Obtain some nice applications of interval groups $N_c(Z_7)$ under addition.

66. Let $G = \{N_o(Z_6)\}$ be a group under addition.
   i. Find order of G.
   ii. Does G contain subgroups?
   iii. Is G simple?
   iv. Verify Lagrange's theorem for G by finding all subgroups of R.
   v. Is G a Smarandache strong group?

67. Let $G = \{[a, b] \mid a, b \in Z_{11} \setminus \{0\}\}$ be a group under multiplication.
   i. Find order of G.
   ii. Find subgroups of G.
   iii. Find automorphism group of G.

68. Prove $G = \{[a, b] \mid a, b \in Z_{11}\}$ is not a group under multiplication.

69. Find some nice applications of rings using natural class of intervals from $N_c(Z)$.

70. Is $N_{oc}(Z)$ a principal ideal domain?

71. Can $N_c(Z)$ be a unique factorization domain? Justify your claim.

72. Can $N_{co}(R)$ be a principal ideal domain?

73. Determine some nice properties enjoyed by $N_{oc}(Q)$.

74. Distinguish the rings Z and $N_o(Z)$.

75. What is the basic difference between the rings $Z_7$ and $N_c(Z_7)$?

76. Show $N_c(Z)$ has zero divisors.

77. Is $R = N_c(Z_2) \times N_o(Z_3) = \{([a, b], (a', b')) \mid a, b, \in Z_2$ and $a'$ $b' \in Z_3\}$ a ring? (direct product of rings $(N_c(Z_2)$ and $N_o(Z_3))$.

78. Let $P = \{([a_1, b_1], [a_2, b_2] [a_3, b_3]) \mid a_i, b_i \in N_c(Z_8); 1 \leq i \leq 3\}$ be a ring.
    i. Find order of P.
    ii. Find ideals of P.
    iii. Find subrings of P.
    iv. Is P a S-ring?

79. Can interval rings constructed using natural class of intervals be a principal ideal domain?

80. Can a ring of order 43 exists using natural class of intervals?

81. Does every ring $S = \{N_c(Z_n)\}$ ($n < \infty$) contain minimal and maximal ideals?

82. Prove or disprove ring using natural class of intervals has ideals!

83. Does there exists vector spaces of finite dimension built using intervals which are not linear algebras?

84. Is $V = \left\{ \begin{bmatrix} (a_1, b_1) \\ (a_2, b_2) \\ (a_3, b_3) \end{bmatrix} \right\}$ $(a_i, b_i) \in N_o(Z_{13}); 1 \leq i \leq 3\}$ a vector space over the field $Z_{13}$ a linear algebra?
    i. Find number of elements in V.
    ii. Is V finite dimensional?
    iii. Find a basis of V.

    iv. Find atleast 2 subspaces of V.

85. Let V = {all 4 × 4 interval matrices with entries from $N_c(Z_7)$} be a vector space over $Z_7$. W = {all 2 × 8 interval matrices with entries from $N_{oc}(Z_7)$} be a vector space over $Z_7$. Let M = {all linear transformation of V to W}. Is M a vector space over $Z_7$? What is the order of M? (That is number of elements in M).

86. Let V = {3 × 3 interval matrices from $N_c(Z_3)$} be a vector space over $Z_3$, the field.
    i. Find a basis of V.
    ii. Is V a linear algebra?
    iii. Find subspace of V.
    iv. What is order of V? (No. of elements in V).
    v. Find a linear operator on V which is one to one.
    vi. Find a linear operator on V which is not invertible.
    vii. If T = {all linear operators on V}. Is T a vector space over $Z_3$?

87. Let V = {$N_c(Z_{27})$} be a ring, I = {[0, a] | a ∈ $Z_{25}$} ⊆ V be an ideal of V.
    i. Find V/I.
    ii. What is the order of V/I?
    iii. What is the order of V?
    iv. Find zero divisors in V.
    v. Can V have S-zero divisors?
    vi. Can V have S-idempotents?
    vii. Can V have subrings which are not ideals?
    viii. Is V a S-ring?
    ix. Can V have S-ideals?

88. Let S = {$N_o(Z_{43})$} be a ring, I = {(a, 0) | a ∈ $Z_{43}$} be an ideal of S. Answer all the questions (i) to (ix) proposed in problem (87).

89. Let G = {$N_o(Z_{41})$} be a group under addition.
    i. Find order of G.
    ii. Can G have subgroup?

    iii. Find the order of (7, 8) and (13, 16), (11, 2) and (14, 3) in G.
    iv. Find automorphism of G.
    v. Is $H = \{(0, a) \mid a \in Z_{41}\}$ a subgroup of G?

90. Determine some special properties enjoyed by the group $(G, +) = G_1 \times G_2 \times G_3$ where $G_1 = N_o(Z_3)$, $G_2 = N_c(Z_5)$ and $G_3 = N_{oc}(Z_{12})$;
    i. What is the order G?
    ii. Find atleast 5 subgroups of G.
    iii. Find the order of $x = \{((0,2), [4,2], (6,10])\}$ and $y = \{((1,2), [3,4], (10,3])\}$ in G.
    iv. Find a nontrivial automorphism on G with nontrivial kernel.

91. Give an example of a ring S using natural class of intervals which has S-zero divisors.

92. Give an example of a ring S using natural class of intervals which has no S-zero divisors.

93. Can the ring $S = \{N_o(Z_{240})\}$ have S-idempotents?

94. What is the order of the ring $R = N_o(Z_{12}) \times N_c(Z_5)$?

95. Let $R = R_1 \times R_2 \times R_3 \times R_4$ where $R_1 = \{N_c(Z_{13})\}$; $R_2 = \{N_o(Z_{10})\}$, $R_3 = \{N_{oc}(Z_7)\}$ and $R_4 = \{N_{co}(Z_6)\}$ are rings.
    i. Find the order of R.
    ii. Define a homomorphism from R to R with nontrivial kernel.
    iii. Can R have subrings which are not ideals?
    iv. Can R have S-ideals?
    v. What is the condition for R to be a S-ring?
    vi. Obtain zero divisors and S-zero divisors if any in R.

96. Find the eigen values, eigen vectors and the characteristic equation of the matrix;

$$M = \begin{pmatrix} (0,8) & (5,1) & (3,1) & (2,4) \\ (-4,1) & (6,2) & (1,5) & (6,9) \\ (8,-4) & (7,0) & (8,3) & (4,0) \\ (5,0) & (0,2) & (3,5) & (0,8) \end{pmatrix},$$ where the entries are from $N_o(R)$.

97. Let $P = \begin{pmatrix} (0,2) & (-7,1) & (3,0) \\ 0 & (1,2) & (4,2) \\ 0 & (0,0) & (3,1) \end{pmatrix}$ be a matrix with entries from $N_o(R)$. Find eigen values and eigen vectors associated with P.

98. Let $M = \begin{pmatrix} [0,8] & [0,2] & 0 & 0 & 0 \\ [6,9] & 0 & [7,1] & [8,1] & 9 \\ [7,0] & [1,4] & 0 & 7 & 0 \\ [1,2] & 0 & [8,2] & [3,1] & -6 \\ [-3,6] & [9,1] & 4 & 1 & 0 \end{pmatrix}$ be a 5 × 5 matrix from $N_c(R)$.
   i. Find eigen values and eigen vectors of M.
   ii. Is M invertible?
   iii. Is M diagonalizable?

99. Prove any of the classical theorems in linear algebra for the linear algebra's built using natural class of intervals.

100. Let $V = \{N_c(Z_7)\}$ be a vector space over $Z_7$.
   i. Find dimension of V.
   ii. Find the number of elements in V.
   iii. Give a basis of V.
   iv. Find subspaces of V.
   v. Find a nontrivial linear operator on V.

101. Let $V = (N_o(Q))$ be a vector space over the field Q.
   i. Give a basis of V.
   ii. What is the dimension of V?

    iii. Can V have proper subspaces?

102. Let $V = \{N_c(Z_3)\} \times \{N_o(Z_3)\} \times \{N_{co}(Z_3)\}$ be a vector space over the field $Z_3$.
    i. Find the number of elements in V.
    ii. Find a basis for V.
    iii. What is dimension of V?
    iv. Give some subspaces of V over $Z_3$.
    v. Define an invertible linear operator on V.

103. Prove $N_o(Z^+ \cup \{0\})$ is a semiring and not a ring.

104. Is $N_o(Q^+ \cup \{0\})$ a semifield? Justify.

105. Can $N_{co}(3Z^+ \cup \{0\})$ be a semifield? Justify.

106. Give some interesting properties about semirings built using intervals.

107. Why $S = N_{co}(5Z^+ \cup \{0\})$ is not a semifield?

108. Is $S = \{N_c(2Z^+ \cup \{0\})\} \times \{N_c(5Z^+ \cup \{0\})\}$ a semiring?

109. Can the semiring $S = N_c(R^+ \cup \{0\})$ have zero divisors? Justify! Is S a semifield?

110. Is $S = \{N_c(Z_{20})\}$ a semiring? Justify your claim.

111. Can $R = \{N_c(Z_{40})\}$ be a semifield?

112. Can we have a semifield using the class of natural intervals?

113. Let $P = \{([a_1, a_2], [a_2, b_2], \ldots, [a_8, b_8]) \mid [a_i, b_i] \in N_c (R^+ \cup \{0\}), 1 \leq i \leq 8\}$ be a semiring.
    i. Find subsemiring of P.
    ii. Is P a semifield?
    iii. Is P a strict semiring?
    iv. Can P have ideals?

114. Let $S = \left\{ \begin{bmatrix} [a_1,b_1] & [a_4,b_4] & [a_7,b_7] \\ [a_2,b_2] & [a_5,b_5] & [a_8,b_8] \\ [a_3,b_3] & [a_6,b_6] & [a_9,b_9] \end{bmatrix} \middle| [a_i, b_i] \in N_c (Z^+ \cup \{0\}) \right\}$ be a interval semiring.
   i. Is S a strict semiring?
   ii. Can S have subsemiring?
   iii. Can S have zero divisors?
   iv. Is S a S-semiring?
   v. Can S be a semifield?
   vi. Is S a commutative semiring?
   vii. Can S have S-subsemirings?

115. Obtain some interesting applications of interval semirings built using the natural class of intervals.

116. Let $V = \left\{ \sum_{i=0}^{\infty} [0, a_i) x^i \middle| [a_i, b_i) \in N_{co} (R^+ \cup \{0\}) \right\}$ be a semiring.
   i. Is V a strict semiring?
   ii. Can V have zero divisors?
   iii. Can V be a semifield?
   iv. In V a S-semiring?
   v. Can V have anti zero divisors?

117. Let V = {5 × 5 interval matrices with intervals from $N_{co}(Z^+ \cup \{0\})$} be a semiring.
   i. Find ideals if any, in V.
   ii. Prove V is noncommutative.
   iii. Can V have zero divisors and antizero divisors?
   iv. Can V have left ideals which are not right ideals?
   v. Is V a S-semiring?
   vi. Find right ideals which are not left ideals in V.

118. Let $V = \left\{ \sum_{i=0}^{\infty} [a_i, b_i) x^i \middle| [a_i, b_i) \in N_{co} (Z_5) \right\}$ be a vector space over the field $Z_5$.

    i. Find a basis of V.
    ii. Find subvector spaces of V
    iii. Is V a linear algebra over V?

119. Let $V = \left\{ \sum_{i=0}^{8} (a_i, b_i) x^i \mid (a_i, b_i) \in N_o (Z_{13}) \right\}$ be a vector space over the field $Z_{13}$.
    i. Find a basis of V.
    ii. Find the number of elements in V.
    iii. Can V be a linear algebra?
    iv. Find subvector spaces of V.

120. Let $V = \left\{ \sum_{i=0}^{3} [a_i, b_i] x^i \mid [a_i, b_i] \in N_c (R) \right\}$ be a vector space over the field R.
    i. Prove V is finite dimensional.
    ii. Find a basis for V.
    iii. Prove V is of infinite order.
    iv. Find a linear operator T on V which is invertible.
    v. Can V be a linear algebra?
    vi. Find atleast 3 vector subspaces of V.
    vii. Find a projection of a vector subspace W of V onto V.

121. Can $V = \left\{ \sum_{i=0}^{7} (a_i, b_i] x^i \mid (a_i, b_i] \in N_{oc} (Z_2) \right\}$ be a vector space over the field $Z_2$?
    i. Find the number of elements in V.
    ii. Find a basis of V over $Z_2$.
    iii. What is the dimension of V over $Z_2$?
    iv. Can V have subvector spaces?

122. Let $W = \left\{ \sum_{i=0}^{17} (a_i, b_i) x^i \mid (a_i, b_i) \in N_o (Z_{17}) \right\}$ be a vector space over the field $Z_{17}$. V = {all 3 × 2 interval matrices with intervals $(a_i, b_i) \in N_o(Z_{17})$} be a vector space over the field $Z_{17}$.

i. Find subspaces of V and W.
ii. Define a linear transformation T from V to W so that T is invertible.
iii. Define a linear transformation T from W to V so that T is not invertible.
iv. Find a basis of V over $Z_{17}$.
v. Find a basis of W over $Z_{17}$.
vi. Can V as a vector spaces be embedded in the vector space W?

123. Let V = {all 5 × 5 interval matrices with entries from $N_c(Z_3)$} be a linear algebra over the field $Z_3$.
   i. Find the number of elements in V.
   ii. Find a basis for V over $Z_3$
   iii. Find sublinear algebras of V.
   iv. Find an invertible linear operator on V.
   v. Can V be written as a direct sum of subspaces?
   vi. Find a linear operator on V which is non invertible.

124. Suppose V = {$N_c(Z_{12})$} be a Smarandache vector space of type II over the S-ring $Z_{12}$.
   i. Find a basis of V.
   ii. Find the number of elements in V.
   iii. Find subspaces if any in V.

125. Let W = {$N_c(Z_{15})$} be a Smarandache vector space of type II over the S-ring $Z_{15}$.
   i. Find an invertible linear operator on W.
   ii. What is the dimension of W over $Z_{15}$?

126. Give an example of a module over a ring R which is not a Smarandache vector space over the ring R.

127. Obtain some interesting properties enjoyed by the Smarandache vector spaces of type II.

128. Is V = ($N_o(Q) \times N_o(Q)$) a Smarandache vector space of type II over the S-ring Q × Q?

129. Can V = ($N_c(Z_{12}) \times N_c(Z_{12}) \times N_c(Z_{12})$) be a Smarandache vector space of type II over the S-ring $N_c(Z_{12})$?

130. Can V = ($N_c(Z_{19})$) be a Smarandache vector space of type II over the S-ring F = $N_c(Z_{19})$?
   i. If so what is the dimension of V over F.
   ii. Is V a S-linear algebra?
   iii. Find a basis for V over F.

131. Let V = {$3 \times 5$ interval matrices with intervals of the form [$a_i, b_i$] $\in N_c(Z_{21})$} be a S-vector space of type II over the S-ring $N_c(Z_{21})$.
   i. Find a basis of V.
   ii. Find the number of elements in V.
   iii. Find atleast 3 subspaces of V over $N_c(Z_{21})$.
   iv. Find an invertible linear operator on V.
   v. If $V_1$ is a vector space over the S-ring $Z_{21}$, find the similarities and differences between V and $V_1$.

132. Let $\left\{ \sum_{i=0}^{8} (a_i, b_i] x^i \mid (a_i, b_i] \in N_{oc}(Z_{15})\right\}; 0 \le i \le 8\}$ be a S-vector space of type II over the S-ring $N_{oc}(Z_{15})$.
   i. Find a basis of V over $N_{oc}(Z_{15})$.
   ii. Find the number of elements in V.
   iii. Can V have subspaces?
   iv. If V is made into a S-vector space of type II over $Z_{15}$ study questions (i), (ii) and (iii).
   v. Is V a vector space over the field F = {0, 5, 10} $\subseteq Z_{15}$?

133. Let V = {all $10 \times 5$ interval matrices with entries from $N_c(Z_7)$} be a S-vector space of type II over the S-ring $N_c(Z_7)$ and W = {all $5 \times 5$ interval matrices with entries from $N_c(Z_7)$} be a S-vector space of type II over $N_c(Z_7)$ the S-ring.
   i. Find a basis of V and W over $N_c(Z_7)$.
   ii. Find the number of elements in V.
   iii. What is the dimension of V over $N_c(Z_7)$?
   iv. Find a basis of W over $N_c(Z_7)$.

- v. Find a linear transformation from V to W.
- vi. If T = {set of all linear transformations from V to W}; Is T a vector space over $N_c(Z_7)$?
- vii. Suppose V is a vector space over $Z_7$ and W a vector space over $Z_7$ find a linear transformation from V to W.

134. Let V = {All 3 × 3 interval matrices with entries from $N_c(Z_{49})$}. Is V a S-vector space of type II over $N_c(Z_{49})$? Justify your claim.

135. Determine any nice property enjoyed by these special types of S-vector spaces.

136. Suppose V be a linear algebra of interval matrices over a field F and V be treated only as a vector space over F. Is there a difference between the number of base elements in general. Justify your claim by examples.

137. Let $V = \left\{ \begin{pmatrix} [a_1,b_1] & [a_2,b_2] \\ [a_3,b_3] & [a_4,b_4] \end{pmatrix} \middle| [a_i, b_i] \in N_c(Z_2); 1 \leq i \leq 4 \right\}$ be a S-vector space of type II over $N_c(Z_2)$.
   - i. Find a basis of V over $N_c(Z_2)$.
   - ii. What is dimension of V over $N_c(Z_2)$?
   - iii. Find subspaces of V.
   - iv. Is V a S-linear algebra?

138. Suppose S = {(a, b) | 0 ≤ a, b ≤ 1} = {(a, b) | (a, b) ∈ $N_o$ ([0,1])} be natural fuzzy a semigroup under multiplication.
   - i. Find ideals of S.
   - ii. Can S have subsemigroups which are not ideals?
   - iii. Can S have zero divisors?
   - iv. Can S have idempotents?
   - v. Can S be a S-semigroup?
   - vi. Can S hae S-subsemigroup?
   - vii. Prove S is of infinite order.

139. Let S = {[a, b] | [a, b] ∈ $N_c$ ([0,1]} be a fuzzy semigroup under min operation.

    i. Find fuzzy subsemigroups in S.
    ii. Can S have S-ideals?
    iii. Can we have the concept of zero divisors in S?
    iv. Is S a S-semigroup?
    v. Obtain any interesting property enjoyed by S.

140. Let $S = \{([a_1, b_1], [a_2, b_2], [a_3, b_3]) \mid [a_i, b_i] \in N_c([0,1])\}$ be a fuzzy semigroup with multiplication as the operation on it.
    i. Find zero divisors in S.
    ii. Can S have S-zero divisors?
    iii. Is S a S-semigroup?
    iv. What is order of S?
    v. Can S have fuzzy subsemigroups which are not ideals?
    vi. Can S have S-idempotents?
    vii. Is S a free semigroup?

141. Let $V = \left\{ \sum_{i=0}^{\infty} (a_i, b_i) x^i \,\middle|\, (a_i, b_i) \in N_o(Q) \right\}$ be a linear algebra over Q.
    i. Find a basis of V over Q?
    ii. Find subspaces of V over Q.
    iii. Find a linear operator on V which is noninvertible.
    iv. If T = {collection of all linear operators on V}, is T a linear algebra over Q?

142. Suppose V in problem (141) is a S-vector space of type II over $N_o(Q)$ study problems (i) to (iv) in problem (141).

143. Let M = {all 3 × 2 interval matrices with entries from $N_c(Z_5)$}, be a vector space over $Z_5$.
    i. Find a basis of M.
    ii. Find the order of M.
    iii. Find subspaces of M.
    iv. Find a non invertible linear operator of M.
    v. Is M better than V = {all 3 × 2 matrices with entries from $Z_5$} over $Z_5$?
    vi. Is $V \subseteq M$?

144. Let V = {N_c (Z_{42})}, a S-vector space of type II over $Z_{42}$.
    i. What is dimension of V over $Z_{42}$?
    ii. Can V be a S-linear algebra over $Z_{42}$?
    iii. Can V have subspaces?
    iv. Give a nontrivial linear operator on V with nontrivial kernel.
    v. Can V be written as a direct sum of subspaces?

145. Let $M = \begin{bmatrix} 8 & (3,0) & (9,1) & 0 \\ (1,2) & 0 & (4,1) & 2 \\ 0 & (5,1) & 0 & (1,2) \\ (3,9) & 0 & 7 & 0 \end{bmatrix}$ be a interval matrix with entries from $N_o(R)$?
    i. Find eigen interval values of M.
    ii. Find the associated interval vectors of M.
    iii. Can M be diagonalizable?
    iv. Is M invertible?
    v. Find det M.

146. Suppose V = {N_c (Z_{19})} be a vector space over $Z_{19}$. Can on V be defined a normal operator?

147. Can Spectral theorem be extended for linear algebras built using intervals defined over a field Q?

148. Can the notion of inner product space be extended to interval linear algebras?

149. Prove for every finite interval linear algebra (vector space) we have a sublinear algebra (subvector space) which satisfies the Generalized Cayley Hamilton Theorem.

150. On similar lines prove every interval linear algebra has a sublinear algebra of finite dimension which satisfies cyclic decomposition theorem.

151. Let V = {$N_c(R)$} be a vector space over the field R of reals. Clearly W = R ⊆ $N_c(R)$ so R is a vector space of dimension one over R.
    i. Find a basis of V.
    ii. Is V finite dimensional?
    iii. Find subspaces of V.
    iv. Find linear operator T on V which preserves the subspace W.

152. Let V = ($N_c(Q) \times N_c(Q) \times N_c(Q) \times N_c(Q)$) be a vector space over Q.
    i. Does W = (Q×Q×Q×Q) ⊆ V satisfy restricted primary decomposition theorem?
    ii. Does W satisfy cyclic decomposition theorem?
    iii. Does W satisfy Generalized Cayley Hamilton Theorem?
    iv. Can on W a diagonalizable operator be determined?

153. Obtain some interesting results on fuzzy rings defined using $N_c([0,1])$.

154. Can a fuzzy group be constructed using $N_o([0,1])$?

155. Let V = {($[a_1, b_1], [a_2, b_2], …, [a_8, b_8]$) where $[a_i, b_i] \in N_{co}([0,1])$; $1 \le i \le 8$}.
    i. Can V be a fuzzy semigroup under usual product?
    ii. Can V be a fuzzy semigroup under max operation?
    iii. Can V be a fuzzy group under max operation?

156. Define fuzzy semigroup using $N_c$ (0,1).

157. Let R = ($N_c(Z)$, +, ×) be a ring. Define a map η : R → [0,1] by η (0 0) = (1, 1), η(a,b) = (1/a, 1/b) if a ≠ 0, b ≠ 0. Is (R, η) a fuzzy ring with suitable operations?

158. Obtain some interesting properties about neutrosophic natural interval semigroups $N_o(\langle Z \cup I \rangle)$.

159. Let W = {$N_{oc}(Z_5 I)$} be a pure neutrosophic semigroup under multiplication;

- i. Find order of W.
- ii. Find ideals if any in W.
- iii. Can W have subsemigroups which are not ideals?
- iv. Find zero divisors in W.
- v. Can W have S-zero divisors?
- vi. Determine units of W.
- vii. Are these units S-units?
- viii. Prove W cannot have idempotents or S-idempotents.

160. Let $S = \{N_c(Z_8I)\}$ be a pure neutrosophic semigroup under multiplication modulo 8. Study the questions (i) to (viii) given in problem (159).

161. Let $M = \{N_c(ZI)\}$ be a group under addition. Find some interesting related properties enjoyed by M.

162. Let $P = \{N_o(RI)\}$ be a group under addition. Compare this group with M in problem 161 and determine the difference between them.

163. Let $T = \{N_c(<Z_9 \cup I>)\}$ be a semigroup under multiplication, elements of $N_c$ are of the form $[a + bI, c + dI]$ where a, b, c, d $\in Z_9$ with $I^2 = I$.
  - i. Find the order of T.
  - ii. Find subsemigroups of T.
  - iii. Can T have subsemigroup which are not ideals?
  - iv. Can T have S-ideals?
  - v. Prove T has zero divisors.
  - vi. Can T have S-idempotents?
  - vii. Can T have S-semigroups?
  - viii. Is T a S-semigroups?
  - ix. Can T have S-units?
  - x. Can T have S-antizero divisors?
  - xi. Prove order of every element is of finite order.

164. Let $S = \{N_o(<Z_{12} \cup I>)\}$ be semigroup under multiplication.
  - i. What is the order of S?
  - ii. Find ideals in S.
  - iii. Is every ideal in S a S-ideal?

iv. Does S have subsemigroups which are not S-ideals?
v. Can S have S-zero divisors?
vi. Find S-units in S.
vii. What is the neutrosophic order of $(2 + 5I, I + 1)$ in S?
viii. Does the order or neutrosophic order of every element divide the order of S?
ix. Does the order of ideals of the semigroup S divide the order of S?

165. Let $V = \{N_c(Z_6I)\}$ be a semigroup under multiplication.
   i. What is the order of V?
   ii. Find subsemigroups in V.
   iii. Is V a S-semigroup?
   iv. Does V have S-ideal?
   v. Can V have S-zero divisors?
   vi. Can V have idempotents which are not S-idempotents?

166. Let $M = \{N_c(\langle Z_6 \cup I \rangle)\}$ be a semigroup (a) study questions (i) to (vi) mentioned in problem (165). (b) Compare V and M.

167. Let $R = \{N_c(Z), +, \times\}$ be a ring. Define a map $\eta : R \to [0,1]$ so that $(R, \eta)$ is a fuzzy interval ring.

168. Let $P = \{N_o(Q), +, \times\}$ be a ring. Find $\eta: P \to [0,1]$ so that $(P, \eta)$ is a fuzzy ring. How many fuzzy rings can be constructed using P?

169. Let $R = \{$all $5 \times 5$ interval matrices with entries from $N_c(Z)\}$ be a ring under matrix addition and multiplication. Define $\eta : R \to [0,1]$ so that $(R, \eta)$ is a fuzzy ring.
   i) Find fuzzy subrings of R. Can R have fuzzy ideals?

170. Let $S = \{N_c (Z^+ \cup \{0\}), +, \times\}$ be a semiring. Define $\eta : S \to [0,1]$ so that $(S, \eta)$ is a fuzzy semiring.

171. Let $S = \{N_{co}(R^+ \cup \{0\}), +, \times\}$ be a semiring. Define $\eta : S \to [0,1]$ to make $(S, \eta)$ is a fuzzy semiring.

i) How many fuzzy semirings can be constructed using S? Does every map $\eta : S \to [0,1]$ make $(S, \eta)$ a fuzzy semiring?

172. Let $G = \{N_c(Z_{10})\}$ be a group under addition. Define $\eta : G \to [0,1]$ so that $(G, \eta)$ is a fuzzy group.

173. Let $G = \{N_c(R^+)\}$ be a group under multiplication. $\eta : G \to [0,1]$ be a map such that $(G, \eta)$ is a fuzzy group.

174. Let $W = \{N_{oc}(Q^+)\}$ be a group under multiplication. $\eta : W \to [0,1]$ be a map such that $(W, \eta)$ is a fuzzy group.

175. Let $T = \{N_{co}(Z)\}$ be a group under addition. Define $\eta : T \to [0,1]$, be a map, then $(T, \eta)$ is a fuzzy group.

176. Let $T = \{N_{co}(Z_8)\}$ be a group under addition. Define $\eta: T \to [0,1]$ so that $(T, \eta)$ is a fuzzy group.

177. Let $G = \{N_c(Z_{11} \setminus \{0\})\}$ be a group under multiplication. Define $\eta : G \to [0,1]$ so that $(G, \eta)$ is a fuzzy group.

178. Let $W = \{N_o(Z_{13} \setminus \{0\})\}$ be a group under multiplication. Define $\eta : W \to [0,1]$ so that $(W, \eta)$ is a fuzzy group.

179. Let $M = \{N_c(Z_{29} \setminus \{0\})\}$ be a group under multiplication. Define $\eta : M \to [0,1]$ so that $(M, \eta)$ is a fuzzy group.

180. Let $G = \{N_c(Z_5)\}$ be a group under addition. Define $\eta$ so that $(G, \eta)$ is a fuzzy group.

181. Let $G = N_c(\langle Z_5 \cup I \rangle)$ be a neutrosophic group under addition.
   i. Find the order of G.
   ii. Find subgroup of G.

182. Let $G = N_c(ZI)$ be a semigroup under multiplication.
   i. Find ideals of G.

    ii. Is G a S-semigroup?
    iii. Does S contain S-ideals?

183. Obtain some interesting properties enjoyed by pure neutrosophic algebraic structures.

184. Obtain some interesting properties associated with $N_c(ZI)$.

185. Prove using $N_c(\langle R^+ \cup I \rangle \cup \{0\})$, S-semirings can be built. Determine some striking properties enjoyed by these S-semirings.

186. Let V = {all 3 × 3 interval matrices with entries from $N_c(Z_5I)$}.
    i. Prove V is a group under addition modulo 5.
    ii. Prove V under multiplication is only a semigroup.
    iii. Prove V is a non commutative ring.
    iv. What is order of V?
    v. Can V as a ring have S-ideals?
    vi. Can V have subrings which are not ideals?

187. Let V = {all 3 × 2 neutrosophic interval matrices with entries from $N_c(\langle Q \cup I \rangle)$} be a vector space of neutrosophic interval matrices over the field Q.
    i. Find a basis for V.
    ii. Is V finite dimensional?
    iii. Find a linear operator on V which is invertible.
    iv. Can V have subspaces?
    v. Can $V = \cup W_i$, ($W_i$'s are subspaces of V)?

188. Let V = {all 2 × 2 neutrosophic interval matrices with entries from $N_o(\langle Z_{19} \cup I \rangle)$} be a vector space of neutrosophic interval matrices over $Z_{19}$ = F, the field of characteristic 19.
    i. Find a basis of V.
    ii. Find the number of elements in V.
    iii. Prove V is a linear algebra of neutrosophic interval matrices over $Z_{19}$.

iv. Find T = {set of all linear operators from V to V}. Is T a linear algebra over F = $Z_{19}$?
v. Does V have subspaces?
vi. Can V be written as a direct sum of subvector spaces?

189. Let M = {set of all 2 × 7 interval matrices with entries from $N_{oc}$ (<$Z^+\cup I$> ∪ {0})} be a semivector space of neutrosophic interval matrices over the semifield F = $Z^+ \cup \{0\}$ of type I.
    i. Find a basis of M.
    ii. Is M finite dimensional?
    iii. Prove M is not a semilinear algebra.
    iv. Find subvector spaces of M.
    v. What is the dimension of M over F?

190. If in the problem (188) the field $Z_{19}$ is replaced by $Z_{19}I$; study the problem (i) to (vi) what are the differences by replacing $Z_{19}$ by $Z_{19}I$?

191. Describe some interesting properties enjoyed by the class of fuzzy neutrosophic open intervals.

192. Let V = {all 2 × 2 interval matrices built using $N_c$(<[0, I] ∪ [0, 1]>)} be a collection of fuzzy neutrosophic closed intervals.
    i. Define algebraic operations on V so that V is a semigroup.
    ii. Can V be made into a semivector space?
    iii. What is the richest algebraic structure that can be constructed using V?

193. Let W = $\left\{ \begin{bmatrix} a_1 \\ a_2 \\ a_3 \\ a_4 \\ a_5 \\ a_6 \end{bmatrix} \middle| a_i \in N_c(Z_6); i = 1, 2, 3, 4, 5, 6 \right\}$ be a group under addition.

    i. Find order of W.

    ii. Find subgroups of W.
   iii. Find the quotient group using any subgroup of your choice.
   iv. Find the automorphism group of W.
    v. Define a homomorphism from W to W so that kernel is nontrivial.

194. Let M = {set of all 5 × 5 neutrosophic interval matrices with entries from $N_o(Z_3I)$} be a semigroup.
    i. Find order of M.
    ii. Is M commutative?
   iii. Is M a S-semigroup?
   iv. Can M have S-ideals?
    v. Does M have subsemigroups which are not ideals?
   vi. Can M have S-zero divisors?
   vii. Can M have idempotents?
  viii. Is every zero divisors in M a S-zero divisors?

195. Let P = {set of all 4 × 4 neutrosophic interval matrices with entries from $N_{oc}(\langle Z_4 \cup I \rangle)$} be a ring.
    i. Find the order of P.
    ii. Is P commutative?
   iii. Find ideals in P.
   iv. Is every subring of P an ideal of P? Justify.
    v. Is P a S-ring?
   vi. Find zero divisors and S-zero divisors of P.
   vii. Does P have idempotents which are not S-idempotents?
  viii. Find a homomorphism on P with a non-trivial kernel.

196. Let M = $\left\{ \sum_{i=0}^{\infty} (a_i, b_i] x^i \mid (a_i, b_i] \in N_{oc}(\langle Z_6 \cup I \rangle) \right\}$ be a ring.
    i. Find ideals in M.
    ii. Is every ideal on M a principal ideal?
   iii. Can M have minimal ideals?
   iv. Can M have subrings which are not ideals?
    v. Is M a S-ring?
   vi. Can M have zero divisors?

197. Derive some interesting properties enjoyed by neutrosophic polynomial ring with coefficients from $N_o (\langle Z \cup I \rangle)$.

198. Find maximal ideals in $N_{oc}(\langle Q \cup I \rangle)$.

199. Can $N_{oc}(RI)$ have ideals?

200. Can $N_c(\langle R \cup I \rangle)$ have ideals?

201. Let $M = \{N_o(\langle R^+ \cup I \rangle) \cup \{0\}\}$ be a semiring. Is M a S-semiring? Prove M is not a semifield.

202. Let $R = N_{oc}(\langle Q \cup I \rangle)$ be a ring. $I = \{(0,a] \mid a \in \langle Q \cup I \rangle\} \subseteq R$.
   i. Is I an ideal of R.
   ii. Find R/I.
   iii. Can R have ideal?
   iv. Is R a S-ring?

# FURTHER READING

# INDEX

**S**

# ABOUT THE AUTHORS

**Dr.W.B.Vasantha Kandasamy** is an Associate Professor in the Department of Mathematics, Indian Institute of Technology Madras, Chennai. In the past decade she has guided 13 Ph.D. scholars in the different fields of non-associative algebras, algebraic coding theory, transportation theory, fuzzy groups, and applications of fuzzy theory of the problems faced in chemical industries and cement industries. She has to her credit 646 research papers. She has guided over 68 M.Sc. and M.Tech. projects. She has worked in collaboration projects with the Indian Space Research Organization and with the Tamil Nadu State AIDS Control Society. She is presently working on a research project funded by the Board of Research in Nuclear Sciences, Government of India. This is her $54^{th}$ book.

On India's 60th Independence Day, Dr.Vasantha was conferred the Kalpana Chawla Award for Courage and Daring Enterprise by the State Government of Tamil Nadu in recognition of her sustained fight for social justice in the Indian Institute of Technology (IIT) Madras and for her contribution to mathematics. The award, instituted in the memory of Indian-American astronaut Kalpana Chawla who died aboard Space Shuttle Columbia, carried a cash prize of five lakh rupees (the highest prize-money for any Indian award) and a gold medal.
She can be contacted at vasanthakandasamy@gmail.com
Web Site: http://mat.iitm.ac.in/home/wbv/public_html/
or http://www.vasantha.in

**Dr. Florentin Smarandache** is a Professor of Mathematics at the University of New Mexico in USA. He published over 75 books and 150 articles and notes in mathematics, physics, philosophy, psychology, rebus, literature.

In mathematics his research is in number theory, non-Euclidean geometry, synthetic geometry, algebraic structures, statistics, neutrosophic logic and set (generalizations of fuzzy logic and set respectively), neutrosophic probability (generalization of classical and imprecise probability).  Also, small contributions to nuclear and particle physics, information fusion, neutrosophy (a generalization of dialectics), law of sensations and stimuli, etc. He can be contacted at smarand@unm.edu